\documentclass[12pt, ShowComments]{template}

\usepackage{filecontents}

\def\crossref#1{\cite[\cref{partI:#1}]{hennionholsteinrobaloI}}

\DeclareRobustCommand{\subtitle}[2]{\vspace{.5\baselineskip}\\#1: #2}

\newsavebox{\subdiagram}

\begin{document}
\todo{Commts.\\ON}

\setcounter{page}{1}
\normalem

\date{}
\title{Gluing invariants of Donaldson--Thomas type\subtitle{Part II}{Matrix Factorizations}}

\authorbenjamin
\authorjulian
\authormarco

\begin{abstract}This paper is a follow-up to \cite{hennionholsteinrobaloI}. Let $X$ be a a $(-1)$-shifted symplectic derived Deligne--Mumford stack. Thanks to the Darboux lemma of Brav--Bussi--Joyce, $X$ is locally modeled by derived critical loci of a function $f$ on a smooth scheme $U$.
In this paper we study the gluing of the locally defined $2$-periodic (big) dg-categories of matrix factorizations $\MF^\infty(U,f)$.
We show that these come canonically equipped with a structure of a $2$-periodic crystal of categories (\ie an action of the dg-category of $2$-periodic $D$-modules on $X$) compatible with a relative Thom--Sebastiani theorem expressing the equivariance under the  action of quadratic bundles. 

  As our main theorem we show that the locally defined categories $\MF^\infty(U,f)$ can be glued along $X$ as a sheaf of crystals of 2-periodic dg-categories ``up to isotopy'', under the prescription of
  orientation data controlled by three obstruction classes.
  This result generalizes the gluing of the Joyce's perverse sheaf of vanishing cycles and partially answers conjectures by Kontsevich--Soibelman and Toda in motivic Donaldson--Thomas theory.
\end{abstract}

\maketitle
\renewcommand{\subtitle}[2]{ -- #1}

\tableofcontents

\section{Introduction}

This paper is a follow-up to \cite{hennionholsteinrobaloI} where we introduced the moduli of Darboux coordinates of a $(-1)$\nobreakdashes-shifted symplectic derived Deligne--Mumford stack $X$ which parametrizes local presentations of $X$ as the derived critical locus of a function $f$ on a smooth formal scheme $\formalU$.
Singularity invariants of such pairs $(\formalU, f)$ -- such as the Milnor number $\mu_f$ and the perverse sheaf of vanishing cycles $\PJoyce_{\formalU,f}$ -- are then defined for any choice of Darboux coordinates without ambiguity. These invariants satisfy two fundamental properties:
\begin{enumerate}
  \item A form of isotopy invariance: For any family of automorphisms $(\phi_t)_{t \in \basefield}$ of a pair $(\formalU, f \colon \formalU \to \affineline{})$ fixing the critical locus, the induced family of automorphisms of the local invariant is constant in $t$, \ie $\forall t\in \basefield ~ \phi_t^\ast=\phi_0^*$. 
  Milnor numbers trivially satisfy this (they have no automorphisms).
  In the case of vanishing cycles, this amounts to the fact the pullback functor of perverse sheaves $\Perv(X) \to \Perv(X \times \affineline{\basefield})$ is fully faithful.
  \item A quadratic stabilization formula provided by a Knörrer periodicity theorem.
\end{enumerate}
The main result from \cite{hennionholsteinrobaloI} shows that those fundamental properties are sufficient for gluing local singularity invariants into a global invariant of $X$ (up to a twist), thus recovering results of Behrend \cite{MR2600874} (\crossref{examplemilnornatural}) and  Brav--Bussi--Dupont--Joyce--Szendroi \cite{MR3353002} (\crossref{naturaltransformationperversesection}).

To resolve the twist, one then needs additional orientation data, which in the case of perverse sheaves amounts to a square root of the canonical bundle of $X$.

\subsection{In this paper.} We are interested in a conjecture of Kontsevich and Soibelman \cite{MR2851153,MR2681792} (see also Toda's \cite{Toda2023}).
In essence, we study the gluing of the locally defined $2$-periodic dg-categories of matrix factorizations $\MF^\infty(\formalU,f)$ (see \cite{MR570778,douglas2001d,MR2101296,MR2437083,MR2824483,MR3121870}).

 
Unfortunately, contrary to the case of Milnor numbers and vanishing cycles, we could not prove the required form of isotopy invariance for matrix factorizations (see  \cref{naiveA1invarianceofMF}  below).

Nevertheless, the locally defined categories $\MF^\infty(\formalU,f)$ can be glued along $X$ as a sheaf of dualizable crystals of 2-periodic dg-categories and dg-functors ``up to isotopy'' (\ie we enforce $\affineline{}$-isotopy invariance of dg-functors).

\begin{theoremnonumber}[{see \cref{theoremgluingMF}}]\label{theoremB}
  Let $X$ be a $(-1)$\nobreakdashes-shifted symplectic derived Deligne--Mumford stack.
  Then the locally defined $2$\nobreakdashes-periodic dg-categories of matrix factorizations $\MF^\infty(\formalU,f)$ -- associated to local Darboux models of $X$ -- can be glued up to isotopies of dg-functors, provided the following orientation data is given:
  \begin{itemize}
    \item a continuous function $X \to \quot{\integers}{2}$ controlling the parity of the dimension of the local Darboux models -- corresponding to a trivialization of the zero class $\beta_1 = 0 \in \cohomology^1(X, \quot{\integers}{2})$;
    \item a square root of the canonical bundle $\canonicalbundle_X \coloneqq \det(\cotangent_X)$ -- corresponding to a trivialization of a class $\beta_2 \in \cohomology^2(X, \quot{\integers}{2})$;
    \item a spinorial theory\footnote{See \cref{defcategoricalorientation}} on $X$ (akin to a shifted Pin structure) -- corresponding to a trivialization of a class $\beta_3 \in \cohomology^3(X, \quot{\integers}{2})$.
  \end{itemize}
  The space of such choices forms a gerbe over $X$, over which such a gluing always canonically exists.
\end{theoremnonumber}

This result generalizes the gluing of the Joyce's perverse sheaf of vanishing cycles and partially answers the conjectures by Kontsevich--Soibelman \cite{MR2851153,MR2681792} and Toda \cite{Toda2023} in motivic Donaldson--Thomas theory.
Note that the obtained invariant is easily seen to be non-trivial in general: For instance it allows one to recover Behrend's function (and hence the numerical Donaldson--Thomas invariants) -- see \cref{invariantnontrivial}.

In order to prove this theorem, we will extend the definition of matrix factorizations and the proof of their important properties (base change, Thom--Sebastiani equivalences and Knörrer periodicity) to a very general context.
This will allow us to characterize the aforementioned orientation data in terms of matrix factorizations of quadratic forms and their relation to Clifford algebras.
The three classes $\beta_1$, $\beta_2$ and $\beta_3$ will then arise naturally, corresponding to (a shift of) the Stiefel--Whitney classes governing trivializations of Clifford algebras (up to Morita equivalence).

\begin{remark}
  In the first arxiv version of \cite{hennionholsteinrobaloI}, \cref{theoremB} was formulated as a gluing in \emph{noncommutative motives}.
  The current formulation of \cref{theoremB} is strictly stronger than any form of motivic gluing.
  However, some care is needed as the theorem
  does not belong to the world of compactly generated categories (even though large categories of matrix factorizations are compactly generated) since the constructions involved do not preserve compact objects.
  Thus our result does not directly fit the definition of noncommutative motives of \cite{MR3281141,MR2822869, tabuada-notes, MR3070515}.  
  One could instead use the recent theory of localizing invariants for dualizable dg-categories \cite{efimov2024k,efimov2025localizinginvariantsinverselimits} to deduce motivic results.
  We will not pursue this in the current paper.
\end{remark}

It remains an open question whether or not matrix factorizations satisfy isotopy invariance (so that we could glue on the nose and not just ``up to isotopies'').
Let us state this question explicitly:
\begin{openquestion}\label{naiveA1invarianceofMF}
  Let $\formalU$ is a smooth formal scheme equipped with a function $f \colon \formalU \to \affineline{}$. Assume that $0$ is the only critical value, that $\Crit(f)$ is a scheme and that $\Crit(f)_\red \simeq \formalU_\red$.
  Consider $1$-parameter families $(\phi_t)_{t \in \basefield}$ of automorphisms of $\formalU$ fixing the function $f$ and such that for every $t$, the automorphism $\phi_t$ induces the identity on the derived critical locus $\dCrit(f)$ (with its $(-1)$-symplectic structure).
  
  Is the induced functor $\phi_t^! \colon \MF(\formalU, f) \to \MF(\formalU, f)$ independent of $t$ for every such family of automorphisms $(\phi_t)_{t \in \basefield}$?
\end{openquestion}
A positive answer would, after some $\infty$-categorical work, imply the gluing of dg-categories of matrix factorization, without needing to force isotopy invariance.

From derived deformation theory, one can prove the question has a positive answer for formal isotopies (\ie $1$-parameter families indexed by $\fAffLine$).
As pointed to us by Bertrand Toën, Artin's algebraization then implies isotopy invariance of matrix factorizations in the case of isolated singularities (such as \crossref{ambiguitymorphisms}).

\subsection{Acknowledgments}\!\!\!\!\footnote{\tiny The first and third authors were supported by the grant ANR-17-CE40-0014, the second author by the Deutsche Forschungsgemeinschaft (DFG, German Research Foundation) through EXC 2121 ``Quantum Universe'' -- project number 390833306 -- and SFB 1624 ``Higher structures, moduli spaces and integrability'' -- project number 506632645.
}
We thank Dario Beraldo, Frédéric Déglise, Jean Fasel, Niels Feld, Tudor Pădurariu, Tony Pantev, Hyeonjun Park, Massimo Pippi, Pavel Safronov and Bertrand Toën for several helpful discussions.

\subsection{Outline}
To aid navigation, we outline here the contents of each section. 

After setting notations in \cref{notations}, \cref{sectionsheavesofcats} reviews and extends some of the foundational aspects of categories of Ind-coherent sheaves on derived inf-schemes developed in \cite{GR-I, GR-II}, setting the stage for \cref{MFsection}, where, in order to simultaneously accommodate smooth formal schemes (our LG-pairs) and de Rham stacks (responsible for flat connections and categorical crystal structures), we establish the framework of matrix factorization categories on derived inf-schemes. Moreover, we refine Preygel's Thom-Sebastiani and Knörrer periodicity results from \cite{MR3121870} in this context. All contents in this section  are independent from our previous work \cite{hennionholsteinrobaloI}.

In \cref{formalsectionMFfromdarboux} we start combining the results of this paper with the results from  \cite{hennionholsteinrobaloI}.
In \cref{axiomatizationofMFinvariants} we assemble several forms of functorial behavior of matrix factorizations and their symmetric monoidal compatibilities established in \cref{MFsection}. In \cref{matrixfactorizationasasectionofstack} we exhibit $\MF^\infty$ as a morphism of stacks from the Darboux stack -- $\Darbstack_X$ -- of a $(-1)$-shifted symplectic Deligne-Mumford stack $X$ to the stack of sheaves of 2-periodic dg-categories on $X_\derham$ which are compatible with the action of quadratic bundles with a flat connections via the Thom-Sebastiani theorem.  Altogether, \cref{formalsectionMFfromdarboux} is concerned with formally constructing a functor that is clear in a naive sense and may be skipped on a first reading.

In \cref{sectiongluing} we finally give the proof of \cref{theoremgluingMF}. 
In \cref{sectionA1invarianceofMFrevisited} we construct the $\affineline{}$-isotopic version of the functor of matrix factorizations, factoring through the quotient of $\DA$ by the action of  $\QA$. In \cref{sectionhomotopygroupscliff} we explain why there are only three obstruction classes  and finally in  \cref{sectioncategoricalorientationdataandpin} we describe them and state and prove the main theorem.

\subsection{Conventions and notations}
\label{notations}
We adopt the notations in \cite{hennionholsteinrobaloI} but for the reader's convenience we briefly revise the most important in this section.

Throughout the paper we work exclusively over $\basefield$.

\begin{notations}{Unless mentioned otherwise, we use the notations of \cite{lurie-htt, lurie-ha} for higher categories. In this work the difference between big and small $\infty$-categories plays an important role. In order not to overload the notations, we prefer to use a single notation for both sizes. In doing so, we make the effort to specify the context in the text. Therefore, we denote:} \label{notationhighercategories}
  \item by $\inftygpd$ the $(\infty,1)$\nobreakdashes-category of small (resp. big) $\infty$\nobreakdashes-groupoids;
  \item by $\inftycats$ for the $(\infty,1)$\nobreakdashes-category of small (resp. big) $(\infty,1)$\nobreakdashes-categories, just called \icategories from now on;
  \item \label{notationsymmoncats} by $\monoidalcats$ the $(\infty,1)$-category of symmetric monoidal \icategories;
  \item by $\calC^\maximalgroupoid$ for the maximal $\infty$\nobreakdashes-groupoid of an \icategory $\calC$;
  \item by $\dgMod_\basefield$ for the \icategory of chain complexes of $\basefield$\nobreakdashes-modules, obtained by inverting quasi-isomorphisms.
  We will use \emph{cohomological} convention.

  \item \label{relativetensor} If $\calC^{\otimes}$ is a symmetric monoidal \icategory and $A$ an algebra object of $\calC$, we denote by $\Mod_A(\calC)^{\otimes_A}$ the \ioperad of $A$-module objects of \cite[3.3.3.9]{lurie-ha}.
  Under the assumption that $\calC^{\otimes}$ admits geometric realizations and $\otimes$ preserves them in each variable separately, $\Mod_A(\calC)^{\otimes_A}$ is a symmetric monoidal category and $\otimes_A$ coincides with the relative bar construction \cite[ 4.5.2.1]{lurie-ha}\footnote{More specifically, we only need to require the existence of those geometric realizations appearing in relative tensor products -- see \cite[Remark 2.5.14]{zbMATH07506911}}.
\end{notations}

\begin{notations}[Derived Geometry -- see also \crossref{notationderivedgeometry}]{We assume the reader is familiar with derived algebraic geometry at the level of the introductory surveys \cite[\S 1-4]{MR3285853} and \cite[\S 1-4]{Pantev2021}. For conventions:}\label{notationderivedgeometry}
  \item Derived affine rings are modeled by $\cdga$, the \icategory of cohomologically negatively graded commutative differential graded algebras of finite presentation over $\basefield$.

  \item $\dSt_\basefield \subseteq \PSh(\dAff_\basefield, \inftygpd)$ denotes the \icategory of (small or big ) derived stacks over $\basefield$, with étale hyperdescent.

  \item \label{descriptionsofslicecategoriesofstacks}For any $X \in \dSt_\basefield$, we denote by $\dSt_X \coloneqq (\dSt_{\basefield})_{/X}$ the slice \icategory of stacks over $X$. As in \crossref{descriptionsofslicecategoriesofstacks}, we have equivalences
  \[
    \dSt_X \simeq \Sh(\dAff_X, \inftygpd) \simeq \Fun^{\mathcal{R}}( \dSt_X^\op, \inftygpd)
  \]
  where the last denotes the \icategory of limit preserving functors.

  \item For $X \in \dSt_\basefield$ we denote by $\DQCoh(X)$ its derived \icategory of quasi-coherent sheaves \cite[\S 4.1]{MR3285853}.
  If $X=\Spec(A)$ then $\DQCoh(X) = \dgMod_A$.
\end{notations}

\begin{reminders}[Formal stacks and formal completions]{%
    We consider formal derived stacks and formal completions using the formalism of \cite[\S 2.1]{MR3653319}.}
  \label{formalstacks}
  \item\label{remindersonderhamstack} For a derived stack $X$, we denote by $X_{\deRham}$ its de Rham stack with functor of points $X_\deRham(A) \coloneqq X(\pi_0(A)_\red)$.
  By adjunction, we get a functorial morphism $X \to X_\deRham$.
  The de Rham stack plays a crucial role, for its ties to flat connections and $\calD$\nobreakdashes-modules and for its relation to formal completions:
  \begin{remindersinner}
    \item\label{sheavesoverderhamareDmodules} When $X$ is a scheme, $\DQCoh(X_\derham)$ is equivalent to the \icategory of complexes of left D-modules \cite[Part I, Chapter 4]{GR-II}.\item\label{derhamandformalcompletion} Given a map $f \colon Z=\Spec(A) \to X$, the formal completion $\widehat{Z}$ of $X$ along $Z$ coincides with the derived fiber product in $\dSt_\basefield$
    \[
      \begin{tikzcd}
        \widehat{Z} \tikzcart \ar{d} \ar{r} &
        X \ar{d} \\
        Z_{\deRham} \ar{r}[swap]{f_{\deRham}} &
        X_{\deRham}
      \end{tikzcd}
    \]
    (\cf \cite[Prop.\,2.1.8]{MR3653319}). Since the cotangent complex of $X_{\deRham}$ vanishes, the map $\widehat{Z} \to X$ is formally étale.

    \item\label{etaleformalcompletionisnoop} The above specializes to: if $f \colon X \to Y$ is formally étale, then $X \to \widehat X$ is an equivalence: \ie we have a Cartesian square
    \[
      \begin{tikzcd}
        X \tikzcart \ar{d} \ar{r} &
        Y \ar{d} \\
        X_{\deRham} \ar{r}[swap]{f_{\deRham}} &
        Y_{\deRham}.
      \end{tikzcd}
    \]
  \end{remindersinner}
  \item As in \cite[\S 2.1]{MR3653319}, throughout this paper, formal stacks will be seen as a full subcategory of the \icategory of derived stacks $\dSt_\basefield$.

  \item \label{notationderhaminvariantreduced}For a derived stack $X \in \dSt_\basefield$ we set $X_\red $ its reduced substack. When $X=\Spec(A)$ is an affine derived scheme, we have $X_\red=\Spec(A_\red)$.
  The assignment $X\mapsto X_\red$ defines a left adjoint to the functor $(-)_{\deRham}$ and there are canonical equivalences $(X_\red)_\deRham\iso X_\deRham$ and $X_\red \iso (X_\deRham)_\red$.
  See also \crossref{formalstacks}.
\end{reminders}

\begin{notations}[Étale site -- see  \crossref{notationetale} for more details and references.]{}\label{notationetale}
  \item When $X$ is a derived stack, we denote by $X_\et$ the small étale site of $X$, with objects given by maps $S \to X$ which are étale (meaning, representable, locally of finite presentation and the relative cotangent complex $\cotangent_{S/X}$ vanishes).
  \item For a presentable \icategory $\calC$, we denote by
  \[
    \Sh(X_\et,\calC)\subseteq \PSh(X_\et, \calC) \coloneqq \Fun(X_\et^\op, \calC)
  \]
  the full subcategory spanned by \emph{hypercomplete} sheaves.

  \item We denote by $X_\et^\daff \subseteq X_\et$ the full subcategory spanned by étale maps $S \to X$ with $S$ an \emph{affine} derived scheme. Whenever $X$ admits an affine étale atlas, the inclusion $X_\et^\daff \subseteq X_\et$ induces an equivalence on sheaves
  \[
    \Sh(X_\et^\daff, \calC) \iso \Sh(X_\et, \calC).
  \]
\end{notations}

\subsection{Reminders on rigid monoidal structures}

\begin{notations}{}\label{notationsmodulesincategories}
  \item \label{definitionPrlStmonoidalrelative} Let $\PrLStmonoidal[\basefield] \coloneqq \Mod_{\dgMod_{\basefield}}(\PrLSt)^{\otimes}$ denote the symmetric monoidal \icategory of presentable $\basefield$\nobreakdashes-linear \icategories with colimit-preserving functors.
  Following \cite[4.8.1.24]{lurie-ha} it admits relative tensor products in the sense of \cref{relativetensor}.
  \item \label{notationpresentablelinearcategoriesoverring}For a presentable $\basefield$\nobreakdashes-linear symmetric monoidal \icategory $\calC^{\otimes} \in \CAlg(\PrLSt[\basefield])$, we denote by $\PrLStmonoidal[\calC]$ the symmetric monoidal \icategory of $\calC$\nobreakdashes-modules in $\PrLSt[\basefield]$. Whenever $\calC=\dgMod_A$ with $A$ an $\Einfinity$\nobreakdashes-ring over $\basefield$, we shorten the notation to $\PrLStmonoidal[A]$.
\end{notations}

We will require the following basic properties of rigid monoidal categories.
References on the subject -- \cite[I-1 \S9]{GR-I}, \cite[Appendix D.7]{lurie-sag}
or \cite{ramziLocallyRigidinftycategories2024} -- use slightly different definitions.
Those differences do not impact the properties we will really need here.
To fix the necessary conventions, we will rely on \cite[I-1 \S9]{GR-I}.
\begin{reminders}{\label{remindersrigid}Recall from \cite[I-1 \S9]{GR-I}:}
  \item Recall from \cite[I-1 Def.\,9.1.2]{GR-I} that $\calC \in \CAlg(\PrLSt[\basefield])$ is rigid if its monoidal unit is compact and if the right adjoint to the tensor functor $\calC \otimes_\basefield \calC \to \calC$ preserves colimits and is $\calC \otimes_\basefield \calC$-linear.
  \item\label{QCohisrigid} For any cdga $A$, the $\basefield$-linear category $\dgMod_A$ is rigid. Further, for any qcqs derived scheme $S$, the \icategory $\DQCoh(X)$ is rigid -- see \cite[I-3 \S3.7 and \S3.5]{GR-I}.
  \item\label{rigiddualizableagree} If $\calC \in \CAlg(\PrLSt[\basefield])$ is rigid, then $\calM \in \PrLSt[\calC]$ is dualizable over $\calC$ if only if it is dualizable over $\basefield$.
  In that case, the $\calC$- and $\basefield$-linear duals are canonically equivalent. See \cite[I-1 Prop.\,9.4.4]{GR-I}.
\end{reminders}

\section{Sheaves of categories}\label{sectionsheavesofcats}

In order to relate matrix factorizations and the Darboux stack of \cite{hennionholsteinrobaloI}, we will need a suitable stack of which $\MF^\infty$ will provide a section.
This will be the stack of \emph{quasi-coherent hypercomplete presentable 2-periodic categories} in \cref{subsectionsheavesbigcategories}. Later on in \cref{matrixfactorizationasasectionofstack} we construct $\MF^\infty$ as a section of this stack.

This section we will first need to establish some base change formulas for Ind-coherent sheaves, which we describe now in \cref{subsectionremindersonindcoh}.

\subsection{Reminders on Ind-coherent sheaves on derived inf-schemes}
\label{subsectionremindersonindcoh}
In this section we collect several facts about Ind-coherent sheaves on derived inf-schemes, necessary for the analysis of matrix factorization categories in \cref{MFsection}. In particular we prove \cref{basechangeIC-Dmod} and \cref{basechangeIC-redeqetale} below.

\subsubsection{Ind-coherent Sheaves and derived inf-schemes}

\begin{reminders}{We now revise some key notions introduced in \cite{GR-I,GR-II}:}
  \item \label{definitioninfschemes} We denote by $\dInfSch_\basefield$ the \icategory of \emph{derived inf-schemes} in the sense of \cite[I-2 Def.\,3.1.2]{GR-II}: These are derived stacks locally of finite type which admit a deformation theory and whose reduction is a scheme.
  In particular, these include de Rham stacks $X_\derham$ and formal completions $\widehat{X}^Y$ of derived schemes \cite[I-2 Ex.\,3.1.3]{GR-II}. See \cref{formalstacks}.

  \item \label{dinfschemesstableunderfinitelimite} The inclusion $\dInfSch_\basefield\subseteq \dSt_\basefield$ is stable under finite limits \cite[I-2, after Def.\,3.1.2]{GR-II}.

  \item \label{infschemesascolimits} Every derived inf-scheme $\formalV$ can be written as a colimit of derived schemes $\formalV \simeq \colim_\alpha V_\alpha$ where the transition morphisms are closed immersions and reduced equivalences \cite[Part I, Chap 2, Cor.\,4.1.4]{GR-II}.

  \item \label{IndCohassymmetricmonoidalfunctoronindinfschemes}We denote by $\Corr(\dInfSch_\basefield)$ its $(\infty,1)$\nobreakdashes-category of correspondences between derived inf-schemes. By \cite[Part I, Chap. 3, Cor.\,6.1.3]{GR-II} the ind-coherent six-functors formalism yields a \emph{symmetric monoidal} $(\infty,1)$\nobreakdashes-functor
  \[
    \IndCoh \colon \Corr(\dInfSch_\basefield^{\times}) \to \PrLStmonoidal[\basefield],
  \]
  where the tensor structure on the LHS is induced from the Cartesian one on $\dInfSch_\basefield$ by \cite[III-9 \S2.1.3]{GR-I}. By the universal property of correspondences (\cf \cite[Chap. 9, Prop. 3.1.5]{GR-I}), $\IndCoh$ maps a correspondence $\smash{X \from^p Y \to^q Z}$ to
  \[
    \IndCoh(X) \lra^{q_*p^!} \IndCoh(Z).
  \]
  and for every Cartesian square of derived inf-schemes
  \[
    \begin{tikzcd}[row sep=scriptsize]
      X \ar{r}{p} \ar{d}[swap]{q} \tikzcart & W \ar{d}{f} \\
      Y \ar{r}{g} &Z
    \end{tikzcd}
  \]
  we have a base change equivalence
  \[
    p_\ast q^! \simeq f^! g_\ast
  \]
  \item \label{indcohcompactlygenerated} By \cite[Part I, Chap. 3, Cor.\,3.2.2]{GR-II}, for every derived inf-scheme $Y$, the ind-coherent \icategory $\IndCoh(Y)$ is \emph{compactly generated} \cite[5.5.7.1]{lurie-htt} by the image of $i_\ast \colon \IndCoh(Y_\red)^\omega \subseteq \IndCoh(Y_\red) \to \IndCoh(Y)$, where $\IndCoh(Y_\red)^\omega$ denotes the full subcategory of compact objects\rlap{.}\footnote{In this paper, compact will always mean $\omega$-compact.} By adjunction, the functor $i^! \colon \IndCoh(Y) \to \IndCoh(Y_\red)$ is conservative \cite[Part I, Chap 1, Lem.\,5.4.3]{GR-I}.

  \item \label{lowerstarcompact} For any proper map $p \colon Y \to Y'$ of derived inf-schemes, the pushforward functor $p_\ast \colon \IndCoh(Y) \to \IndCoh(Y')$ preserves compact objects.
  Indeed $(p_\red)_\ast \colon \IndCoh(Y_\red) \to \IndCoh(Y'_\red)$ preserves compact objects (\cf \cite[Lem.\,3.3.5 and Cor.\,3.3.6]{gaitsgoryrozenblyum_indcoh}) and we conclude using (d).

  \item \label{uppershriekcompact} For an eventually coconnective morphism $p \colon Y \to Y'$ between derived schemes, the induced pullback functor $p^! \colon \IndCoh(Y') \to \IndCoh(Y)$ preserves compact objects (\cf \cite[Lem.\,7.1.2]{gaitsgoryrozenblyum_indcoh}).
  If $p \colon Y \to Y'$ is schematic and eventually coconnective (between derived inf-schemes), the map $Y \times_Y' Y'_\red \to Y_\red'$ is eventually coconnective.
  It then follows from base-change and (d) that $p^!$ preserves compact objects in that case.

  \item \label{reminderequivalenceICQCderham} Following \cite[Cor.\,5.7.4]{gaitsgoryrozenblyum_indcoh} we have a \emph{symmetric monoidal natural transformation} of \ifunctors on derived schemes almost of finite type
  \[
    \DQCoh^\ast \to \IndCoh^!.
  \]
  Moreover, by \cite[\S2.4.1, Prop.\,2.4.4]{gaitsgoryrozenblyum:crystals}, for any derived stack $X$ this natural transformation yields an equivalence $\DQCoh(X_\derham) \simeq \IndCoh(X_\derham)$.

  \item \label{basechangeQcohIndCoh} By \cite[Prop.\,7.5.7]{gaitsgoryrozenblyum_indcoh}, if $f \colon Y \to Y'$ is a smooth map of derived schemes, the functor $f^! \colon \IndCoh(Y') \to \IndCoh(Y)$ induces an equivalence
   \[
    \DQCoh(Y) \otimes_{\DQCoh(Y')} \IndCoh(Y') \simeq \IndCoh(Y)
  \]
\end{reminders}

\begin{lemma}\label{indcohhyperdescent}
  The functor $\IndCoh \colon \dAff_\basefield^\op \to \PrLSt[\basefield]$ (with $!$\nobreakdashes-pullbacks) satisfies smooth hyperdescent.
  \begin{proof}
    This is a simply adaptation of the proof of \cite[Prop.\,8.3.6]{gaitsgoryrozenblyum_indcoh} to hyperdescent.
    Fix $S_\bullet \colon \Delta_+^\op \to \dAff_\basefield$ a smooth hypercovering.
    By \cref{basechangeQcohIndCoh}, the coaugmented cosimplicial diagram $\IndCoh(S_\bullet)$ with $\IndCoh^!$-functoriality is equivalent to the diagram
    \[
      \DQCoh(S_\bullet) \otimes_{\DQCoh(S_{-1})} \IndCoh(S_{-1}) \colon \Delta_+ \to \PrLSt[\basefield].
    \]
    where $\DQCoh(S_\bullet) $ has $\DQCoh^\ast$-functoriality.
    As $\DQCoh(S_{-1})$ is rigid over $\basefield$ and by \cite[Part I, Chap. 3, \S 6.2]{GR-II} $\IndCoh(S_{-1})$ is self-dualizable over $\basefield$, it follows from \cref{rigiddualizableagree} that $\IndCoh(S_{-1})$ is dualizable over $\DQCoh(S_{-1})$.
    In particular, the functor $-\otimes_{\DQCoh(S_{-1})} \IndCoh(S_{-1}) $ preserves limits, and the result thus follows from hyperdescent for $\DQCoh$ \cite[Cor.\,6.13]{lurie:dagvii}.
  \end{proof}
\end{lemma}

\subsubsection{Base change of categories of Ind-Coherent sheaves}

\begin{lemma}\label{basechangeIC-Dmod}
  Let $\formalU \to \formalV \in \dInfSch_\basefield$ and denote by $\calX$ the formal neighborhood of $\formalU$ in $\formalV$ defined as the fiber product in \cref{derhamandformalcompletion} . The induced functor
  \[
    \IndCoh(\formalV) \otimes_{\IndCoh(\formalV_\derham)} \IndCoh(\formalU_\derham) \to \IndCoh(\calX)
  \]
  is an equivalence in $\PrLSt[\basefield]$.
  \begin{proof}
    By \cite[formula (2.9)]{beraldo:center}, the result holds if $\formalV$ is a derived scheme.
    For the general case of a derived inf-scheme $\formalV$, as in \cref{infschemesascolimits}, write $\formalV$ as a colimit of derived schemes $\formalV \simeq \colim_\alpha V_\alpha$ where the transition morphisms are closed immersions and reduced equivalences.
    Using \cite[Part I - Chap.\,3 - \S4, Prop.\,4.1.3 and Cor.\,4.3.5]{GR-II}, we have $\IndCoh(\formalV) \simeq \colim \IndCoh_*(V_\alpha)$ in $\PrLSt[\basefield]$, and this description is compatible with the action of $\IndCoh(\formalV_\derham) \simeq \IndCoh(V_{\alpha,\derham})$ by base change (\cf \cref{IndCohassymmetricmonoidalfunctoronindinfschemes}) along the Cartesian square
    \[
      \begin{tikzcd}[row sep=scriptsize]
        \formalV_\alpha \tikzcart \ar{r}\ar{d}[swap]{\alpha\to\beta}& \formalV_\alpha\times \formalV_\derham \ar{d}{(\alpha\to\beta)\times \Id_{\formalV_\derham}}\\
        \formalV_\beta \ar{r}& \formalV_\beta\times \formalV_\derham\rlap{.}
      \end{tikzcd}
    \]
    The result follows.
  \end{proof}
\end{lemma}

\begin{corollary}\label{basechangeIC-redeqetale}
  Let $\formalU \to \formalV \from \calW \in \dInfSch_\basefield$ and denote by $\calX$ the fiber product $\formalU \times_\formalV \calW$\footnote{By \cref{dinfschemesstableunderfinitelimite} this is again a derived inf-scheme.}.
  Assume that the morphism $\formalU \to \formalV$ is a reduced equivalence and that $\calW \to \formalV$ is formally étale.
  The induced functor
  \[
    \IndCoh(\formalU) \otimes_{\IndCoh(\formalV)} \IndCoh(\calW) \to \IndCoh(\calX)
  \]
  is an equivalence.
\end{corollary}
\begin{proof}
  As $\calW \to \formalV$ is formally étale, as in \cref{etaleformalcompletionisnoop}, we have Cartesian squares
  \[
    \begin{tikzcd}[row sep=small]
      \calX \ar{r} \ar{d} \tikzcart & \calW \ar{d} \ar{r} \tikzcart & \calW_\derham \ar{d} \\
      \calU \ar{r} & \calV \ar{r} & \calV_\derham\mathrlap{{} \simeq \formalU_\derham.}
    \end{tikzcd}
  \]
  Using \cref{basechangeIC-Dmod} on both the right and outer squares of this diagram, we get
  \begin{multline*}
    \IndCoh(\formalU) \otimes_{\IndCoh(\formalV)} \IndCoh(\calW)
    \simeq \IndCoh(\formalU) \otimes_{\IndCoh(\formalV)} \left(\IndCoh(\formalV) \otimes_{\IndCoh(\formalV_\derham)} \IndCoh(\calW_\derham)\right)
    \\
    \simeq \IndCoh(\formalU)\otimes_{\IndCoh(\formalV_\derham)} \IndCoh(\calW_\derham)
    \simeq \IndCoh(\calX).\qedhere
  \end{multline*}
\end{proof}

\begin{corollary}\label{basechangeIC-overDmodgeneral}
  Let $\formalU \in \dInfSch_\basefield$ and $\formalV \to \formalU_\derham \in \dInfSch_\basefield$.
  The box product functor
  \[
    \IndCoh(\formalU) \otimes_{\IndCoh(\formalU_\derham)} \IndCoh(\formalV) \to \IndCoh\left(\formalU \times_{\formalU_\derham} \formalV \right)
  \]
  is an equivalence.
\end{corollary}
\begin{proof}
  Observe first the equivalence $\formalU \times_{\formalU_\derham} \formalV \simeq \formalV_\derham \times_{\formalU_\derham \times \formalV_\derham} (\formalU \times \formalV)$.
  We apply \cref{basechangeIC-Dmod} to the morphism $\formalV_\red \to \formalU \times \formalV$ and find
  \[
    \IndCoh\left(\formalU \times_{\formalU_\derham} \formalV\right) \simeq \IndCoh(\formalU \times \formalV) \otimes_{\IndCoh(\formalU_\derham \times \formalV_\derham)} \IndCoh\left(\formalV_\derham \right).
  \]
  The result follows as $\IndCoh$ is symmetric monoidal over $\basefield$ (\cf \cref{IndCohassymmetricmonoidalfunctoronindinfschemes}).
\end{proof}

\subsection{Sheaves of 2-periodic presentable dg-categories and 1-affineness.}
\label{subsectionsheavesbigcategories}

\begin{reminders}{}
  \item \label{definitionofhypercompletemodule} Let $A \in \cdga$ and $\calC \in \PrLSt[A]$. Following \cite[Def.\,6.8]{lurie:dagvii}\footnote{See also \cite[Def.\,2.2]{zbMATH06084026}.} we say that $\calC$ is (étale) \emph{hypercomplete} if the functor $\cdga[A] \to \PrLSt$ mapping $B$ to $\calC \otimes_A B \coloneqq \calC \otimes_{\dgMod_A} \dgMod_B$ satisfies étale hyperdescent. We denote by $\PrLStHyp[A] \subset \PrLSt[A]$ the full subcategory of hypercomplete $A$\nobreakdashes-linear \icategories.

  \item \label{dualizableishypercomplete} It follows from \cite[Cor.\,6.13]{lurie:dagvii} that every \emph{dualizable} object $\calC \in \PrLSt[A]$ is hypercomplete.
  In particular, for any $B \in \cdgaunbounded_A$  the $A$-linear category $\dgMod_B \in \PrLSt[A]$ is hypercomplete. We denote by $\PrLStdual[A] \subset \PrLSt[A]$ the full subcategory of dualizable $A$\nobreakdashes-linear \icategories. By definition it is closed under tensor products.

  \item \label{algebraobjectsinhypercompletecategories} The full subcategory $\PrLStHyp[A]$ is not necessarily closed under tensor products in $\PrLSt[A]$.
  Nevertheless, we get a sub-$\infty$-operad $\PrLStHypOperad[A] \subseteq \PrLStmonoidal[A]$.

  \item \label{hypercompletestableunderlimits} By \cite[Rem.\,7.3]{lurie:dagxi}, hypercomplete $A$-linear \icategories are stable under limits.
  \item \label{cechdescentsheavescats} The \ifunctor $A \mapsto \PrLSt[A]$
  satisfies étale \v Cech descent (\cite[Thm.\,5.13]{lurie:dagxi}).
  \item \label{hyperdescentsheavescats} Hypercomplete \icategories are stable under base change. Moreover, the induced functor with values in big \icategories
  \begin{equation}\label{functorhypercompletemodules}
    \cdga \to \inftycats \hspace{1cm} A \mapsto \PrLStHyp[A]
  \end{equation}
  satisfies étale hyperdescent as a straightforward adaptation of the proof of \cite[Thm.\,7.5]{lurie:dagxi}).
\end{reminders}

\begin{notations}{Throughout this paper we will describe $2$\nobreakdashes-periodic \icategories as linear \icategories over the free algebra with an invertible element of cohomological degree $2$: }
  \item \label{notationbeta} We denote by $\kbeta$ the free cdga over $\basefield$ with a generator $\beta$ in \emph{cohomological} degree $2$.
  As an object of $\dgMod_\basefield$, we have $\kbeta \simeq \bigoplus_{n\geq 0} \basefield[2n] \simeq \prod_{n\geq 0} \basefield[2n]$.

  \item \label{notationbetainversed} We denote by $\kbetaLaurent$ the cdga obtained from $\kbeta$ by inverting $\beta$.
  As an object of $\dgMod_\basefield$, we have $\kbetaLaurent \simeq \bigoplus_{n \in \integers}\basefield[2n]$.

  The \icategory of \emph{presentable 2-periodic dg-categories} is the category of modules $ \PrLSt[\kbetaLaurent]$.
  \item\label{notationsmodulesincategoriesbeta}For any presentable $\basefield$\nobreakdashes-linear \icategory $\calC$, we set $\calC \ofbeta \coloneqq \calC \otimes_\basefield \dgMod_{\kbeta} \in \PrLSt[\kbeta]$ and $\calC\ofbetaLaurent \coloneqq \calC \otimes_\basefield \dgMod_{\kbetaLaurent} \in \PrLSt[\kbetaLaurent]$.
  In particular, for any $S \in \dInfSch_\basefield$, we have the symmetric monoidal \icategories $\IndCoh(S)\ofbeta$ and $\IndCoh(S)\ofbetaLaurent$ of $2$-periodic ind-coherent complexes.
\end{notations}

We now discuss the notion of hypercompleteness in the non-affine case:

\begin{definitions}{Let $Y$ be a derived stack.}\label{notionofhypercompletenonaffine}
  \item \label{hypercompletenonaffinecase}
  A $\DQCoh(Y)$\nobreakdashes-module $\calC \in \PrLSt[\DQCoh(Y)]$ is called \emph{hypercomplete} if for any $\Spec A \to X$, the pulled back category $\calC \otimes_{\DQCoh(Y)} \dgMod_A \in \PrLSt[A]$ is hypercomplete in the sense of \cref{definitionofhypercompletemodule}.
  We denote by $\PrLStHyp[\DQCoh(Y)]$ the full subcategory of $\PrLSt[\DQCoh(Y)]$ spanned by hypercomplete modules.

  \item \label{hypercompletenonaffinecase2periodic} A $2$\nobreakdashes-periodic $\DQCoh(Y)$\nobreakdashes-module is a $\DQCoh(Y)\ofbetaLaurent$\nobreakdashes-module.
  We say it is \emph{hypercomplete} if the underlying $\DQCoh(Y)$\nobreakdashes-module is hypercomplete in the sense of \cref{hypercompletenonaffinecase}.
  We denote by
  \[
    \PrLStHyp[\DQCoh(Y)\ofbetaLaurent] \subset \PrLSt[\DQCoh(Y)\ofbetaLaurent]
  \]
  the full subcategory of hypercomplete $2$\nobreakdashes-periodic modules.
  In the case where $Y = \Spec A$ is affine, we set $\PrLStHyp[A\ofbetaLaurent] \coloneqq \PrLStHyp[\DQCoh(Y)\ofbetaLaurent]$.
  \item \label{notationdualizable2permodules}We denote by
  \[
    \PrLStdual[\DQCoh(Y)\ofbetaLaurent] \subset \PrLStHyp[\DQCoh(Y)\ofbetaLaurent] \subset \PrLSt[\DQCoh(Y)\ofbetaLaurent]
  \]
  the full subcategory of dualizable modules (note that by \cref{dualizableishypercomplete}, every dualizable module is hypercomplete).
  In the case where $Y = \Spec A$ is affine, we set $\PrLStdual[A\ofbetaLaurent] \coloneqq \PrLStdual[\DQCoh(Y)\ofbetaLaurent]$.
  \item \label{notationforinvertible} Denote by
   \[
  \PrLStinvertible[\DQCoh(Y)\ofbetaLaurent] \subset \PrLStdual[\DQCoh(Y)\ofbetaLaurent]
  \]
  the full symmetric monoidal subcategory spanned by those dualizable modules that are invertible.
\end{definitions}

\begin{lemma}\label{twoperdgcatdescent}
  The functors $\cdga \to \inftycats$ with values in big \icategories mapping $A$ to $\PrLStHyp[A\ofbetaLaurent]$ and to $\PrLStdual[A\ofbetaLaurent]$ satisfy étale hyperdescent.
  The functor mapping $A$ to $\PrLSt[A\ofbetaLaurent]$ satisfies \v Cech étale descent.
\end{lemma}
\begin{proof}
  Let us start by observing that using \cref{algebraobjectsinhypercompletecategories}, the \ifunctor \cref{functorhypercompletemodules} lifts to an \ifunctor with values in the \icategory of $\infty$-operads $\inftyoperads$
  \begin{equation}\label{sheafvalueioperads}
    \cdga \to \inftyoperads \hspace{1cm} A \mapsto \PrLStHypOperad[A]
  \end{equation}
  satisfying étale hyperdescent (\cf \cref{hyperdescentsheavescats}).
  Now for any cdga $A$, the $A$-linear \icategory $\dgMod_{A\ofbetaLaurent}$ is dualizable \footnote{By \cref{rigiddualizableagree} and the fact it is compactly generated.} and hence hypercomplete (\cf \cref{dualizableishypercomplete}).
  The functor $\calA \colon A \mapsto \dgMod_{A\ofbetaLaurent}$ is a commutative algebra section of the stack \cref{sheafvalueioperads}.
  The result then follows from \cref{hyperdescentsheavescats,cechdescentsheavescats} by taking the stacks of $\calA$\nobreakdashes-modules, and from the fact that locally dualizable modules are globally dualizable.
\end{proof}

\begin{definitions}{}\item \label{definitionofgoodstackofsheavesofcategories}
  We denote by $\shcatsqcohtwoper$ the (limit preserving) functor obtained from \cref{hypercompletenonaffinecase,twoperdgcatdescent}:
  \begin{align*}
    \shcatsqcohtwoper \colon \dSt_\basefield^\op & \to \inftycats \\ Y &\mapsto \lim_{\Spec A \to Y} \PrLStdual[A\ofbetaLaurent].
  \end{align*}
\item  \label{definitionofgoodstackofsheavesofcategoriesinvertible} Similarly, we denote by $\Azumayatwoper$ the full substack of $\shcatsqcohtwoper$ spanned by  invertible modules (\cf \cref{notationforinvertible}). Notice that by definition of invertible objects, this has a natural structure of abelian group stack.
\end{definitions}

\begin{reminders}{1-affineness:}\label{1affinessreminder}
  \item \label{definition1affine} Let $Y$ be a derived stack. Following \cite{gaitsgory_affineness}, we say that $Y$ is \emph{$1$\nobreakdashes-affine} if the natural pullback functor $\PrLSt[\DQCoh(Y)] \to \lim_{\Spec A \to Y} \PrLSt[A]$ is an equivalence.
  \item \label{YdR1affine} By \cite[Thm.\,2.6.3]{gaitsgory_affineness}, the de Rham stack of any derived inf-scheme is $1$\nobreakdashes-affine.
  \item \label{qcqsscheme} By \cite[Thm.\,2.1.1]{gaitsgory_affineness} every quasi-compact quasi-separated derived scheme is $1$\nobreakdashes-affine.
  \item \label{product1affine} By \cite[Cor.\,3.2.8]{gaitsgory_affineness} the Cartesian product of 1-affine derived stacks is 1-affine.
\end{reminders}

\begin{lemma}\label{1affinebeta}
  Let $Y$ be a $1$\nobreakdashes-affine derived stack. The pullback functors
  \begin{enumerate}[label={\textrm{\upshape{(\alph*)}}}]
    \item $\PrLStHyp[\DQCoh(Y)] \to \lim_{\Spec A \to Y} \PrLStHyp[A]$,
    \item $\PrLSt[\DQCoh(Y)\ofbetaLaurent] \to \lim_{\Spec A \to Y} \PrLSt[A\ofbetaLaurent]$,
    \item $\PrLStHyp[\DQCoh(Y)\ofbetaLaurent] \to \lim_{\Spec A \to Y} \PrLStHyp[A\ofbetaLaurent]$ and
    \item \label{1affinebetadualizable}$\PrLStdual[\DQCoh(Y)\ofbetaLaurent] \to \shcatsqcohtwoper(Y) \coloneqq \lim_{\Spec A \to Y} \PrLStdual[A\ofbetaLaurent]$
    \item \label{1affinebetainvertible} $\PrLStinvertible[\DQCoh(Y)\ofbetaLaurent] \to \Azumayatwoper(Y) \coloneqq \lim_{\Spec A \to Y} \PrLStinvertible[A\ofbetaLaurent]$
  \end{enumerate}
  are equivalences.
\end{lemma}
\begin{proof}
  The assertion (a) stems from \cref{hypercompletenonaffinecase} while (b) and (c) follow by taking modules in (a).
  Assertion (d) follows from \cite[Cor.\,1.4.3 and Prop.\,1.4.5]{gaitsgory_affineness} and (e) follows from (d).
\end{proof}

\begin{observation}\label{observationwhichonesaresymmetricmonoidal}
	Notice that (b), (d) and (e) in \cref{1affinebeta} are symmetric monoidal equivalences.
\end{observation}

\begin{lemma}\label{smoothoverregularisregular}
  Let $Y$ be a $1$-affine derived stack and let $\phi \colon W \to Y$ be a qcqs schematic morphism.
  Then $W$ is $1$-affine. If moreover $\phi$ is smooth, then the canonical functor
  \[
    \IndCoh(Y) \otimes_{\DQCoh(Y)} \DQCoh(W) \to \IndCoh(W)
  \]
  is an equivalence.
\end{lemma}
\begin{proof}
  $1$-affineness follows from \cite[Cor.\,3.2.7]{gaitsgory_affineness}.
  For any derived affine scheme $S$ over $Y$, we set $W_S \coloneqq W \times_Y S$.
  The stack $W_S$ is a qcqs derived scheme and the projection $W_S \to S$ is smooth by assumption.
  Using \cref{basechangeQcohIndCoh} we get a functorial equivalence $\IndCoh(W_S) \simeq \IndCoh(S) \otimes_{\DQCoh(S)} \DQCoh(W_S)$.
  As $Y$ is $1$-affine, we have $\DQCoh(W_S) \simeq \DQCoh(W) \otimes_{\DQCoh(Y)} \DQCoh(S)$ by \cite[Prop.\,3.1.9]{gaitsgory_affineness}.
  We find
  \[
    \IndCoh(W_S) \simeq \IndCoh(S) \otimes_{\DQCoh(Y)} \DQCoh(W).
  \]
  Now, \cref{remindersrigid} imply $\DQCoh(W_S)$ is dualizable over $\DQCoh(S)$ for any $S \to Y$ affine.
  It follows from \cite[Prop.\,1.4.5 and Cor.\,1.4.3]{gaitsgory_affineness} that $\DQCoh(W)$ is dualizable over $\DQCoh(Y)$.
  As a consequence, writing $W \simeq \colim_{S \to Y} W_S$ (with $S$ affine), we get using $!$-pullbacks
  \begin{multline*}
    \IndCoh(W) \simeq \lim_{S \to Y} \IndCoh(W_S) \simeq \lim_{S \to Y} \left(\IndCoh(S) \otimes_{\DQCoh(Y)} \DQCoh(W)\right)
    \\
    \simeq \left(\lim_{S \to Y} \IndCoh(S)\right) \otimes_{\DQCoh(Y)} \DQCoh(W) \simeq \IndCoh(Y) \otimes_{\DQCoh(Y)} \DQCoh(W).\qedhere
  \end{multline*}
\end{proof}

\begin{corollary}\label{basechangeisconservative}
  Let $Y$ be a $1$-affine derived stack and let $\phi \colon W \twoheadrightarrow Y$ be a schematic epimorphism of stacks.
  The base-change functor
  \[
    - \otimes_{\DQCoh(Y)\ofbetaLaurent} \DQCoh(W)\ofbetaLaurent \colon \PrLStHyp[\DQCoh(Y)\ofbetaLaurent] \to \PrLStHyp[\DQCoh(W)\ofbetaLaurent]
  \]
  is conservative.
  If moreover $\phi$ is smooth, then the base-change functor
  \[
    - \otimes_{\IndCoh(Y)\ofbetaLaurent} \IndCoh(W)\ofbetaLaurent \colon \PrLStHyp[\IndCoh(Y)\ofbetaLaurent] \to \PrLStHyp[\IndCoh(W)\ofbetaLaurent]
  \]
  is conservative as well.
\end{corollary}
\begin{proof}
  It follows straightforwardly from \cref{1affinebeta,smoothoverregularisregular} and the fact that $\shcatsqcohtwoper$ (\cf \cref{definitionofgoodstackofsheavesofcategories}) is a stack.
\end{proof}

\section{Matrix factorizations}
\label{MFsection}

In this \cref{sectiondefinitionMF} we extend the construction of matrix factorization categories from  \cite{MR3121870} to the setting of derived inf-schemes. In (\cref{subsectionrelativethomsebastiani})  we  prove a relative version of the Thom-Sebastiani theorem, which we use in \cref{subsection-actionquadraticbundles} to construct the action of quadratic bundles and Kn\"orrer periodicity.

\subsection{Big \texorpdfstring{($\MF^\infty$) }{}and small \texorpdfstring{($\MF$) }{}matrix factorizations on inf-schemes}\label{sectiondefinitionMF}
The goal of this section is to extend the construction of the category of matrix factorizations of \cite{MR3121870} to derived inf-schemes. We recall this construction and establish its functorialities.

\begin{notations}{}
  \item We denote by $\LGdinf \coloneqq \dInfSch_{/\affineline{}}$ the \icategory of derived inf-schemes $\formalU$ over $\basefield$ equipped with a function $f \colon \formalU \to \affineline{\basefield}$.
  \item \label{notationforB} Let $\bbB \in \dAff_\basefield$ be the derived abelian\footnote{Indeed, because $(\affineline{},+)$ is abelian, Eckmann–Hilton implies $\bbB$ also is.} affine group scheme $\bbB \coloneqq \Omega_0\affineline{} = \Spec(\basefield[\varepsilon])$ (with $\varepsilon$ in cohomological degree $-1$).
\end{notations}

\begin{construction}\label{definitionoffunctorPreMFoninfschemes} Combining \cite[5.2.6.12, 5.2.6.15, 5.2.6.28, 5.2.6.29]{lurie-ha} and the fact $\bbB$ is abelian, the extraction of the derived fiber at zero gives a functor with values in left modules \cite[4.5.1.4]{lurie-ha}:
  \[
    \zerolocus\colon \dSt_{/\affineline{}} \to \Mod_{\bbB}(\dSt_\basefield)\hspace{1cm} ( U \to \affineline{}) \mapsto U_0 \coloneqq U\fiberproduct{\affineline{}}0
  \]
  which, by \cref{dinfschemesstableunderfinitelimite}, restricts to an \ifunctor
  \[
    \zerolocus\colon \dInfSch_{/\affineline{}} \to \Mod_{\bbB}(\dInfSch_\basefield).
  \]
  On the other hand, the symmetric monoidal inclusion
  \[
    \hrm \colon \dInfSch_\basefield^{\times} \hookrightarrow \Corr\left(\dInfSch_\basefield^{\times}\right) \hspace{2cm} (X \to Y) \mapsto
    (X \from^\Id X \to Y)
  \]
  maps the abelian group object $\bbB$ to a commutative algebra object in correspondences $\hrm(\bbB)$.
  We get:
  \[
    \begin{tikzcd}[column sep=normal]
      \LGdinf\ar{r}{\zerolocus}& \Mod_{\bbB}(\dInfSch_\basefield)\ar{r}{\hrm} & \Mod_{\hrm(\bbB)}\left(\Corr\left(\dInfSch_\basefield\right)\right).
    \end{tikzcd}
  \]
  We now observe that $\hrm \circ\, \zerolocus$ satisfies all the Beck-Chevalley prerequisites of the universal property of correspondences \cite[Part III-Chap.\,7 Thm.\,3.2.2]{GR-I}:
  \begin{enumerate}
    \item For every map $i \colon U \to V$ in $\LGdinf$ the induced map of $\bbB$-modules $U_0 \to V_0$ seen as a map of correspondences via $\hrm$, admits an $\hrm(\bbB)$-linear adjoint, given by the opposite correspondence
          \[
            \begin{tikzcd} V_0 & U_0 \ar{l} \ar{r}{\Id} & U_0\rlap{.}\end{tikzcd}
          \]
          The $\hrm(\bbB)$-linear structure of this adjoint stems from the Cartesian square
          \[
            \begin{tikzcd}
              \tikzcart U_0 \times \bbB \ar{r}{\Id\times i} \ar{d}[swap]{\text{action}} & V_0 \times \bbB \ar{d}{\text{action}}
              \\
              U_0 \ar{r}{i} & V_0\rlap{.}
            \end{tikzcd}
          \]
    \item The fact that the Beck-Chevalley transformation associated to a Cartesian square in $\LGdinf$ is an equivalence, is an immediate consequence of the fact $\zerolocus$ preserves fiber products.
  \end{enumerate}

  Therefore, we obtain a canonical extension
  \[
    \begin{tikzcd}[row sep=scriptsize]
      \LGdinf\ar{r}{\zerolocus}\ar[hookrightarrow]{d} & \Mod_{\bbB}(\dInfSch_\basefield)\ar{r}{\hrm} & \Mod_{\hrm(\bbB)}\left(\Corr\left(\dInfSch_\basefield\right)\right)\rlap{.}
      \\
      \Corr\left(\LGdinf\right) \ar[dashed, out=0, in=-160]{urr}
    \end{tikzcd}
  \]
  Finally, by composing with the symmetric monoidal functor $\IndCoh$ of \cref{IndCohassymmetricmonoidalfunctoronindinfschemes}, we define the \ifunctor of pre-matrix factorizations as:
  \[
    \begin{tikzcd}[column sep=small]
      \PreMF \colon \Corr\left(\LGdinf\right) \ar{r} & \Mod_{\hrm(\bbB)}\left(\Corr\left(\dInfSch_\basefield\right)\right) \ar{r}{\IndCoh} & \Mod_{\IndCoh(\bbB)}\left(\PrLSt[\basefield]\right)\rlap{,}
    \end{tikzcd}
  \]
  where the symmetric monoidal structure on $\IndCoh(\bbB)$ is a convolution given by the pushforward along the multiplication $\bbB \times \bbB \to \bbB$
  \[
    \IndCoh(\bbB) \otimes \IndCoh(\bbB) \simeq \IndCoh(\bbB \times \bbB) \to \IndCoh(\bbB)\rlap{.}
  \]
\end{construction}

\begin{lemma}[{see \cite[Prop.\,3.1.4]{MR3121870}}]\label{ICBiskbeta}
  The \icategory $\IndCoh(\bbB)$ equipped with its convolution monoidal structure is equivalent to $\dgMod_{\kbeta}^{\otimes}$ with the $\kbeta$\nobreakdashes-linear tensor structure of \cref{notationbeta}.
  In particular we have a symmetric monoidal equivalence with respect to the relative tensor products (\cref{relativetensor}):
  \[
    \Mod_{\IndCoh(\bbB)}\left(\PrLSt[\basefield]\right)^{\otimes_{\IndCoh(\bbB)}} \simeq \PrLStmonoidal[\kbeta].
  \]
\end{lemma}

\begin{definition}\label{definitionoffunctorMFfunctoroninfschemes}
  We define the \ifunctor of \emph{matrix factorizations} as the composition
  \[
    \MF^\infty \colon
    \begin{tikzcd}[column sep=2.6em]
      \Corr\left( \LGdinf\right) \ar{r}{\PreMF} & \Mod_{\IndCoh(\bbB)}\left(\PrLSt[\basefield]\right) \simeq \PrLSt[\kbeta] \ar{rr}{- \otimes_\kbeta \kbetaLaurent} && \PrLSt[\kbetaLaurent].
    \end{tikzcd}
  \]
\end{definition}

\begin{example}
  \label{MFzerofunction}
  Let $S \in \dInfSch_\basefield$ equipped with the zero function. Then, $\zerolocus(S,0) \simeq S\times \bbB$ and by \cref{IndCohassymmetricmonoidalfunctoronindinfschemes} we have an equivalence of \icategories
  \[
    \MF^\infty(S,0)=\IndCoh(S\times \bbB)\otimes_{\IndCoh(\bbB)}\kbetaLaurent \simeq \left(\IndCoh(S)\otimes_{\basefield}\IndCoh(\bbB)
    \right)\otimes_{\IndCoh(\bbB)}\kbetaLaurent \simeq \IndCoh(S)\ofbetaLaurent
  \]
  (\cf \cref{notationsmodulesincategoriesbeta}).
\end{example}

\begin{remark}
  We will construct in \cref{sectionconvolutionproductsaxioms} a symmetric monoidal enhancement of the functor $\MF^\infty$ -- \cf \cref{thomsebastianiMF}.
\end{remark}

\begin{warning}\label{warningsinglecriticalvalue}
  The construction of $\MF^\infty$ in \cref{definitionoffunctorMFfunctoroninfschemes} is only partially the right one, as it captures only what happens over the zero fiber.
  To obtain a correct definition of $\MF^\infty$ one should really consider a direct sum of the different matrix factorization categories $\MF^\infty(\formalU, f-\lambda)$ for all critical values $\lambda$ of $f$. See for instance \cite[Thm.\,4.1.3-(iii)]{MR3121870}.
  We will however focus on the case where the sole critical value is $0$.
\end{warning}

\begin{notation}\label{notationLGinfzero}
  We denote by $\LGdinfzero\subseteq \LGdinf$ the full sub \icategory spanned by functions $f$ vanishing on the reduced critical locus $\Crit(f)_\red$, \ie $f_{|_{\Crit(f)_\red}}=0$. In other words, zero is the only critical value. From now on, we will only consider the functor of matrix factorizations $\MF^\infty$ of \cref{definitionoffunctorMFfunctoroninfschemes} restricted to $\LGdinfzero$. See \crossref{sectionfixingexactstructureequalsconstant} for the relation with fixing exact structures on the derived critical loci.
\end{notation}

\begin{lemma}\label{compactgenerationMF}
  Let $(\formalU,f) \in \LGdinfzero$. Then the \icategory $\MF^\infty(\formalU,f)$ is compactly generated.
  \begin{proof}
    \cref{indcohcompactlygenerated} implies the \icategory $\PreMF(\formalU,f) \coloneqq \IndCoh( \zerolocus(\formalU,f))$ of \cref{definitionoffunctorPreMFoninfschemes} is \emph{compactly generated} by the image of $i_\ast \colon \IndCoh(\zerolocus(\formalU,f)_\red)^\omega \subseteq \IndCoh(\zerolocus(\formalU,f)_\red) \to \IndCoh(\zerolocus(\formalU,f))$.
    Since both the group multiplication map $\bbB\times \bbB \to \bbB$ and the action $\bbB\times \zerolocus(\formalU,f) \to \zerolocus(\formalU,f)$ are proper maps (obtained by pullback from the closed immersion of $0$ in $\affineline{}$), by \cref{lowerstarcompact} the induced pushforward functors $\IndCoh(\bbB) \otimes \IndCoh(\bbB) \simeq \IndCoh(\bbB \times \bbB) \to \IndCoh(\bbB)$ and
    \[
      \IndCoh(\bbB) \otimes \PreMF(\formalU,f) \simeq \IndCoh(\bbB\times \zerolocus(\formalU,f)) \to \PreMF(\formalU,f) \coloneqq \IndCoh(\zerolocus(\formalU,f))
    \]
    both preserve compact objects.
    As the functor $\IndCoh(\bbB) \simeq \dgMod_{\kbeta} \to \dgMod_{\kbetaLaurent}$ also preserves compact objects, we conclude using \cite[I-1 Cor.\,8.7.4]{GR-I}.
  \end{proof}
\end{lemma}

\begin{definition}\label{smallMF}
  In view of \cref{compactgenerationMF},
  for $(\formalU, f) \in \LGdinfzero$ , we set the \emph{small category of matrix factorizations} $\MF(\formalU, f) \subset \MF^\infty(\formalU, f)$ to be the full subcategory of compact objects.
\end{definition}

\begin{observations}{}
  \item Combining \cite{MR2101296} and \cite{MR3121870}, if $U$ is a smooth affine scheme, the category $\MF(U,f)$ of \cref{smallMF} is equivalent to the more hands-on construction as the idempotent completion of the 2-periodic dg-category whose objects are diagrams $r \colon E \rightleftarrows F :s$ with $E, F$ finitely generated projective modules on $U$ and $rs = f\Id_F$ and $sr = f\Id_E$ and whose complexes of morphisms are the natural ones.

  \item\label{rmkfunctorialitiesMF} \cref{definitionoffunctorPreMFoninfschemes} and \cref{definitionoffunctorMFfunctoroninfschemes} establish various functorialities: for any $\phi \colon (\formalU,f) \to (\calV,g)$, we have functors
  \[
    \begin{tikzcd}
      \MF^\infty(\calV,g) \ar[shift left]{r}{\phi^!} & \MF^\infty(\formalU,f) \ar[shift left]{l}{\phi_*}
    \end{tikzcd}
  \]
  satisfying base change.
  Using \cref{uppershriekcompact}, we see that $\phi^!$ restricts to a functor $\phi^! \colon \MF(\calV,g) \to \MF(\formalU,f)$ if $\phi$ is schematic and eventually coconnective.
  Assuming on the other hand $\phi$ to be proper, we get $\phi_* \colon \MF(\formalU,f) \to \MF(\calV,g)$ by \cref{lowerstarcompact}.
\end{observations}

\begin{lemma}\label{MFdescent}
  The restrictions
  \[
    \begin{tikzcd}[row sep=0]
      \LGdinfzero^\op \ar{r} & \Corr\left(\LGdinfzero\right) \ar{r}{\PreMF} & \PrLSt[\kbeta] \rlap{\hspace{1em} and}\\
      \LGdinfzero^\op \ar{r} & \Corr\left(\LGdinfzero\right) \ar{r}{\MF^\infty} & \PrLSt[\kbetaLaurent]
    \end{tikzcd}
  \]
  satisfy \v Cech $h$\nobreakdashes-descent (and thus \v Cech fppf descent as well) and smooth hyperdescent.

  \begin{proof}
    By \cref{indcohhyperdescent}, the functor $\IndCoh$ satisfies smooth hyperdescent.
    Moreover, by \cite[I-3 Prop.\,3.3.3]{GR-II}, it satisfies \v{C}ech descent along proper covers of derived inf-schemes.
    It follows that $\IndCoh$ satisfies \v Cech $h$\nobreakdashes-descent for derived inf-schemes, and thus so does $\PreMF$.
    Finally, the functor $\MF^\infty \coloneqq \PreMF \otimes_{\kbeta} \kbetaLaurent$ then also satisfies both \v Cech $h$\nobreakdashes-descent and smooth hyperdescent, since $\dgMod_{\kbetaLaurent}$ is dualizable in $\PrLSt[\kbeta]$ (\cf \cref{remindersrigid}).
  \end{proof}
\end{lemma}

\subsection{Pro-nilpotent functions}

\begin{notation}
  We denote by $\fAffLine \in \dInfSch_\basefield$ the formal completion of $\affineline{\basefield}$ at $0$;
\end{notation}

\begin{observation}
  The étale map $\fAffLine \to \affineline{}$ induces an equivalence of derived abelian group schemes $\Omega_0 \fAffLine \simeq \bbB$ (with $\bbB$ as in \cref{notationforB}).
\end{observation}

\begin{lemmalist}{\label{lemmapronilpotent}
    Let $Y \in \dInfSch_\basefield$ and $f \colon Y \to \affineline{}$ a function.
    The following $\infty$\nobreakdashes-groupoids are equivalent and are either empty or contractible:}
  \item The $\infty$\nobreakdashes-groupoid of nullhomotopies of $Y \to \affineline{} \to \affineline{\deRham}$.
  \item The $\infty$\nobreakdashes-groupoid of nullhomotopies of $Y_\red \to Y \to \affineline{}$.
  \item The $\infty$\nobreakdashes-groupoid of factorizations of $Y \to \affineline{}$ through $\smash{\fAffLine} \subset \affineline{}$.
\end{lemmalist}
\begin{proof}
  The $\infty$\nobreakdashes-groupoids (a) and (b) are equivalent by adjunction, and the equivalence with (c) follows from the fiber sequence $\fAffLine \to \affineline{} \to \affineline{\deRham}$ (see for instance \cite[Construction 2.3.1]{bhatt2022prismatic}).
  Since $Y_\red$ is a scheme by assumption, the $\infty$\nobreakdashes-groupoid (b) is either empty or a single point and the result follows.
\end{proof}

\begin{definition}\label{definitionpronilpotent}
  In the situation of \cref{lemmapronilpotent}, we will say that the function $f$ (or the inf-schematic $\LG$\nobreakdashes-pair $(Y,f)$) is \emph{pro-nilpotent} if the equivalent $\infty$\nobreakdashes-groupoids of \cref{lemmapronilpotent} are non-empty.
  The full subcategory of $\LGdinfzero$ spanned by pro-nilpotent functions is canonically identified with the slice category $\dInfSch_{\fAffLine}$.
\end{definition}

\begin{remark}
  If $\formalU \to \fAffLine$ is a pro-nilpotent function on $\formalU$ a derived inf-scheme, then its derived zero fiber $\zerolocus(\formalU, f)$ together with the action of the group $\bbB$ of \cref{definitionoffunctorPreMFoninfschemes} completely determines the function $f$ via the equivalences $\formalU \simeq \quot{\zerolocus(\formalU,f)}{\bbB} \to \quot{\ast}{\bbB} \simeq \fAffLine$ (in suitably considered formal moduli problems).
\end{remark}

\begin{example}\label{pronilpotentexample}
  If $U$ is a smooth $\basefield$-scheme, any function $f \colon U \to \affineline{}$ such that $\Crit(f)_\red \subset (f^{-1}(0))_\red$ (\ie $f_{|_{\Crit(f)_\red}}=0$) induces, by passing to formal completions both on $U$ and in $\affineline{}$, a pro-nilpotent function on the formal neighborhood of $\Crit(f)$ in $U$, $\widehat{f} \colon \formalU \coloneqq \widehat{U}^{\Crit(f)} \to \fAffLine$. Both ends of $\widehat{f}$ are inf-schemes by \cref{definitioninfschemes}.
\end{example}

The following lemma is due to \cite{MR2735755}.
\begin{lemma}\label{lem:mfcompletion}
  Let $(U,f) \in \LGdinfzero$ be an $\LG$-pair with $U$ a smooth scheme.
  Denote by $f \colon \widehat U \to \fAffLine$ the completion around $\Crit(f)$.
  Then restriction along the inclusion $\widehat U \to U$ induces equivalences
  \[
    \MF^{\infty}(U,f) \to^\sim \MF^{\infty}(\widehat U, f)
    \textrm{ and }
    \MF(U,f) \to^\sim \MF(\widehat U, f).
  \]
\end{lemma}
\begin{proof}
  Consider the derived fiber products $U_0 = U \times_{\affineline{}} \{0\}$ and $\widehat U_0 = \widehat U \times_{\fAffLine} \{0\}$ and $i\colon \widehat U_0 \to U_0$ the inclusion.
  Notice that $\widehat U_0$ is also the completion of $U_0$ in $U$, as completion (\cf \cref{derhamandformalcompletion}) commutes with finite limits. Now, consider the composition
  \[
    \IndCoh_{\Crit(f)}(U_0) \to \IndCoh(U_0) \to \IndCoh(\widehat U_0)
  \]
  where the first map is the canonical inclusion of Ind-coherent sheaves with support (see for instance \cite[\S 7.4.4]{gaitsgoryrozenblyum_dgindschemes}) and the second map is the $!$-pullback along the inclusion of the formal completion $\widehat U_0 \to U_0$.
  The inclusion functor $\IndCoh_{\Crit(f)}(U_0) \hookrightarrow \IndCoh(U_0)$ is the kernel of the $\IndCoh(\bbB)$-linear functor $\IndCoh(U_0) \to \IndCoh((U \smallsetminus \Crit(f))_0)$ and is thus $\IndCoh(\bbB)$-linear itself.
  The restriction functor $\IndCoh(U_0) \to \IndCoh(\widehat U_0)$ is also $\IndCoh(\bbB)$-linear.
  We thus get functors
  \[
    \MF^\infty_{\Crit(f)}(U,f) \coloneqq \IndCoh_{\Crit(f)}(U_0) \otimes_\kbeta \kbetaLaurent \to \MF^\infty(U,f) \to \MF^\infty(\widehat U, f).
  \]
  By \cite[Prop.\,4.1.6]{MR3121870} the first map is an equivalence.
  The composition is an equivalence as $\IndCoh_{\Crit(f)}(U_0) \to \IndCoh(\widehat U_0)$ is (\cf \cite[Prop.\,7.4.5]{gaitsgoryrozenblyum_dgindschemes} -- $\Crit(f)$ is by construction defined by finitely many polynomials).
  Thus $\MF^\infty(U,f) \to \MF^\infty(\widehat U, f)$ is an equivalence as announced.
\end{proof}

In view of the above lemma, we will restrict ourselves to pro-nilpotent functions throughout the rest of the paper.

\subsection{Convolution monoidal structure and Thom--Sebastiani on \texorpdfstring{$\MF^\infty$}{MF∞} over \texorpdfstring{$\basefield$}{ℂ}}\label{sectionconvolutionproductsaxioms}

In this section we construct a symmetric monoidal enhancement of the functor $\MF^\infty$ of \cref{definitionoffunctorMFfunctoroninfschemes} -- \cf \cref{thomsebastianiMF}.
In order to achieve this we need to describe first a general construction of a convolution symmetric monoidal structure.

\begin{construction}\label{convolutionmonoidalstructure}
  Let $\calC$ be an \icategory with finite products. We construct a lax symmetric monoidal structure on the \ifunctor
  \[
    \calC_{/-} \colon \calC \to \inftycats,\hspace{1cm}X\mapsto \calC_{/X},\hspace{1cm} f \colon X \to Y\mapsto f\circ - \colon \calC_{/X} \to \calC_{/Y}
  \]
  given on a pair of objects $Y_1,Y_2 \in \calC$, by the product functor
  \[
    \calC_{/Y_1} \times \calC_{/Y_2} \to \calC_{/Y_1 \times Y_2} \hspace{1.2cm} (X_1 \to Y_1, X_2 \to Y_2) \mapsto (X_1 \times X_2 \to Y_1 \times Y_2).
  \]

  To exhibit this symmetric lax monoidal structure in a more precise way, one must provide the map of \ioperads
  \[
    \begin{tikzcd}[row sep=small, column sep=small]
      \calC^{\times} \ar{dr}[swap]{q} \ar{rr} && \inftycats^{\times} \ar{dl} \\
      & \pFSets &
    \end{tikzcd}
  \]
  implementing it, where $\calC^{\times}$ is as in \cite[Const.\,2.4.1.4]{lurie-ha}. By \cite[2.4.1.7]{lurie-ha}, the data of such a map of \ioperads corresponds to the data of a \emph{lax Cartesian structure}
  \begin{equation}
    \label{laxcartesianstructure}
    \calC^{\times} \to \inftycats
  \end{equation}
  in the sense of \cite[2.4.1.1]{lurie-ha}. It is therefore enough to construct the corresponding coCartesian fibration classifying \cref{laxcartesianstructure} and to the verify it satisfies the property of a lax Cartesian structure.
  Since the evaluation map $\ev_1 \colon \Fun(\Delta^1,\calC) \to \calC$ commutes with Cartesian products, by \cite[2.4.1.8]{lurie-ha} it extends to a map of \ioperads
  \[
    \begin{tikzcd}[row sep=small, column sep=small]
      \Fun(\Delta^1,\calC)^{\times} \ar{dr}[swap]{p} \ar{rr}{\ev_1^{\times}} && \calC^{\times} \ar{dl}{q} \\
      & \pFSets\rlap{.} &
    \end{tikzcd}
  \]
  We claim that $\ev_1^{\times}$ is in fact a coCartesian fibration, classifying the required map \cref{laxcartesianstructure}. Indeed, by definition (\cf \cite[Construction 2.4.1.4]{lurie-ha}), we have a Cartesian square of \icategories
  \[
    \begin{tikzcd}[row sep=small, column sep=small]
      \Fun(\Delta^1,\calC)^{\times} \ar{d}\ar[hookrightarrow]{r}\tikzcart &\ar{d} \Fun(\Delta^1, \calC^{\times})\\
      \pFSets\ar[hookrightarrow]{r}{\Id}& \Fun(\Delta^1, \pFSets)
    \end{tikzcd}
  \]
  and $\ev_1^{\times}$ identifies with the composition with the evaluation map at $1$, $e$
  \begin{equation}\label{compositioncocartesian}
    \begin{tikzcd}
      \Fun(\Delta^1,\calC)^{\times} \ar[hookrightarrow]{r}{} & \Fun(\Delta^1, \calC^{\times})\ar{r}{e} & \calC^{\times}
    \end{tikzcd}
  \end{equation}
  By \cite[2.4.7.11]{lurie-htt} $e$ is a coCartesian fibration, and therefore to show that the composition \cref{compositioncocartesian} is coCartesian we can conclude using \cite[\href{https://kerodon.net/tag/02R5}{Prop.\,02R5}]{kerodon}. This is a routine computation using the fact $q$ is a coCartesian fibration and the explicit definition of coCartesian edges in \cite[2.4.1.10-(2)]{lurie-htt}.
\end{construction}

\begin{constructions}{ Let $\calC$ be an \icategory with Cartesian products.}

  \item \label{convolutionmonoidalstructurealgebra}
  Let $A$ be a commutative algebra object in $\calC^{\times}$ with multiplication $A \times A \to A$. Then via the symmetric lax monoidal structure on the functor $\calC_{/-} \colon \calC \to \inftycats$ of \cref{convolutionmonoidalstructure}, the \icategory $\calC_{/A}$ inherits a symmetric monoidal structure given informally by the functor
  \[
    \begin{tikzcd}
      \calC_{/A} \times \calC_{/A} \ar{r}{\times} & \calC_{/A \times A} \ar{r}{m} & \calC_{/A}
    \end{tikzcd}
  \]
  with tensor unit $0 \colon \Spec(\basefield) \to A$. The corresponding coCartesian fibration $\calC_{/A}^{\boxplus} \to \pFSets$ is the fiber product
  \[
    \begin{tikzcd}[row sep=small, column sep=small]
      \calC_{/A}^{\boxplus} \ar{r} \ar{d} & \Fun(\Delta^1,\calC)^{\times} \ar{d}{\cref{compositioncocartesian}}
    \\
      \pFSets\ar{r}{A} & \calC^{\times}
    \end{tikzcd}
  \]
  where $A$ is the section of $\calC^{\times} \to \pFSets$ classifying the commutative algebra object given by $A$. In \crossref{constructionLiouville4marsmonoidal} we already described a particular case of this construction.
  This is a particular instance of \cite[2.2.2.4, 2.2.2.5]{lurie-ha} in the case the monoidal structure is Cartesian.

  \item \label{pullbackmapalgebrasconvolution} If $\phi \colon A \to B$ is a map of commutative algebra objects in $\calC^{\times}$ we can form the Cartesian square
  \[
    \begin{tikzcd}
      \mathcal{P}\ar{r}\ar{d}{\tilde{p}}& \Fun(\Delta^1,\calC)^{\times}\ar{d}{\cref{compositioncocartesian}}\\
      \Delta^1\times \pFSets\ar{r}{\phi}&\calC^{\times}
    \end{tikzcd}
  \]
  to obtain a coCartesian fibration $\tilde{p}$ classifying a symmetric monoidal functor $\phi_\ast \colon \calC_{/A}^{\boxplus} \to \calC_{/B}^{\boxplus}$ induced by composition with $\phi$. Since $\calC$ admits finite limits, $\phi_*$ admits a right adjoint $\phi^\ast$ taking fiber products. By \cite[Cor.\,7.3.2.7]{lurie-ha}, the \ifunctor $\phi^\ast$ is \emph{lax} symmetric monoidal.
\end{constructions}

\begin{observation}\label{laxmonoidalunitfactorizationmodules}
  Let $F \colon \calC^{\otimes} \to \calD^{\otimes}$ be a symmetric monoidal functor between presentably symmetric monoidal \icategories\footnote{In particular, compatible with geometric realizations as in \cref{relativetensor}.} (with units $\unit_\calC$ and $\unit_\calD$) with a lax monoidal right adjoint $G$ and set $B \coloneqq G(\unit_\calD)$. Then, since $F$ commutes with colimits, base change along the counit of the adjunction $F(B)=F(G(\unit_\calD)) \to \unit_{\calD}$ induces a commutative diagram of symmetric monoidal left adjoints
  \[
    \begin{tikzcd}[every label/.append style = {font = \tiny}]
      \calC^{\otimes}\ar{r}{F}\ar{d}[swap]{-\otimes B} & \calD^{\otimes}\ar{dr}[sloped]{\sim} \ar{d}[swap]{-\otimes F(B)}
      \\
      \Mod_{B}(\calC)^{\otimes_B}\ar{r}{F} & \Mod_{F(B)}(\calD)^{\otimes_{F(B)}} \ar{r}[swap]{_{- \otimes_{F(B)}\unit_\calD}} &[2em] \Mod_{\unit_\calD}(\calD)^{\otimes}.
    \end{tikzcd}
  \]
  Passing to right adjoints and evoking \cite[Cor.\,7.3.2.7]{lurie-ha} we obtain a lax symmetric monoidal factorization $\overline{G}$ of $G$
  \[
    \begin{tikzcd}
      \calD^{\otimes}\ar[bend left=15]{rr}{G} \ar{r}[swap]{\overline{G}} & \Mod_{B}(\calC)^{\otimes_B}\ar{r}[swap]{\mathsf{forget}} & \calC^{\otimes}\rlap{.}
    \end{tikzcd}
  \]
\end{observation}

We are now finally ready to start the construction of the lax monoidal structure on the functor $\MF^\infty$.

\begin{construction}\label{zerofiberfunctorlaxmonoidal}
  Let $\calC=\dSt_\basefield$ and consider the abelian group object $A=(\fAffLine,+)$.
  By \cite[1.2.13.8]{lurie-htt} the forgetful functor $\dSt_{/\fAffLine} \to \dSt$ creates colimits and $\dSt_{/\fAffLine}$ is presentable by \cite[5.5.3.10]{lurie-htt}. By \cite[3.4.4.1, 4.8.1.9]{lurie-ha} the symmetric monoidal structure $\dSt_{/\fAffLine}^{\boxplus}$ of \cref{convolutionmonoidalstructurealgebra} is a presentable symmetric monoidal \icategory.
  As the inclusion $0 \colon \Spec(\basefield)=\ast \to \smash{\fAffLine}$ is the unit of the abelian group, in particular it defines a map of commutative algebras in $\dSt^{\times}$ and by \cref{pullbackmapalgebrasconvolution} pullback along $0$ defines a \emph{lax} symmetric monoidal functor encoded by a map of \ioperads
  \[
    0^\ast\colon\dSt_{/\fAffLine}^{\boxplus} \to \dSt_{/\ast}^{\boxplus} \simeq \dSt^{\times}
  \]
  sending the tensor unit $0 \colon \ast \to \fAffLine$ to the derived abelian group scheme $\bbB \coloneqq \Omega_0\fAffLine$ of \cref{notationforB}. \cref{laxmonoidalunitfactorizationmodules} applied to the adjunction $(0_\ast, 0^\ast)$ yields a \emph{lax} symmetric monoidal factorization $\zerolocus^{\boxplus} \coloneqq \overline{0^\ast}$
  \[
    \begin{tikzcd}
      \dSt_{/\fAffLine}^{\boxplus}\ar[bend left=15]{rr}{0^\ast}\ar{r}[swap]{\zerolocus^{\boxplus}}
      & \Mod_{\bbB}(\dSt)^{\otimes_{\bbB}}\ar{r}[swap]{\mathsf{forget}}
      & \dSt^{\times}
    \end{tikzcd}
  \]
  where $\Mod_{\bbB}(\dSt)^{\otimes_{\bbB}}$ is equipped with the relative tensor product $\otimes_{\bbB}$ given by the relative bar construction, \cf \cref{relativetensor}.
\end{construction}

We can finally present the lax monoidal structure on $\MF^\infty$:

\begin{construction}\label{constructionMFasconvolutionmonoidal}
  We proceed in several steps:

  \begin{enumerate}[leftmargin=1em,itemindent=1em,label={(\arabic*)}]
    \item As in \cref{convolutionmonoidalstructurealgebra} the abelian group structure on $\fAffLine$ induces a symmetric monoidal structure $\smboxplus$ on the slice category $\dInfSch_{\fAffLine}$ (\cf \cref{definitionpronilpotent}), informally given by $(\formalU, f) \smboxplus (\calV,g) \coloneqq (\formalU \times \calV, f + g)$. \cref{dinfschemesstableunderfinitelimite} implies $\dInfSch_{\fAffLine}^{\boxplus}\subseteq \dSt_{/\fAffLine}^{\boxplus}$ is a symmetric monoidal full sub-$\infty$-category.

    \item \label{relativeBproducts} Since $\bbB\in  \dInfSch_\basefield$, combining \cref{zerofiberfunctorlaxmonoidal} and \cref{dinfschemesstableunderfinitelimite} we conclude that the $0$\nobreakdashes-locus functor $\zerolocus$ (\cf \cref{definitionoffunctorPreMFoninfschemes}) yields a morphism of \ioperads
          \[
            \zerolocus^{\boxplus} \colon \dInfSch_{\fAffLine}^{\boxplus} \to \Mod_{\bbB}(\dInfSch_\basefield)^{\otimes_{\bbB}} \subset \Mod_{\bbB}(\dSt)^{\otimes_{\bbB}}.
          \]
    \item The symmetric monoidal functor $\hrm$ of \cref{definitionoffunctorPreMFoninfschemes} sends the commutative algebra object $\bbB$ to a commutative algebra object $\hrm(\bbB)$ in correspondences and by passing to module objects we obtain a morphism of \ioperads \footnote{Notice that as the \icategory of correspondences does not have colimits in general: as explained in \cref{relativetensor} we can only guarantee the structure of \ioperad on module objects.}
          \[
            \hrm \colon \Mod_{\bbB}(\dInfSch_\basefield)^{\otimes_{\bbB}} \to \Mod_{\hrm(\bbB)}(\Corr(\dInfSch_\basefield))^{\otimes_{\hrm(\bbB)}}.
          \]

    \item Consider the composite morphism of \ioperads
          \[
            \hrm\circ\, \zerolocus\colon \dInfSch_{\fAffLine}^{\boxplus} \to \Mod_{\hrm(\bbB)}(\Corr(\dInfSch_\basefield))^{\otimes_{\hrm(\bbB)}}.
          \]
          We saw in \cref{definitionoffunctorPreMFoninfschemes} that the underlying \ifunctor factors through correspondences. The \emph{lax monoidal} universal property of correspondences\footnote{The construction of a symmetric monoidal structure on the category of correspondences $\Corr(\calC)$ of \cite[Part III, Ch.\,9 , \S 2]{GR-I} when $\calC^{\otimes}$ is not Cartesian seems to lack a supplementary hypothesis, namely, the tensor product functor $\otimes \colon \calC \times \calC \to \calC$ must preserve fiber products (as a whole, not in each variable separately).
            This is because the construction sending $\calC\mapsto \Corr(\calC)$ is only functorial with respect to functors preserving fiber products.
            In our case, the category $(\dInfSch_{\fAffLine})^{\boxplus}$ satisfies this requirement since $\boxplus$ is given by a fiber product and therefore $\Corr(\dInfSch_{\fAffLine})$ inherits a tensor product $\boxplus$.} of
          \cite[Part III, Ch.\,9 Prop.\,3.2.4]{GR-I} provides an enhancement of this factorization into a morphism of \ioperads
          \[
            \Corr(\dInfSch_{\fAffLine})^{\boxplus} \to \Mod_{\hrm(\bbB)}(\Corr(\dInfSch_\basefield))^{\otimes_{\hrm(\bbB)}}.
          \]
    \item The symmetric monoidal functor $\IndCoh$ of \cref{IndCohassymmetricmonoidalfunctoronindinfschemes}, induces, by passing to module objects, a morphism of \ioperads (using \cref{ICBiskbeta})
          \[
            \IndCoh\colon \Mod_{\hrm(\bbB)}(\Corr(\dInfSch_\basefield))^{\otimes_{\hrm(\bbB)}} \to \Mod_{\IndCoh(\bbB)}\left(\PrLSt[\basefield]\right)^{\otimes_{\IndCoh(\bbB)}} \simeq \PrLStmonoidal[\kbeta]
          \]
          (where the target is in fact a symmetric monoidal \icategory).
  \end{enumerate}
  Finally, composing the operadic morphisms obtained above, we obtain a \emph{lax} symmetric monoidal enhancement of the functor $\PreMF$ of \cref{definitionoffunctorPreMFoninfschemes},
  \[
    \PreMFmon \colon \Corr\left(\dInfSch_{\fAffLine}\right)^{\boxplus} \to \PrLStmonoidal[\kbeta].
  \]
  It in turn induces a \emph{lax} monoidal structure on the \ifunctor $\MF^\infty$ of \cref{definitionoffunctorMFfunctoroninfschemes} implementing matrix factorizations by inverting $\beta$:
  \[
    \MF^{\infty, \boxplus} \colon
    \begin{tikzcd}
      \Corr\left(\dInfSch_{\fAffLine}\right)^{\boxplus} \ar{r}{ \PreMFmon} & \PrLStmonoidal[\kbeta] \ar{rr}{- \otimes_\kbeta \kbetaLaurent} && \PrLStmonoidal[\kbetaLaurent]\rlap{.}
    \end{tikzcd}
  \]
\end{construction}

\begin{observation}\label{observationMFzeroisICsymmetricmonoidal}
  The symmetric monoidal functor $\dInfSch^{\times}_\basefield \to \dInfSch_{\fAffLine}^{\boxplus}$ sending $S$ to $(S,0)$ induces a symmetric monoidal functor at the level of correspondences.
  By construction and uniqueness of the extension to correspondences, we get a commutative square of lax symmetric monoidal functors
  \[
    \begin{tikzcd}[row sep=small, column sep=small]
      \Corr(\dInfSch_\basefield)^{\times}\ar{r}\ar{d}{\IndCoh}& \Corr( \dInfSch_{\fAffLine})^{\boxplus}\ar{d}{\PreMF}\\
      \PrLSt[\basefield]\ar{r}{\otimes\ofbeta}& \PrLSt[\ofbeta]\rlap{.}
    \end{tikzcd}
  \]
  This upgrades the equivalence $\MF^\infty(S,0) \simeq \IndCoh(S)\ofbetaLaurent$ of \cref{MFzerofunction} to a symmetric monoidal equivalence.
\end{observation}

Finally, \cref{constructionMFasconvolutionmonoidal} allows us to extend the Thom--Sebastiani equivalence of Preygel (see \cite[Thm.\,4.1.3]{MR3121870}) to general derived inf-schemes:
\begin{theorem}[Thom--Sebastiani equivalences]\label{thomsebastianiMF}
  The \emph{lax} symmetric monoidal functors $\PreMFmon$ and $\MF^{\infty,\boxplus}$ obtained in \cref{constructionMFasconvolutionmonoidal}
  are \emph{strongly} symmetric monoidal.
\end{theorem}

\begin{proof}
  The base change functor $\PrLSt[\kbeta] \to \PrLSt[\kbetaLaurent]$ is symmetric monoidal.
  It therefore suffices to focus on the case of $\PreMFmon$.
  Fix $f \colon \formalU \to \fAffLine$ and $g \colon \calV \to \fAffLine$ and set $\formalU_0 \coloneqq \zerolocus(\formalU,f)$, $\calV_0 \coloneqq \zerolocus(\calV,g)$ and $(\formalU \times \calV)_0 \coloneqq \zerolocus(\formalU \times \calV,f + g)$.
  The comparison functor arising from the given lax monoidal structure
  \[
    \PreMF(\formalU, f) \otimes_{\kbeta} \PreMF(\calV, g) \to \PreMF\left(\formalU \times \calV, f + g\right)
  \]
  is induced from the $\bbB$\nobreakdashes-bi-equivariant closed immersion
  $i \colon \formalU_0 \times \calV_0 \to \left(\formalU \times \calV\right)_0$.
  Since $i$ can be obtained as the inclusion of the zero locus of $f-g \colon \left(\formalU \times \calV\right)_0 \to \fAffLine$, it is schematic, proper and surjective on $\basefield$\nobreakdashes-points.
  Thus the conditions of \cite[I-3 Prop.\,3.3.3 (b)]{GR-II} are satisfied and we find that $\PreMF(\formalU \times \calV, f+g) = \IndCoh(\left(\formalU \times \calV\right)_0)$ is equivalent to the totalization of $\IndCoh_*(\Nerve(i))$ of $\PreMF$ of the \v Cech nerve of $i$ (with push-forward functors).
  The bi-equivariance of $i$ then identifies this simplicial diagram with the bar complex computing $\PreMF(\formalU, f) \otimes_{\kbeta} \PreMF(\calV, g)$.
\end{proof}

\cref{thomsebastianiMF} has several consequences worth mentioning here:
\begin{corollarylist}{\label{MFinftydualizable}}
  \item For any $(\formalU, f) \in \dInfSch_{\fAffLine}$, the \icategory $\MF^\infty(\formalU,f)$ is dualizable in $\PrLSt[\kbetaLaurent]$, with dual $\MF^\infty(\formalU, -f)$.
  \item For any $(\formalU, f), (\calV, g) \in \dInfSch_{\fAffLine}$, we have
  \[
    \myuline\Fun_{\kbetaLaurent}\left(\MF^\infty(\formalU,f), \MF^\infty(\calV,g)\right) \simeq \MF^\infty\left(\formalU \times \calV,-f \smboxplus g\right).
  \]
\end{corollarylist}
\begin{proof}
  In the category $\Corr\left(\dInfSch_{\fAffLine}^{\boxplus}\right)$, every object $(\formalU, f)$ is dualizable, with dual $(\formalU, -f)$. 
  Evaluation and coevaluation are both given by the correspondence
  \[
    \begin{tikzcd}
      \mathllap{\mathbb{1} ={}} (*,0) & (\formalU, 0) \ar{l} \ar{r}{\Delta} & (\formalU,f) \smboxplus (\formalU,-f).
    \end{tikzcd}\qedhere
  \]
\end{proof}

\begin{warning}
  Fix $(\formalU, f) \in \dInfSch_{\fAffLine}$.
  The derived schemes $\zerolocus(f)$ and $\zerolocus(-f)$ are obviously equivalent, but their natural actions of the group scheme $\bbB$ differ by $-1$.
  In particular, we get
  \[
    \MF^\infty(\formalU, f) \otimes_{\kbetaLaurent} \kbetaLaurent \simeq \MF^\infty(\formalU, -f),
  \]
  where $\kbetaLaurent \to^\sim \kbetaLaurent$ maps $\beta$ to $-\beta$.
  As a $\basefield$-linear \icategory, we thus always have $\MF^\infty(\formalU, -f) \simeq \MF^\infty(\formalU, f)$. The $2$-periodic structures however differ by $-1$.
\end{warning}

Using \cref{rmkfunctorialitiesMF}, \cref{thomsebastianiMF} translates into a similar statement for small categories of matrix factorizations.

\begin{corollary}
  Denote by $\Corr(\dInfSch_{\fAffLine}^{\boxplus})^{\mathrm{ec,pr}} \subset \Corr(\dInfSch_{\fAffLine}^{\boxplus})$ the non-full sub-\icategory containing all the objects but as morphisms only those correspondences
  $
    \smash{X \from^p Y \to^q Z}
  $
  where $p$ is schematic and eventually coconnective and $q$ is proper.
  The construction $\MF$ of \cref{smallMF} extends to a symmetric monoidal \ifunctor
  \[
    \Corr(\dInfSch_{\fAffLine}^{\boxplus})^{\mathrm{ec,pr}} \to \mathsf{Cat}^{\mathrm{St,ic}}_{\infty,\kbetaLaurent}
  \]
  where $\mathsf{Cat}^{\mathrm{St,ic}}_{\infty,\kbetaLaurent}$ denotes the \icategory of small idempotent complete $\kbetaLaurent$\nobreakdashes-linear stable \icategories.\qed
\end{corollary}

\subsection{Interlude: relative convolution}\label{sectionrelativeconvolution}

In order to obtain more structure on the categories of matrix factorizations, we will need further $\infty$-categorical considerations.
The trusting reader can safely skip this subsection.

\subsubsection{Relative convolution and modules}

\begin{observations}{\label{observationfiberproductaslimit}
    Let $\calC^{\amalg}$ be a coCartesian symmetric monoidal \icategory.}

  \item \label{uniquecoalgebracocartesian} By \cite[2.4.3.10]{lurie-ha} the forgetful functor $\CAlg(\calC) \to \calC$ is an equivalence. Its inverse equips an object $X \in \calC$ with the commutative algebra structure given by the canonical map $X \amalg X \to X$.

  \item \label{descriptionofmodulescocartesian} Let $X \in \calC \simeq \CAlg(\calC)$. As an \icategory, $\Mod_X(\calC)$ is equivalent to the slice category $\calC_{X/}$ via the forgetful functor. This equivalence is functorial in both $\calC$ and $X$. See \cite[Lem.\,2.1]{zbMATH07681894}, \cite[Cor.\,2.6.6, Rmk.\,2.6.7]{zbMATH07506911}.
\end{observations}

\begin{lemma}\label{modulesiscocartesian}
  Let $\calC$ be an \icategory with finite colimits and geometric realizations and $\calD \subset \calC$ a full subcategory stable under finite colimits.
  Denote by $\calD^{\amalg} \subset \calC^{\amalg}$ the associated coCartesian symmetric monoidal structures.
  Then, for any $X \in \calD \simeq \CAlg(\calD)$, the \ioperad $p \colon \Mod_X(\calD)^{\otimes_X} \to \pFSets$ of \cref{relativetensor} is again a coCartesian symmetric monoidal \icategory, equivalent to $\calD_{X/}^{\amalg}$.
  \begin{proof}
    Since the inclusion $\calD \to \calC$ preserves finite colimits, it induces compatible symmetric monoidal fully faithful functors $\calD^{\amalg} \hookrightarrow \calC^{\amalg}$ and $\calD_{X/}^{\amalg} \hookrightarrow \calC_{X/}^{\amalg}$.
    In particular, we can reduce to the case $\calD = \calC$.

    As recalled in \cref{relativetensor}, under the current assumptions, the \ioperad $\Mod_X(\calC)^{\otimes_X}$ is a symmetric monoidal \icategory with tensor products given by the relative bar construction.
    More concretely, using the equivalence $\Mod_X(\calC) \simeq \calC_{X/}$ of \cref{descriptionofmodulescocartesian}, if $X \to Z$ and $X \to Y$ are two $X$-modules in $\calC$, according to the formula in \cite[4.4.2.8]{lurie-ha} we have
    \[
      Z \otimes_X Y \simeq \colim_{\Delta^\op} \left(Z \amalg X^{\amalg n} \amalg Y\right).
    \]
    We now observe that this colimit can be expressed as a pushout.
    Let $\horn_2^0$ be the category indexing diagrams of the form $1 \leftarrow 0 \to 2$ and consider the functor $\Psi \colon \horn_2^0 \to \Delta^{\op}$ sending $0 \mapsto [1]$, both $1$ and $2$ to $[0]$ and the two arrows to the two boundary maps.
    Denote by $d \colon \horn_2^0 \to \calC$ the functor classifying the diagram $Z \leftarrow X \to Y$ in $\calC$, the left Kan extension $\LKE(d)$ along $\Psi$ is the simplicial diagram encoding the bar construction.
    Therefore, by definition of left Kan extensions, we have an equivalence
    \[
      Z \amalg_X Y \simeq \colim_{\Delta^\op} \left(Z \amalg X^{\amalg n} \amalg Y\right) \simeq Z \otimes_X Y.
    \]
    In particular, this upgrades the equivalence of \cref{descriptionofmodulescocartesian} to a symmetric monoidal equivalence
    \[
      \Mod_X(\calC)^{\otimes_X} \simeq \calC_{X/}^{\amalg}.\qedhere
    \]
  \end{proof}
\end{lemma}

\begin{reminders}{
    Let $\calC^{\otimes}$ be a symmetric monoidal \icategory.}
  \item Recall that a \emph{commutative bialgebra} object in $\calC$ is an object of
  \[
    \BiCAlg(\calC) \coloneqq \CAlg(\CAlg(\calC)^{\op})^{\op} \simeq \CAlg(\CAlg(\calC^\op)^\op).
  \]
  \item \label{bialgebramonoidalstructure} The box product $\boxtimes$ endows the \ifunctor $\CAlg(\calC)^\op \to \inftycats$ mapping $A$ to $\Mod_A(\calC)$ (with forgetful functors) with a lax symmetric monoidal structure (follows from \cite[Thm.\,4.8.5.16]{lurie-ha} in the presentable case -- for the general case, see \cite[Const.\,2.2.6]{Arpon}).
  In particular, modules over commutative bialgebras come with a symmetric monoidal structure:
  \[
    \begin{tikzcd}[column sep=large]
      \BiCAlg(\calC)^\op \simeq \CAlg(\CAlg(\calC)^{\op}) \ar{r}{\Mod_{-}(\calC)^\odot} & \CAlg(\inftycats) \eqqcolon \monoidalcats.
    \end{tikzcd}
  \]
  For $A \in \BiCAlg(\calC)$, we denote by $\odot$ the resulting symmetric monoidal structure on $\Mod_A(\calC)$.
  Note that the monoidal unit $\unit \in \calC$ has a unique bi-algebra structure and is both initial and final in $\BiCAlg(\calC)$.
  For any $A \in \BiCAlg(\calC)$, we thus have symmetric monoidal \emph{trivial module} and \emph{forgetful functors}
  \[
    \begin{tikzcd}
      \calC^{\otimes} \ar{r}{\mathsf{Triv}} & \Mod_A(\calC)^\odot \ar{r}{\mathsf{forget}} & \calC^{\otimes}.
    \end{tikzcd}
  \]
  \item \label{functorialitybialgebratensor}
  Let $F \colon \calC^{\otimes} \to \calD^{\otimes}$ be a symmetric monoidal functor and $A \in \BiCAlg(\calC)$.
  The induced functor $F \colon \Mod_A(\calC)^{\odot} \to \Mod_{F(A)}(\calD)^{\odot}$ is symmetric monoidal (see \cite[Const.\,2.2.6]{Arpon}).
\end{reminders}

\begin{notation}
  Let $\calC^{\amalg}$ be a coCartesian symmetric monoidal category and $A \in \CoCAlg(\calC) \coloneqq \CAlg(\calC^\op)^\op$.
  We denote by $\boxdot$ the symmetric monoidal structure on $\calC_{A/}$ defined as the opposite of the convolution on $(\calC^\op)^{\boxplus}_{/A}$ of \cref{convolutionmonoidalstructurealgebra}.
  We abusively also call $\boxdot$ the convolution product.
\end{notation}

Finally, we can give an alternative description of the convolution product:
\begin{lemma}\label{alternativedescriptionconvolutionascomodules}
  Let $\calC^{\amalg}$ be a coCartesian symmetric monoidal category.
  \Cref{uniquecoalgebracocartesian} induces an equivalence $\CoCAlg(\calC) \simeq \BiCAlg(\calC)$
  and the equivalence of \cref{descriptionofmodulescocartesian} lifts to an symmetric monoidal equivalence between the two functors
  \[
    \begin{tikzcd}[row sep=-.4em]
      \CoCAlg(\calC)^\op \ar[phantom]{r}{\simeq} &[-2em] \BiCAlg(\calC)^\op \ar{r} & \monoidalcats \\
      & A \ar[mapsto]{r} & \Mod_A(\calC)^\odot \\
      A \ar[mapsto]{rr} && \calC_{A/}^{\boxdot}.
    \end{tikzcd}
  \]
  In particular, for any $A \in \CoCAlg(\calC)$, we have a symmetric monoidal equivalence
  \[
    \Mod_{A}(\calC)^{\odot} \simeq \calC_{A/}^{\boxdot}.
  \]
  \begin{proof}
    By construction, the convolution product arises from the lax monoidal structure on the \ifunctor $\calC_{-/} \colon (\calC^{\amalg})^\op \to \inftycats^{\times}$ sending $X \mapsto \calC_{X/}$, described in \cref{convolutionmonoidalstructure}.
    Similarly, the tensor structure $\odot$ is induced by the lax monoidal structure on $\Mod_{-}(\calC) \colon \CAlg(\calC)^\op \to \inftycats$.
    Since $\calC$ is equipped with the coCartesian monoidal structure, the equivalence of \cref{descriptionofmodulescocartesian} can be promoted to a natural equivalence rendering the commutativity of the diagram
    \[
      \begin{tikzcd}[column sep=0, row sep=scriptsize]
        \calC^\op \ar{rr}[swap]{\sim} \ar{dr}[swap]{\calC_{-/}} && \ar{dl}{\Mod} \CAlg(\calC)^\op
        \\
        & \phantom{{}_\infty}\inftycats\rlap{.} &
      \end{tikzcd}
    \]
    The universal property of the lax structure provided by \cref{convolutionmonoidalstructure} allows to identify the two lax monoidal structures at hand.
    We omit the details.
  \end{proof}
\end{lemma}

\cref{alternativedescriptionconvolutionascomodules} allows us to treat convolution products as relative tensor products of comodules.
We will use this fact to construct relative versions of the convolution product.
\begin{construction}\label{coalgebraonSzero}
  Let $\calC$ be an \icategory with finite colimits and let $A \in \CoCAlg(\calC^{\amalg})$ with counit $0 \colon A \to \unit$ (where $\unit$ is the monoidal unit, a.k.a.\,the initial object of $\calC$).
  It follows from \cref{pullbackmapalgebrasconvolution} that the \ifunctor $0^\ast \colon \calC = \calC_{\unit/} \to \calC_{A/}$ sending $Y \in \calC$ to $A \to \unit \to Y$ defines a symmetric monoidal functor with respect to convolution $0^\ast \colon \calC^{\amalg} = \calC_{\unit/}^{\boxdot} \to \calC_{A/}^{\boxdot}$.
  Using \cref{uniquecoalgebracocartesian}, we get
  \[
    0^\ast \colon \calC \simeq \CAlg(\calC^{\amalg}) \to \CAlg(\calC_{A/}^{\boxdot}).
  \]
  For any $S \in \calC$, we denote by $S^0 \coloneqq 0^*(S) \in \CAlg(\calC_{A/}^{\boxdot})$ its image.
  Note that under the equivalence $\Mod_{A}(\calC)^{\odot} \simeq \calC_{A/}^{\boxdot}$, the symmetric monoidal functor $0^*$ corresponds to the trivial module functor $\mathsf{Triv} \colon \calC^{\amalg} \to \Mod_{A}(\calC)^{\odot}$.

  Assume further that $\calC$ admits all geometric realizations.
  By \cite[4.4.2.9 and 5.5.8.7]{lurie-htt}, the forgetful functor $\calC_{A/} \to \calC$ detects sifted colimits (and pushouts), so that
  $(\calC_{A/})^{\boxdot}$ satisfies the two requirements of \cref{relativetensor}.
  It follows that the \icategory of modules $\Mod_{S^0}(\calC_{A/})$ inherits a symmetric monoidal structure $\otimes_{S^0}$ (the relative tensor product of modules over $S^0$).
  We get an \ifunctor (with base-change functorialities)
  \[
    \begin{tikzcd}[row sep=-.3em]
      \calC \ar{r} & \monoidalcats \\
      S \ar[mapsto]{r} & \Mod_{S^0}(\calC_{A/}^{\boxdot})^{\otimes_{S^0}}\rlap{.}
    \end{tikzcd}
  \]
\end{construction}
\begin{construction}\label{relativeconvoverS}
  Let $\calC \in \inftycats$ have finite colimits and let $A \in \CoCAlg(\calC^{\amalg})$.
  For every $S \in \calC$, we have $S \amalg A \in \CoCAlg(\calC_{S/}^{\amalg})$.
  We get in particular a relative convolution symmetric monoidal category $(\calC_{S/})_{S \amalg A/}^{\boxdot}$.
  Notice that by \cite[\href{https://kerodon.net/tag/018L}{Cor.\,018L}]{kerodon}, we have at the level of the underlying \icategories $(\calC_{S/})_{S \amalg A/} \simeq \calC_{S \amalg A/}$.
  We get an \ifunctor (via pushout functorialities)
  \[
    \begin{tikzcd}[row sep=-.3em]
      \calC \ar{r} & \monoidalcats \\
      S \ar[mapsto]{r} & (\calC_{S/})_{S \amalg A/-}^{\boxdot}\rlap{.}
    \end{tikzcd}
  \]
\end{construction}

\begin{proposition}\label{originofcomodulestructure}
  Let $\calC$ be a \icategory with finite colimits and geometric realizations.
  Let $A \in \CoCAlg(\calC^{\amalg})$.
  The \ifunctors $\calC \to \monoidalcats$ of \cref{coalgebraonSzero,relativeconvoverS} are equivalent.
  In particular, for any $S \in \calC$, we have a symmetric monoidal equivalence
  \[
    \Mod_{S^0}(\calC_{A/}^{\boxdot})^{\otimes_{S^0}} \simeq (\calC_{S/})_{S \amalg A/}^{\boxdot}.
  \]
  If moreover $\calD \subset \calC$ is a full sub-\icategory stable under finite colimits
  and if $A \in \CoCAlg(\calD^{\amalg})$, the above functorial equivalence restricts to a symmetric monoidal equivalence functorial in $S \in \calD$:
  \[
    \Mod_{S^0}(\calD_{A/}^{\boxdot})^{\otimes_{S^0}} \simeq (\calD_{S/})_{S \amalg A/}^{\boxdot}.
  \]
  \begin{proof}
    Applying \cref{modulesiscocartesian} we find $(\calC_{S/})^{\amalg} \simeq \Mod_{S}(\calC^{\amalg})^{\otimes_S}$.
    Therefore, using \cref{alternativedescriptionconvolutionascomodules} and \cite[Thm.\,4.8.4.6]{lurie-ha} and \cite[Rem.\,2.2.8]{Arpon}, we obtain symmetric monoidal equivalences
    \begin{multline*}
      (\calC_{S/})_{S \amalg A/}^{\boxdot}
      \simeq \Mod_{S \amalg A} \left(\calC_{S/}^{\amalg}\right)^{\odot}
      \simeq \Mod_{S \amalg A} \left(\Mod_{S}(\calC^{\amalg})^{\otimes_S}\right)^{\odot}\\
      \simeq \Mod_A(\calC^{\amalg})^\odot \otimes_\calC \Mod_S(\calC)^{\otimes_S}.
    \end{multline*}
    On the other hand, the same results assemble into
    \[
      \Mod_{S^0}(\calC_{A/}^{\boxdot})^{\otimes_{S^0}} \simeq \Mod_{S^0}(\Mod_A(\calC^{\amalg})^\odot)^{\otimes_{S^0}} \simeq \Mod_A(\calC^{\amalg})^\odot \otimes_\calC \Mod_S(\calC)^{\otimes_S},
    \]
    noticing that $S^0$ is the image of $S$ under the (symmetric monoidal) trivial module functor $\calC^{\amalg} \to \Mod_A(\calC^{\amalg})^\odot$.
    As every step is functorial, this equivalence is functorial. The restriction to a full subcategory $\calD \subset \calC$ stable under finite colimits is obvious.
  \end{proof}
\end{proposition}

\begin{remark}
  Let $\calD$ be an \icategory with finite colimits.
  It follows from \cite[\S5.3.6]{lurie-htt} that there exist an \icategory $\calC$ with finite limits and geometric realizations and a fully faithful embedding $\calD \hookrightarrow \calC$ preserving finite colimits.
  We can thus always apply \cref{originofcomodulestructure} on such an \icategory $\calD$ with finite colimits (and the resulting equivalence is independent of the choice of $\calC$).
\end{remark}

\subsubsection{Lax monoidal lax functors}

We now discuss base change with respect to the relative convolution product of \cref{originofcomodulestructure}.
We place ourselves in the following setting:
\begin{setting}\label{settinglaxmonoidallaxfunctors}
  Let $\calC \in \inftycats$ with finite colimits and let $A \in \CoCAlg(\calC^{\amalg})$.
  Let $\calD^{\otimes}$ be a symmetric monoidal \icategory. Assume that $\calD$ admits geometric realizations and that the tensor product preserves them in each variable.
  Let $\sfM^{\otimes} \colon \calC_{A/}^{\boxdot} \to \calD^{\otimes}$ be a symmetric monoidal \ifunctor.
\end{setting}

\begin{observation}\label{laxmonoidalMrelative}
  Assume \cref{settinglaxmonoidallaxfunctors}.
  For every object $S \in \calC$, we get using \cref{originofcomodulestructure} a \emph{lax} symmetric monoidal functor
  \[
    \sfM_S^{\otimes} \colon (\calC_{S/})_{S \amalg A/}^{\boxdot} \simeq \Mod_{S^0}(\calC_{A/}^{\boxdot})^{\otimes_{S^0}} \to \Mod_{\sfM(S^0)}(\calD)^{\otimes_{\sfM(S^0)}}
  \]
  which fails to be strongly symmetric monoidal in general.
  Moreover, if $u \colon T \to S$ is a morphism in $\calC$, a priori all we have is an induced lax natural transformation given by the universal property of the relative bar construction as a colimit:
  \[
    \begin{tikzcd}
      (\calC_{S/})_{S \amalg A/}^{\boxdot} \ar{r}{\sfM_S^{\otimes}} \ar{d}[swap]{- \amalg_{S} T} & \Mod_{\sfM(S^0)}(\calD)\mathrlap{{}^{\otimes_{\sfM(S^0)}}} \ar{d}{- \otimes_{\sfM(S^0)} \sfM(T^0)} \ar[Rightarrow]{dl}
      \\
      (\calC_{T/})_{T \amalg A/}^{\boxdot} \ar{r}[swap]{\sfM_T^{\otimes}} & \Mod_{\sfM(T^0)}(\calD)\mathrlap{{}^{\otimes_{\sfM(T^0)}}}\rlap{.}
    \end{tikzcd}
  \]
\end{observation}

\begin{notations}{\label{notationslaxdiagrams}}
  \item \label{ioperads} We denote by $\inftyoperads$ the \icategory of \ioperads, by $\laxmonoidal \subseteq \inftyoperads$ the \emph{full} subcategory spanned by symmetric monoidal \icategories.
  Note that $\monoidalcats \subset \laxmonoidal$ is a (non-full) sub-\icategory (imposing the functors to be symmetric monoidal).

  \item Let $\laxfunctors$ denote the \icategory whose objects are \ifunctors $\Delta^1 \to \laxmonoidal$, \ie lax symmetric monoidal functors $\calC_0^{\otimes} \to \calC_1^{\otimes}$ (\cf \cref{ioperads}) and morphisms given by lax monoidal natural transformations
  \[
    \begin{tikzcd}[row sep=small]
      \calC_{0}^{\otimes} \ar{r} \ar{d} & \calC_{1}^{\otimes} \ar{d} \ar[Rightarrow, shorten=1ex]{dl}
    \\
      \calD_{0}^{\otimes} \ar{r} & \calD_{1}^{\otimes}\rlap{.}
    \end{tikzcd}
  \]
  See \cite[Thm.\,E]{zbMATH07785229} and \cite[Def.\,3.9]{MR4367222}) for a definition in terms of \emph{Gray tensor products}.

  \item Let $\laxmatters \subseteq \laxfunctors$ denote the non-full subcategory containing all objects but with morphisms where the two vertical functors are symmetric monoidal.
\end{notations}

The next construction assembles relative convolution products with their (lax) base change structure:

\begin{construction}\label{constructiongeneralMFnaturaltransformation}
  In the \cref{settinglaxmonoidallaxfunctors}, the functoriality for modules (see \cite[4.5.3.1]{lurie-ha}) implies the collection of lax monoidal \ifunctors $\sfM^{\otimes}_S$ of \cref{laxmonoidalMrelative} can be assembled to an \ifunctor
  \[
    \begin{tikzcd}[row sep=-.3em, column sep=small]
      \mathllap{\underline{\sfM}^{\otimes} \colon{}} \calC \ar{r} & \laxmatters
      \\
      S \ar[mapsto]{r} & \sfM_S^{\otimes}\rlap{.}
    \end{tikzcd}
  \]
\end{construction}

\subsection{Action, base-change and relative Thom--Sebastiani \texorpdfstring{$\MF^\infty$}{MF∞}}
\label{subsectionrelativethomsebastiani}
We can now extract further consequences of the Thom--Sebastiani \cref{thomsebastianiMF}, by specializing the results of \cref{sectionrelativeconvolution}.

\subsubsection{Categories of matrix factorizations as modules}

\begin{construction}\label{constructionMFtoFunlax}
  We apply \cref{relativeconvoverS} to $\calC = \dInfSch_\basefield^\op$ with $A = \fAffLine \in \CAlg(\dInfSch_\basefield)$ and obtain a functor
  \[
    \begin{tikzcd}[row sep=-.5em]
      \mathllap{\stackofinfschemes^{\boxplus}_{\fAffLine} \colon {}} \dInfSch_\basefield^\op \ar{r} & \monoidalcats \\
      S \ar[mapsto]{r} & \dInfSch_{S \times \smash{\fAffLine}}^{\boxplus} \mathrlap{{}\coloneqq (\dInfSch_{/S})_{/S \times \smash{\fAffLine}}^{\boxplus}.}
    \end{tikzcd}
  \]
  where each $\dInfSch_{S \times \smash{\fAffLine}}^{\boxplus}$ is equipped with a \emph{relative convolution} over $S$ informally given by $(\formalU,f) \smboxplus (\formalV,g) = (\formalU \times_S \formalV, f\boxplus g)$.

  Moreover, applying \cref{constructiongeneralMFnaturaltransformation} to the symmetric monoidal \ifunctor from \cref{constructionMFasconvolutionmonoidal}
  \[
    \MF^{\infty,\boxplus} \colon \calC_{A/}^{\boxdot} \coloneqq \dInfSch_{\fAffLine}^{\boxplus,\op} \to \Corr\left(\dInfSch_{\fAffLine}\right)^{\boxplus} \to \calD^{\otimes} \coloneqq \PrLStmonoidal[\kbetaLaurent]
 \]
  (since by \cref{definitionPrlStmonoidalrelative} $\calD$ verifies the assumptions in \cref{settinglaxmonoidallaxfunctors} ), we obtain an \ifunctor
  \[
    \begin{tikzcd}[row sep=-.3em, column sep=small]
      \mathllap{\myuline{\MF}^\infty \colon{}} \dInfSch_\basefield^\op \ar{r} & \laxmatters
      \\
      S \ar[mapsto]{r} & \MF^{\infty}_S\rlap{.}
    \end{tikzcd}
  \]
  In essence, it encodes a lax notion of morphism $\stackofinfschemes^{\boxplus}_{\fAffLine} \rightsquigarrow \PrLStmonoidal[\MF^\infty(-,0)]$: for any $S \in \dInfSch_\basefield$ we have a lax symmetric monoidal functor
  \[
    \begin{tikzcd}[row sep=-.5em]
      \mathllap{\MF^\infty_S \colon{}} \dInfSch_{S \times \smash{\fAffLine}}^{\boxplus,\op} \ar{r} & \PrLStmonoidal[\MF^\infty(S,0)]
      \\
      (\fAffLine \from^f \formalU \to S) \ar[mapsto]{r} & \MF^\infty(\formalU, f)
    \end{tikzcd}
  \]
  and for any $\xi \colon T \to S \in \dInfSch_\basefield$, a (lax monoidal) natural transformation $\theta_\xi$
  \[
    \begin{tikzcd}
      \dInfSch_{S \times \smash{\fAffLine}}^{\boxplus,\op} \ar{r}{\MF^\infty_S} \ar{d}[swap]{- \underset{S}{\times} T} & \PrLStmonoidal[\MF^\infty(S,0)] \ar{d}{\mathrlap{- \underset{\MF^\infty(S,0)}{\otimes} \MF^\infty(T,0)}} \ar[Rightarrow]{dl}{\theta_\xi}
      \\
      \dInfSch_{T \times \smash{\fAffLine}}^{\boxplus,\op} \ar{r}[swap]{\MF^\infty_T} & \PrLStmonoidal[\MF^\infty(T,0)]\rlap{.}
    \end{tikzcd}
  \]
\end{construction}

\begin{corollary}\label{MFhasactionofDmodules}
  For any $\fAffLine \from^f \formalU \to S \in \dInfSch_\basefield$, the \icategory $\MF^\infty(\formalU, f)$ carries a natural action of $\MF^\infty(S, 0) \simeq \IndCoh(S)\ofbetaLaurent \coloneqq \IndCoh(S) \otimes_\basefield \kbetaLaurent$ (\cf \cref{notationsmodulesincategoriesbeta} and \cref{MFzerofunction}).
  In particular, when $S = \formalU_{\derham}$ and $\formalU \to S$ is the canonical map, the \icategory of 2-periodic D-modules $\DQCoh(\formalU_\derham)\ofbetaLaurent$ naturally acts on $\MF^\infty(\formalU,f)$.
  \begin{proof}
    The first claim is a direct consequence of \cref{constructionMFtoFunlax}. The second one follows from the equivalence $\DQCoh(\formalU_\derham) \simeq \IndCoh(\formalU_\derham)$ of \cref{reminderequivalenceICQCderham}.
  \end{proof}
\end{corollary}

\begin{remark}
  The above \cref{constructionMFtoFunlax} (and the contents of \cref{sectionrelativeconvolution}) can be informally unpacked as follows.
  For any $S \in \dInfSch_\basefield$, the diagonal morphism $S \to S \times S$ endows $(S,0) \coloneqq (0 \colon S \to \fAffLine) \in \dInfSch_{\fAffLine}$ with a coalgebra structure for $\smboxplus$.
  For any $f \colon \formalU \to \smash{\fAffLine}$, a morphism $\formalU \to S$ then amounts to a coaction of $(S,0)$ on $(\formalU, f)$.
  Applying the symmetric monoidal functor $\MF^\infty \colon \dInfSch_{\fAffLine}^\op \to \PrLSt[\kbetaLaurent]$, we get the action of $\MF^\infty(S,0)$ on $\MF^\infty(\formalU, f)$.
  The above construction further establishes the functorialities of this action.
\end{remark}

\subsubsection{Base change of categories of matrix factorizations}

\begin{definition}\label{definitionmorphismhasbasechangeMF}
  We say that a morphism $\xi \colon T \to S \in \dInfSch_\basefield$ has \emph{base change for matrix factorizations} if for any $S \from \formalU \to^f \fAffLine \in \dInfSch_\basefield$ such that $\formalU \to S$ is a reduced equivalence, the induced box product functor
  \[
    \theta_\xi(\formalU,f) \colon \MF^\infty(\formalU,f) \otimes_{\IndCoh(S)\ofbetaLaurent} \IndCoh(T)\ofbetaLaurent \to \MF^\infty\left(\formalU \times_S T, f \smboxplus 0\right)
  \]
  is an equivalence.
\end{definition}

The base-change results of \cref{basechangeIC-redeqetale,basechangeIC-overDmodgeneral} imply similar statements at the level of matrix factorizations categories:

\pagebreak
\begin{lemmalist}{}
  \item \label{basechangeMF-redeqetale} If $\xi \colon T \to S \in \dInfSch_\basefield$ is formally étale, then it has base change for matrix factorizations.
  \item \label{basechangeMF2}Let $Y \in \dInfSch_\basefield$. Any $\xi \colon T \to S \coloneqq Y_\derham \in \dInfSch_\basefield$ has base change for matrix factorizations.
\end{lemmalist}
\begin{proof}
  Fix $S \from \formalU \to^f \fAffLine \in \dInfSch_\basefield$ such that $\formalU \to S$ is a reduced equivalence.
  We apply \cref{basechangeIC-redeqetale} or \cref{basechangeIC-overDmodgeneral} to the formally étale morphism between the zero loci $(\formalU \times_S T)_0 \simeq \formalU_0 \times_S T \to \formalU_0$ to deduce the analogous statement for $\PreMF$.
  The result follows by inverting $\beta$.
\end{proof}

\begin{proposition}\label{basechangeformatrixfactorizations}
  If $\xi \colon T \to S \in \dInfSch_\basefield$ is schematic, smooth and $S$ is $1$-affine, then $\xi$ has base change for matrix factorizations.
\end{proposition}
\begin{proof}
  Denote by $\pi \colon X \coloneqq T \times_S S_\red \to S_\red$ the projection. By assumption $X$ is a scheme and $\pi$ is smooth.
  By the local étale form of smooth maps (see \cite[\href{https://stacks.math.columbia.edu/tag/054L}{054L}]{stacks-project}),
  there exists a Zariski covering $X' \twoheadrightarrow X$ and an étale morphism $a \colon X' \to S_\red \times \bbA^n$ such that the following diagram commutes
  \[
    \begin{tikzcd}[row sep=scriptsize]
      X' \ar{r} \ar{d}[swap]{a} & X \ar{d}{\pi} \\ S_\red \times \bbA^n \ar{r} & S_\red\rlap{.}
    \end{tikzcd}
  \]
  Denote by $\zeta \colon T' \to T$ the Zariski covering corresponding to $X' \to X$ ($X$ and $T$ have the same underlying topological space).
  The morphism $X' \to S_\red \times \bbA^n \to \bbA^n$ extends to the formal thickening $T'$ of $X'$ and yields $b \colon T' \to S \times \bbA^n$.
  We then have a Cartesian diagram
  \[
    \begin{tikzcd}[row sep=scriptsize]
      X' \ar{r} \ar{d}[swap]{a} \tikzcart & T' \ar{d}{b} \\
      S_\red \times \bbA^n \ar{r} & S \times \bbA^n\rlap{.}
    \end{tikzcd}
  \]
  As $a$ is étale (and $S_\red \to S$ is a reduced equivalence), it implies that so is $b$ by the topological invariance of small étale sites.

  Now, if $S$ is 1-affine so is $T$ (\cf \cref{smoothoverregularisregular}) and it follows from \cref{basechangeisconservative,basechangeMF-redeqetale} that $\xi$ has base change for matrix factorizations if and only if $\xi \circ \zeta$ does.
  As $\xi \circ \zeta \colon T' \to S$ factors as $T' \to^b S \times \bbA^n \to S$ with $b$ étale, the result follows from \cref{basechangeMF-redeqetale} and the equivalence $\IndCoh(S \times \bbA^n) \simeq \IndCoh(S) \otimes \IndCoh(\bbA^n)$ of \cref{IndCohassymmetricmonoidalfunctoronindinfschemes}.
\end{proof}

\subsubsection{Relative Thom--Sebastiani equivalences}

By \cref{constructionMFtoFunlax}, the relative version of the functor $\MF^\infty$ is a priori only lax symmetric monoidal.
We prove here strict monoidality in the special case where the base is a de Rham stack (see \cref{relativeTSoverderham} below).

\begin{proposition}\label{relativeTSfollowsfrombasechange}
  Let $Y \in \dInfSch_\basefield$ and $(\formalU, f)$, $(\formalV, g) \in \dInfSch_{\fAffLine \times Y}$.
  Denote by $\formalU_0$ and $\formalV_0$ the $0$-loci of $f \colon \formalU \to \fAffLine$ and $g \colon \formalV \to \fAffLine$ respectively.
  If the box product functor
  \[
    \IndCoh(\formalU_0) \otimes_{\IndCoh(Y)}\IndCoh(\formalV_0) \to \IndCoh\left( \formalU_0 \times_Y \formalV_0 \right)
  \]
  is an equivalence, then so is the box product functor
  \[
    \MF^\infty(\formalU, f) \otimes_{\IndCoh(Y)\ofbetaLaurent} \MF^\infty(\formalV, g) \to \MF^\infty\left(\formalU \times_Y \formalV, f \smboxplus g\right).
  \]
\end{proposition}
\begin{proof}
  For any integer $n$, we compute using \cite[Thm.\,4.8.4.6]{lurie-ha}:
  \begin{align*}
    \IndCoh(\formalU_0) & \otimes_{\IndCoh(Y)} \IndCoh(\bbB_Y)^{\otimes n} \otimes_{\IndCoh(Y)} \IndCoh(\formalV_0)
    \simeq \IndCoh(\formalU_0) \otimes_{\IndCoh(Y)} \IndCoh(\formalV_0) \otimes_{\IndCoh(Y)} \IndCoh(\bbB_Y)^{\otimes n}
    \\& \simeq \IndCoh\left(\formalU_0 \times_Y \formalV_0\right) \otimes_{\IndCoh(Y)} \Mod_{\kbetaLaurent^{\otimes n}}(\IndCoh(Y))
    \simeq \Mod_{\kbetaLaurent^{\otimes n}}\left(\IndCoh\left(\formalU_0 \times_Y \formalV_0\right)\right)
    \\& \simeq \IndCoh\left( \formalU_0 \times_Y \bbB_Y^n \times_Y \formalV_0 \right).
  \end{align*}
  As in the proof of \cref{thomsebastianiMF}, the functor
  \[
    F \coloneqq \boxtimes \colon \PreMF\left(\formalU, q\right) \otimes_{\IndCoh(Y)\ofbeta} \PreMF(\formalV, f) \to \PreMF\left(\formalU \times_Y \formalV, q \smboxplus f\right)
  \]
  thus arises from the augmented simplicial inf-scheme
  \[
    \begin{tikzcd}
      \cdots \ar[shift left=6pt]{r} \ar{r} \ar[shift right=6pt]{r} & \displaystyle
      \formalU_0 \times_Y \bbB_Y \times_Y \formalV_0 \ar[shift left=3pt]{r} \ar[shift right=3pt]{r} & \displaystyle
      \formalU_0 \times_Y \formalV_0 \ar{r}{i} &
      \displaystyle
      \Bigl(\formalU \times_Y \formalV\Bigr)_{\!0}.
    \end{tikzcd}
  \]
  We get that $F$ is an equivalence using \cite[I-3 Prop.\,3.3.3 (b)]{GR-II} as in the proof of \cref{thomsebastianiMF}, as $i$ is proper and surjective on $\basefield$-points.
\end{proof}

\begin{corollary}[Relative Thom--Sebastiani over de Rham stacks]\label{relativeTSoverderham}
  Consider $\formalU \to^f \fAffLine \in \dInfSch_\basefield$ and $\formalU_\derham \from \formalV \to^g \fAffLine \in \dInfSch_\basefield$.
  The box product functor
  \[
    \MF^\infty(\formalU, f) \otimes_{\IndCoh(\formalU_\derham)\ofbetaLaurent} \MF^\infty(\formalV, g) \to \MF^\infty\left(\formalU \times_{\formalU_\derham} \formalV, f \smboxplus g\right)
  \]
  is an equivalence.
  In particular, $\MF^\infty(\formalU, f)$ is dualizable over $\IndCoh(\formalU_\derham)\ofbetaLaurent$, with dual $\MF^\infty(\formalU, -f)$.
\end{corollary}
\begin{proof}
  Apply \cref{basechangeIC-overDmodgeneral,relativeTSfollowsfrombasechange}. The dualizable property follows from the same arguments as in \cref{MFinftydualizable} replacing $(\ast,0)$ by $(\formalU_\derham,0)$ and the convolution by the relative convolution $-\underset{(\formalU_\derham,0)}{\boxplus}-$.
\end{proof}

\begin{remark}
  Since dualizable categories are hypercomplete, \cref{relativeTSoverderham} also implies that the $\DQCoh(\formalU_\derham)\ofbetaLaurent$\nobreakdashes-module $\MF^\infty(\formalU, f)$ is hypercomplete (\cf \cref{hypercompletenonaffinecase2periodic}).
\end{remark}

\subsection{Quadratic bundles and Clifford algebras}
\label{subsection-actionquadraticbundles}
We relate here matrix factorizations of quadratic functions to representations of Clifford algebras.
This will then allow us to extend our relative Thom--Sebastiani \cref{relativeTSoverderham} to any base, when one of the functions involved is a non-degenerate quadratic form.

\subsubsection{Clifford algebras}
Classically, if $(Q,q) \to S$ is a quadratic bundle over a smooth scheme $S$, then its Clifford algebra $\clifford_S(Q,q)$ is the $\quot{\integers}{2}$-graded associative algebra over $S$:
\begin{equation}\label{naiveClifford}
  \quot{\left(\bigoplus Q^{\otimes n}\right)}{\{v \otimes v - q(v)\}} \in \Alg\left(\DQCoh(S)\ofbetaLaurent \right).
\end{equation}
Clifford algebras are naturally related to matrix factorizations of quadratic functions (see \eg \cite[Thm.\,9.3.4]{MR3121870} or \cref{MFQisClifford} below):
\[
  \MF^\infty(Q,q) \simeq \Mod_{\clifford_S(Q,q)} \left( \DQCoh(S)\ofbetaLaurent \right).
\]
We will now see that the notion of Clifford algebras and the above equivalence extend to the general setup of derived inf-schemes.
\begin{construction}\label{constructionofcliffordalgebra}
  Let $S \in \dInfSch_\basefield$ and $(Q,q) \to S$ be a quadratic bundle.
  Denote by $\widehat{Q}$ the formal completion of $Q$ along its $0$-section $S \to Q$, and by $Q_0$ the $0$-locus of the restriction of $q$ to $\widehat Q$.
  We have reduced equivalences
  \[
    \begin{tikzcd}
      S \ar{r}{s_0} \ar["\Id", rounded corners, to path={|- ([yshift=.7em]\tikztotarget.north)[near end]\tikztonodes -- (\tikztotarget)}]{rr} & Q_0 \ar{r} & S.
    \end{tikzcd}
  \]
  The inf-scheme $Q_0$ moreover comes equipped with an action of the relative derived group scheme $\bbB_S \coloneqq \Omega_0 \affineline{} \times S$.
  Via the equivalences of \cite[Part\,II Chap.\,7, Thm.\,3.6.2 and Chap.\,5, Thm.\,1.6.4]{GR-II}, the thickening $Q_0$ corresponds to some Lie algebra $L$ in $\IndCoh(S)$ and $\bbB_S$ to the abelian Lie algebra $\omega_S[-2]$.
  Thus the action $\bbB_S \circlearrowright Q_0$ corresponds to an action $\omega_S[-2] \circlearrowright L$ (in Lie algebras in $\IndCoh(S)$).
  Passing to enveloping algebras (\cf \cite[Part II, Chap. 8, \S 4.2.1]{GR-II}), we find $\mathcal U_S(L) \in \Alg(\IndCoh(S)\ofbeta)$.

  We set $\preclifford_S(Q,q) \coloneqq \mathcal U_S(L) \in \Alg(\IndCoh(S)\ofbeta)$ and
  \[
    \clifford_S(Q,q) \coloneqq \preclifford_S(Q,q) \otimes_{\kbeta} \kbetaLaurent = \mathcal U_S(L) \otimes_{\kbeta} \kbetaLaurent \in \Alg(\IndCoh(S)\ofbetaLaurent).
  \]
\end{construction}
\begin{definition}
  Let $S \in \dInfSch_\basefield$ and $(Q,q) \to S$ be a quadratic bundle.
  The $2$-periodic dg-algebra $\clifford_S(Q,q)$ of \cref{constructionofcliffordalgebra} is called the Clifford algebra of $(Q,q)$ over $S$.
\end{definition}
\begin{remark}
  If $S \in \dInfSch_\basefield$ is such that $\Upsilon_S \colon \DQCoh(S) \to \IndCoh(S)$ is an equivalence (\eg for $S$ a smooth scheme or $S = T_\derham$ for $T$ a scheme), the Lie algebra $L$
  can be computed in $\DQCoh(S)$:
  \[
    L \simeq \tangent_{S/Q_0} \simeq s_0^*\tangent_{Q_0/S}[-1] \simeq \left(
    \begin{tikzcd}[row sep=0]
        \mathrm{deg.} & 0 & 1 & 2 & 3 \\
        & 0 \ar{r} & Q \ar{r}{0} & \structuresheaf_{S} \ar{r} & 0
      \end{tikzcd}\right)
  \]
  with bracket given by $q \colon \Sym^2_S Q \to \structuresheaf_{S}$.
  The action of $\beta$ on $Q_0$ amounts to the multiplication by $1 \in \structuresheaf_{S}[-2] \subset \tangent_{S/Q_0} \subset \mathcal{U}_S(\tangent_{S/Q_0})$. (Note that the action is by $\structuresheaf_S$ instead of $\Upsilon_S(\structuresheaf_S) = \omega_S$ as we are working in $\DQCoh(S)$.)
  Inverting $\beta$, we then observe that $\clifford_S(Q,q) \coloneqq \mathcal{U}_S(\tangent_{S/Q_0})[\beta^{-1}] \in \Alg(\DQCoh(S)\ofbetaLaurent)$ coincides with the usual definition of \cref{naiveClifford}.
\end{remark}

\begin{proposition}\label{MFQisClifford}
  Let $S \in \dInfSch_\basefield$ and $(Q,q) \to S$ be a quadratic bundle and denote by $s \colon (S,0) \to (Q,q)$ the zero section. Then, there is a natural $\IndCoh(S)\ofbetaLaurent$\nobreakdashes-linear (\resp $\IndCoh(S)\ofbeta$\nobreakdashes-linear) equivalence commuting with the forgetful functor
  \[
  	\begin{tikzcd}[row sep=small, column sep={between origins,6em}]
      \MF^\infty(Q,q) \ar{rr}{\sim} \ar[end anchor={[xshift=-1.2em]}]{dr}[swap,very near start]{s^!} && \ar[end anchor={[xshift=1.3em]}]{dl}[very near start]{\forget} \Mod_{\clifford_S(Q,q)}(\IndCoh(S)\ofbetaLaurent)
    \\
      & \MF^\infty(S,0) = \IndCoh(S)\ofbetaLaurent &
    \end{tikzcd}
  \]
  (and similarly $\PreMF(\widehat Q,q) \simeq \Mod_{\preclifford_S(Q,q)}(\IndCoh(S)\ofbeta)$).

  \begin{proof}
    By \cite[Part\,II Chap.\,7, Prop.-Const.\,5.1.2]{GR-II}, \cref{constructionofcliffordalgebra} and the fact $\mathcal{U}_S(L)$ is in $\Alg(\IndCoh(S)\ofbeta)$, we have
    \[
      \PreMF(\widehat Q,q) \coloneqq \IndCoh(Q_0) \simeq \Mod_{\mathcal{U}_S(L)}(\IndCoh(S)) \simeq \Mod_{\mathcal{U}_S(L)}(\IndCoh(S)\ofbeta)
    \]
    commuting with the forgetful functors. Inverting $\beta$, we find
    \[
      \MF^\infty(Q,q) \simeq \MF^\infty(\widehat Q,q) \simeq \Mod_{\clifford_S(Q,q)}(\IndCoh(S)\ofbetaLaurent).\qedhere
    \]
  \end{proof}
\end{proposition}

\subsubsection{Quadratic Thom--Sebastiani and Knörrer periodicity}

\begin{theorem}[{see also \cite[Thm.\,9.1.7]{MR3121870}}]\label{relativeThomSebastianiMF}
  Let $Y \in \dInfSch_\basefield$ and let $(Q,q) \to Y$ be a (non-degenerate) quadratic bundle.
  Denote by $\widehat Q$ the formal neighborhood of the $0$-section of $Y \to Q$.
  Then the following holds:
  \begin{enumerate}[beginpenalty=10000, label={\textrm{\upshape{(\alph*)}}}, ref={\textrm{\upshape{(\alph*)}}}]
    \item\label{assertionrelativeTS}(Relative Thom--Sebastiani for quadratic bundles)
    Let $(\formalU, f) \in \dInfSch_{Y \times \fAffLine}$ and set $\widehat Q^\formalU \coloneqq \widehat Q \times_{Y} \formalU$ endowed with the function $q \smboxplus f \colon \widehat Q^\formalU \to \fAffLine$.
    The exterior product functor
          \[
            \boxtimes \colon \MF^\infty(\widehat Q,q) \otimes_{\IndCoh(Y)\ofbetaLaurent} \MF^\infty(\formalU, f) \to \MF^\infty\left(\widehat Q^\formalU, q \smboxplus f\right)
          \]
          is an equivalence.
    \item\label{assertionKnorrer}(Knörrer periodicity) Let $i \colon L \hookrightarrow Q$ be a maximal isotropic subbundle. Denote by $\widehat L$ the formal neighborhood of the $0$-section $Y \to L$ and by $p \colon \widehat L \to Y$ the structural projection, the correspondence $\smash{(\widehat Q,q) \from^i (\widehat L,0) \to^{p} (Y, 0)}$ induces an $\IndCoh(Y)\ofbetaLaurent$\nobreakdashes-linear equivalence
          \[
            \begin{tikzcd}
              \MF^\infty(\widehat Q,q) \ar{r}{p_*i^!} & \MF^\infty(Y, 0) \simeq \IndCoh(Y)\ofbetaLaurent.
            \end{tikzcd}
          \]
  \end{enumerate}
\end{theorem}
\begin{examples}{The above theorem applies in any of the following cases:}
  \item If $Y$ is any derived scheme, thus improving on \cite[Thm.\,9.1.7(i)]{MR3121870},
  \item More generally, if $Y$ is the formal completion of a scheme along a closed subset.
  \item If $Y = S_\derham$ for some scheme $S$.
\end{examples}

\begin{proof}[Proof of \cref{relativeThomSebastianiMF}]
  \cref{MFQisClifford} and \cite[Thm.\,4.8.4.6]{lurie-ha} imply the equivalence:
  \begin{multline*}
    \IndCoh(\widehat{Q}_0) \otimes_{\IndCoh(Y)} \IndCoh(\formalU_0)
    \simeq \Mod_{\preclifford_Y(Q,q)}(\IndCoh(Y)) \otimes_{\IndCoh(Y)} \IndCoh(\formalU_0)
    \\
    \simeq \Mod_{\preclifford_Y(Q,q)}(\IndCoh(\formalU_0))
    \simeq \IndCoh\left( \widehat{Q}_0 \times_Y \formalU_0 \right).
  \end{multline*}
  Assertion \ref{assertionrelativeTS} thus follows from \cref{relativeTSfollowsfrombasechange}.

  We now focus on assertion \ref{assertionKnorrer} and fix a maximal isotropic $L \hookrightarrow Q$.
  The correspondence $\smash{(\widehat Q,q) \from^i (\widehat L,0) \to^p (Y,0)}$ naturally lives over $Y$.
  The induced functor $p_*i^!$ is thus indeed $\IndCoh(Y)\ofbetaLaurent$\nobreakdashes-linear.
  Setting $\theta_Y \coloneqq p_* i^!$, we get by writing $Y$ as a colimit of derived affine schemes $S$ as in \cref{infschemesascolimits}
  \[
    \begin{tikzcd}[row sep=small]
      \mathllap{\MF^\infty(\widehat Q, q) \simeq{}} \Mod_{\clifford_Y(Q,q)}(\IndCoh(Y))
      \ar[phantom, start anchor=center, end anchor=center]{d}[sloped]{\simeq} \ar{r}{\theta_Y} & \IndCoh(Y)\ofbetaLaurent \ar[phantom, start anchor=center, end anchor=center]{d}[sloped]{\simeq}
      \\
      \displaystyle \lim_{S \to Y} \Mod_{\clifford_S(Q \times_Y S,q)}(\IndCoh(S)) \ar{r}{\lim \theta_S} &
      \displaystyle \lim_{S \to Y} \IndCoh(S)\ofbetaLaurent\rlap{.}
    \end{tikzcd}
  \]
  where the commutativity of the diagram follows from base change.

  We can thus reduce to the case where $Y$ is affine.
  Up to replacing $Y$ by an étale cover, we can further assume both $(Q,q)$ and the maximal isotropic $L \hookrightarrow Q$ are trivial: \ie $L \simeq L^0 \times Y$ is a trivial vector bundle and $(Q,q) \simeq (Q^0, q^0) \times Y$ is the hyperbolic bundle $(L^0 \oplus L^{0\dual}) \times Y$.
  Using \cref{thomsebastianiMF}, the functor of interest becomes
  \[
    \begin{tikzcd}[row sep=tiny]
      \MF^\infty(\widehat Q,q) \ar{r}{p_*i^!} \ar[phantom]{d}[sloped]{\simeq} &[2em] \IndCoh(Y)\ofbetaLaurent \ar[phantom]{d}[sloped]{\simeq}
    \\
      \displaystyle \MF^\infty(\widehat Q^0,q^0) \otimes_{\kbetaLaurent} \IndCoh(Y)\ofbetaLaurent \ar{r}{(p_*i^!) \otimes \Id} & \displaystyle \MF^\infty(*, 0) \otimes_{\kbetaLaurent} \IndCoh(Y)\ofbetaLaurent\rlap{.}
    \end{tikzcd}
  \]
  It is thus enough to treat the case where $Y = \Spec \basefield$.
  This case now follows from \cite[\S5.3]{MR2824483}: the functor $p_*i^!$ is right adjoint to $i_* p^*$, which is the functor considered in loc. cit.
\end{proof}

\begin{corollary}\label{MFquadselfinverse}
  Let $Y \in \dInfSch_\basefield$ be $1$\nobreakdashes-affine.
  Let $(Q,q) \to Y$ be a non-degenerate quadratic bundle and denote by $\widehat Q$ the formal neighborhood of its $0$-section.
  Then the \icategory $\MF^\infty(\widehat Q,q)$ is invertible as a $\IndCoh(Y)\ofbetaLaurent$\nobreakdashes-module.
  Moreover, the choice of a square root of $-1$ in $\basefield$ induces an $\IndCoh(Y)\ofbetaLaurent$\nobreakdashes-linear equivalence $\MF^\infty(\widehat Q,q)^{-1} \simeq \MF^\infty(\widehat Q,q)$.
\end{corollary}

\begin{proof}
  Consider the obvious diagonal maximal isotropic $Q \hookrightarrow (Q \oplus Q, q \smboxplus -q)$.
  Applying \cref{relativeThomSebastianiMF} (\ref{assertionrelativeTS} then \ref{assertionKnorrer}), we get a (natural) trivialization
  \[
    \MF^\infty(\widehat Q,q) \otimes_{\IndCoh(Y)\ofbetaLaurent} \MF^\infty(\widehat Q,-q) \simeq \MF^\infty(\widehat Q \smboxplus_Y \widehat Q,q \smboxplus -q) \simeq \IndCoh(Y)\ofbetaLaurent.
  \]
  It follows that $\MF^\infty(\widehat Q,q)$ is indeed invertible with inverse $\MF^\infty(\widehat Q,-q)$.
  We then observe that a choice of a square root of $-1$ identifies $(Q,-q)$ with $(Q,q)$, concluding our proof.
\end{proof}

\section{\texorpdfstring{$\MF^\infty$}{MF∞} as a map from the Darboux stack}\label{formalsectionMFfromdarboux}

In this section we start combining the results of this paper with the results from our previous work \cite{hennionholsteinrobaloI}.

\subsection{Functorialities of  \texorpdfstring{$\MF^\infty$}{MF∞} on derived inf-schemes}\label{axiomatizationofMFinvariants}
The goal of this section is purely formal. We construct an \ifunctor
\[\myuline{\MF}^{\infty, \acts, \twotorsion}_{\simplex-\constant\times \derham}\colon \Delta^\op \times \dInfSch_\basefield^\op \to \Fun(\Delta^1, \catmod)\]
(see \cref{factorizationthrough2torsion2} below),  assembling  the several results about matrix factorization categories of derived inf-schemes proven in \cref{MFsection} (functoriality, base change, Thom-Sebastiani, Kn\"orrer and so on), together with the simplicial coherences encoding isotopy relations.
The trusting reader not interested in those $\infty$-categorical details can safely skip this section.

\subsubsection{Packaging the action of quadratic bundles and Thom-Sebastiani}
\label{ThomSebastianiPackage}
In this subsection we obtain an \ifunctor $\myuline{\MF}^{\infty, \acts}$ -- see \cref{constructionMFwithactionfunctorialityTS} below -- encoding the action of quadratic bundles and the Thom-Sebastiani theorem for matrix factorization categories on derived inf-schemes, assembling the results of  \cref{MFsection}.

\begin{construction}\label{mapofstackquadmonoidal}
  Recall from \crossref{definitionBO} the abelian monoid stack in groupoids classifying non-degenerate quadratic bundles \[\B\orthogonal^{\oplus} \colon \dSt_\basefield^\op \to \CAlg(\inftygpd)\subseteq \monoidalcats.\]
  Consider its restriction $\B\orthogonal_{|_\dInf}^{\oplus}$ along the inclusion $\dInfSch_\basefield\subseteq \dSt_\basefield$.
  Since the formal completion along the zero section of a quadratic bundle on a derived inf-scheme is again a derived inf-scheme with a pro-nilpotent function, taking formal completions along the zero section defines a symmetric monoidal morphism of stacks $\dInfSch^\op \to \monoidalcats$ (see also \crossref{orthogonaltoliouv})
  \[
    \B\orthogonal_{|_\dInf}^{\oplus} \to \stackofinfschemes^{\boxplus}_{\fAffLine}
  \]
  (\cf \cref{constructionMFtoFunlax}) corresponding to an \ifunctor
  \[
    \begin{tikzcd}[row sep=-.3em]
      \mathllap{\myuline{\rmA} \colon{}} \dInfSch^\op \ar{r} & \Fun(\Delta^1, \monoidalcats)
      \\
      S \ar[mapsto]{r} & \left(((Q,q) \to S) \mapsto (\widehat Q \to^q \fAffLine \times S) \right)\rlap{.}
    \end{tikzcd}
  \]
\end{construction}

\begin{definition}\label{Azumayaprestack2periodic}
  Consider symmetric monoidal prestack given by the target of $\myuline{\MF}^\infty$ in \cref{constructionMFtoFunlax}. It is given by
  \[
  \begin{tikzcd}[row sep=-.3em]
    \dInfSch_\basefield^\op \ar{r}{\myuline{\MF}^\infty} & \laxmatters \ar{r}{\ev_1} & \monoidalcats
    \\
    S \ar[mapsto]{rr} && \PrLStmonoidal[\MF^\infty(S,0)]\rlap{.}
  \end{tikzcd}
  \]
  We define the prestack of \emph{2-periodic} invertible modules $\PrLStinvertible[\MF^\infty(-,0)] \colon \dInfSch_\basefield^\op \to \monoidalcats$ by setting $\PrLStinvertible[\MF^\infty(S,0)] \subseteq \PrLStmonoidal[\MF^\infty(S,0)]$ as the full subcategory spanned by invertible $\MF^\infty(S,0)$-modules.
  By construction $\PrLStinvertible[\MF^\infty(S,0)]$ has a symmetric monoidal structure which is grouplike.
  The inclusion $\PrLStinvertible[\MF^\infty(S,0)]\subseteq \PrLStmonoidal[\MF^\infty(S,0)]$ defines a symmetric monoidal morphism of symmetric monoidal prestacks which we view as an \ifunctor $\myuline{\rmC} \colon \dInfSch_\basefield^\op \to \Fun(\Delta^1, \monoidalcats)$.
\end{definition}

The next construction puts together the equivariant structure on $\myuline{\MF}^\infty$ with respect to the action of quadratic bundles. In order to accommodate future generalizations (for instance, the case of small categories of matrix factorizations not treated in this paper) we give here a general axiomatic framework that assembles the main actors and isolates the key properties of the functor $\MF^\infty$.

\begin{construction}\label{assembleQuadAndMF}
  Consider the following data:
  \begin{enumerate}[(D1)]
    \item A symmetric monoidal \icategory $\calD^{\otimes}$ and a symmetric monoidal functor
    \[
      \sfM^{\otimes} \colon \dInfSch_{\fAffLine}^{\boxplus,\op} \to \calD^{\otimes}
    \]
    as in \cref{settinglaxmonoidallaxfunctors} and $\myuline{\sfM}^{\otimes} \colon \dInfSch_\basefield^\op \to \laxmatters$
    the output of
    \cref{constructiongeneralMFnaturaltransformation} obtained by passing to module-objects.

    \item $\myuline{\rmA}^{\op} \colon \dInfSch^\op \to \Fun(\Delta^1, \monoidalcats)$ the functor of \cref{mapofstackquadmonoidal} classifying $\B\orthogonal_{|_\dInf}^{\oplus,\op} \to \stackofinfschemes^{\boxplus,\op}_{\fAffLine}$

  \end{enumerate}
  and set the notations:
  \begin{enumerate}[(N1)]
    \item
    $\myuline{\calD}^{\otimes} \colon \dInfSch_\basefield^\op \to \monoidalcats$ the target of the natural transformation $\myuline{\sfM}^{\otimes}$, sending $S\mapsto \Mod_{\myuline{\sfM}(S,0)}(\calD)^{\otimes}$ as in \cref{laxmonoidalMrelative}.

    \item $\myuline{\calD}^\invertible \colon \dInfSch_\basefield^\op \to \monoidalcats$ the symmetric monoidal full sub-prestack of $\myuline{\calD}^{\otimes}$ given on $S$ by the full subcategory of  $\myuline{\calD}^{\otimes}(S)$ spanned by tensor-invertible $\sfM(S,0)$-modules.
    We see the inclusion $\myuline{\calD}^\invertible \subseteq \myuline{\calD}^{\otimes}$ as an \ifunctor
    \[
    \myuline{\rmC} \colon \dInfSch_\basefield^\op \to \Fun(\Delta^1, \monoidalcats).
    \]
  \end{enumerate}
  We now specify two axioms:
  \begin{enumerate}[label=(A{{\arabic*}})]
    \item The composition with the lax morphism of functors associated to $\myuline{\sfM}^{\otimes}$
    \begin{equation}\label{restrictiontoquadratic}
      \begin{tikzcd}[row sep=small]
        \B\orthogonal_{|_\dInf}^{\oplus,\op}\ar{r}{\myuline{\rmA}^\op}& \stackofinfschemes^{\boxplus,\op}_{\fAffLine}\ar{r}{\myuline{\sfM}^{\otimes}} &\myuline{\calD}^{\otimes}
      \end{tikzcd}
    \end{equation} is strongly symmetric monoidal and the lax natural transformation structure induced from that of $\myuline{\sfM}^{\otimes}$ is strict.

    \item The composition \eqref{restrictiontoquadratic} factors through $\myuline{\calD}^\invertible \subseteq \myuline{\calD}^{\otimes}$, \ie
    \[
      \begin{tikzcd}[row sep=small, column sep=small]
        \ar[dashed]{d}\B\orthogonal_{|_\dInf}^{\oplus,\op}\ar{r}{\myuline{\rmA}^\op} & \stackofinfschemes^{\boxplus,\op}_{\fAffLine} \ar{d}{\myuline{\sfM}^{\otimes}}
      \\
        \myuline{\calD}^\invertible \ar[hookrightarrow]{r}{\myuline{\rmC}}&\myuline{\calD}^{\otimes}\rlap{.}
      \end{tikzcd}
    \]
  \end{enumerate}
  Let us reformulate this: by construction, $\myuline{\rmA}^\op$, $\myuline{\sfM}^{\otimes}$ and $\myuline{\rmC}$ can be assembled as part of a commutative diagram (using \cref{notationslaxdiagrams})
  \[
    \begin{tikzcd}[row sep=scriptsize, column sep=tiny]
      \Fun(\Delta^1, \monoidalcats) \ar{d}{\ev_1} & \dInfSch_{\basefield}^\op \ar{d}{\myuline{\sfM}^{\otimes}} \ar{l}[swap]{\myuline{\rmA}^\op} \ar{r}{\myuline{\rmC}} & \ar{d}{\ev_1}\Fun(\Delta^1, \monoidalcats)
    \\
      \monoidalcats & \laxmatters
      \ar{l}[swap]{\ev_0} \ar{r}{\ev_1} & \monoidalcats
    \end{tikzcd}
  \]
  Axiom (A1) guarantees that the composition
  \[
    \begin{tikzcd}[column sep=small, row sep=scriptsize]
      \dInfSch_{\basefield}^\op \ar[shorten <=-.2em]{d}{(\myuline{\rmA}^\op,\myuline{\sfM}^{\otimes},\myuline{\rmC})}
    \\
      \Fun(\Delta^1, \monoidalcats) \fiberproduct{\monoidalcats} \laxmatters \fiberproduct{\monoidalcats} \Fun(\Delta^1, \monoidalcats) \ar[shorten <=-.8em]{d}[pos=0.1]{\circ \times \Id}
    \\
      \laxmatters \fiberproduct{\monoidalcats} \Fun(\Delta^1, \monoidalcats)
    \end{tikzcd}
  \]
  lands in the subcategory $\Fun(\Delta^1, \monoidalcats) \times_{\monoidalcats} \Fun(\Delta^1, \monoidalcats) $ and axiom (A2) guarantees the existence of a lift of the composition above
  \[
    \begin{tikzcd}[column sep=small, row sep=scriptsize]
      & \ar[dashed, rounded corners=2em, "", to path={-| (\tikztotarget)[near start]\tikztonodes}]{dl} \dInfSch_{\basefield}^\op \ar[shorten <=-.2em]{d}
    \\
      \Fun(\Delta^2, \monoidalcats) \ar{r}{\partial_1 \times \partial_0} & \Fun(\Delta^1, \monoidalcats) \times_{\monoidalcats} \Fun(\Delta^1, \monoidalcats)\rlap{.}
    \end{tikzcd}
  \]
  The categorical equivalences $\Delta^1\times \Delta^1 \simeq \Delta^2 \amalg_{\Delta^1} \Delta^2$ and $\Delta^2 \simeq \horn^1_2$ allow us compute the fiber product
  \[
  \begin{tikzcd}[column sep=-0.5in,row sep=0.1in]
  	&  \Fun(\Delta^1, \monoidalcats) \fiberproduct{\monoidalcats} \laxmatters \fiberproduct{\monoidalcats} \Fun(\Delta^1, \monoidalcats) \ar[shorten <=-.8em]{d}[pos=0.1]{\circ \times \Id}\\
  	  \Fun(\Delta^2, \monoidalcats) \ar{r}{\partial_1 \times \partial_0} & \laxmatters \times_{\monoidalcats} \Fun(\Delta^1, \monoidalcats)\rlap{.}
  \end{tikzcd}
  \]
  and finally exhibit the data of $(\myuline{\rmA},\myuline{\sfM}^{\otimes},\myuline{\rmC})$ as an \ifunctor 
  \[
    \begin{tikzcd}
      \dInfSch_{\basefield}^\op\ar{r} &\calE\subseteq \Fun(\Delta^1, \laxmatters)
    \end{tikzcd}
  \]
  where $\calE$ is the non-full subcategory spanned by:
  \begin{itemize}
    \item objects corresponding to \ifunctors $\Delta^1\times \Delta^1\to\laxmonoidal$ classifying homotopy commutative diagrams
    \[
      \begin{tikzcd}[row sep=small]
        \calC_{00}^{\otimes}\ar{r}\ar{d} & \calC_{01}^{\otimes}\ar{d}
        \\
        \calC_{10}^{\otimes}\ar{r}& \calC_{11}^{\otimes}
      \end{tikzcd}
    \]
    where all arrows are strongly symmetric monoidal except $\calC^{\otimes}_{10} \to \calC_{11}^{\otimes}$ which may only be lax symmetric monoidal;
    \item morphisms are given by diagrams
    \[
      \begin{tikzcd}[row sep=small, column sep=small]
        \ar{dr}[swap]{\rho^{\otimes}_{00}} \calC_{00}^{\otimes} \ar{rr} \ar{dd} && \calC_{01}^{\otimes} \ar{dd} \ar{dr}{\rho^{\otimes}_{01}} &
      \\
        & \calV_{00}^{\otimes} \ar[crossing over]{rr} && \calV_{01}^{\otimes} \ar{dd}
      \\
        \calC_{10}^{\otimes} \ar{rr} \ar{dr}{\rho^{\otimes}_{10}} && \calC_{11}^{\otimes}\ar[Rightarrow]{dl} \ar{dr}{\rho^{\otimes}_{11}}&
      \\
        & \calV_{10}^{\otimes} \ar{rr} \ar[from=uu, crossing over] && \calV_{11}^{\otimes}
      \end{tikzcd}
    \]
    where $\rho^{\otimes}_{00}$, $\rho^{\otimes}_{01}$, $\rho^{\otimes}_{10}$, $\rho^{\otimes}_{11}$ are strongly symmetric monoidal and all the faces are homotopy commutative except the bottom face which comes together with the data of a natural transformation.
  \end{itemize}
  Finally, recall from \cite[4.8.3.19, 2.4.2.5,4.2.1.13,]{lurie-ha} the \icategory $\catmod$ classifying pairs $(\calC^{\otimes}, \calM)$ where $\calC^{\otimes}$ is a symmetric monoidal \icategory and $\calM$ is a left $\calC^{\otimes}$-module in $\inftycats^{\times}$. The equivalence \[\Fun(\Delta^1, \monoidalcats)\simeq \CAlg(\catmod) \] of \cite[3.4.1.3]{lurie-ha} composed with the forgetul functor \[\CAlg(\catmod) \to \catmod\]
   sends a symmetric monoidal functor $F \colon \calC^{\otimes} \to \calV^{\otimes}$ to the pair $(\calC^{\otimes}, \calV)$ with $\calV$ seen as a left $\calC^{\otimes}$-module via $F$. This construction induces an \ifunctor
  \[
    \begin{tikzcd}[column sep=small]
      \myuline{\sfM}^{\acts} \colon \dInfSch_\basefield^\op\ar{r} & \calE \ar{r} & \Fun^\lax(\Delta^1, \catmod)
    \end{tikzcd}
  \]
  where $\Fun^\lax(\Delta^1, \catmod)$ denotes the \icategory whose:
  \begin{itemize}
    \item objects are given by \ifunctors $\Delta^1 \to \catmod$ classifying pairs of functors $(F, \phi) \colon (\calC_0^{\otimes}, \calM_0) \to (\calC_1^{\otimes}, \calM_1)$ where $F$ is a symmetric monoidal functor and $\phi \colon \calM_0 \to \calM_1$ is an \ifunctor together with an equivariant structure with respect to the actions;
    \item morphisms are given by
    \[
      \begin{tikzcd}[row sep=small]
        (\calC_{0}^{\otimes}, \calM_0)\ar{r}{(F_1,\phi_1)}\ar{d}& (\calC_{1}^{\otimes}, \calM_1)\ar{d}\\
        (\calV_{0}^{\otimes}, \calN_0)\ar{r}{(F_2,\phi_2)}& (\calV_{1}^{\otimes}, \calN_1)
      \end{tikzcd}
    \]
    where the two vertical maps are in $\catmod$, together with homotopies rendering the commutativity of the square of symmetric monoidal functors, a lax natural transformation $\eta$
    rendering the lax commutativity of the functors on underlying modules
    \[
    \begin{tikzcd}[row sep=small]
      \calC_{0}^{\otimes}\ar{r}\ar{d}& \calC_{1}^{\otimes}\ar{d}\\
      \calV_{0}^{\otimes}\ar{r}& \calV_{1}^{\otimes}
    \end{tikzcd}\hspace{0.5cm}\begin{tikzcd}[row sep=small]
    \calM_{0}\ar{r}\ar{d} & \ar[Rightarrow]{dl}{\eta}\calM_{1}\ar{d}\\
    \calN_{0}\ar{r}& \calN_{1}
  \end{tikzcd}\]
  and a compatible lax natural transformation rendering the equivariant compatibility of the two $\calC_0^{\otimes}$ actions on $\calN_1$
  \[
    (\calC_0^{\otimes}\to \calC_1^{\otimes}\to \calV_1^{\otimes} \acts \calN_1) \Rightarrow (\calC_0^{\otimes}\to \calV_0^{\otimes}\to \calV_1^{\otimes} \acts \calN_1).
  \]
  \end{itemize}
\end{construction}

\begin{construction}\label{constructionMFwithactionfunctorialityTS}
  We apply \cref{assembleQuadAndMF} to the \ifunctor $\sfM^{\otimes}=\MF^{\infty,\boxplus}$ of \cref{constructionMFtoFunlax}. In this case, the grouplike prestack $\myuline{\calD}^\invertible$ coincides with the  prestack  $\PrLStinvertible[\MF^\infty(-,0)]$ of \cref{Azumayaprestack2periodic} and
  \begin{itemize}
    \item Axiom (A1) is a consequence of the relative Thom--Sebastiani for quadratic bundles of \cref{relativeThomSebastianiMF}\ref{assertionrelativeTS};
    \item Axiom (A2) follows from \cref{relativeThomSebastianiMF}\ref{assertionrelativeTS} and \cref{MFquadselfinverse}.
  \end{itemize}
  As an output, we obtain an \ifunctor
  \[
    \begin{lrbox}{\subdiagram}
      \begin{tikzcd}[row sep=tiny, column sep=scriptsize, cramped]
        \dInfSch_{S \times \smash{\fAffLine}}^{\op} \ar{r}{\myuline{\MF}^\infty_S} & \PrLSt[\MF^\infty(S,0)]
      \\
        \B\orthogonal^{\oplus,\op}_{|_\dInf}(S) \ar[phantom,start anchor=center, end anchor=center]{u}{\circlearrowleft} \ar{r}{\myuline{\MF}^\infty_S} & \PrLStinvertible[\MF^\infty(S,0)] \ar[phantom,start anchor=center, end anchor=center]{u}{\circlearrowleft}
      \end{tikzcd}
    \end{lrbox}
    \begin{tikzcd}[column sep=small, row sep=tiny]
      \mathllap{\myuline{\MF}^{\infty, \acts} \colon {}} \dInfSch_\basefield^\op \ar{r} & \Fun^\lax(\Delta^1, \catmod)
    \\
      S \ar[mapsto]{r} & \left[\usebox{\subdiagram}\right]
    \end{tikzcd}
  \]
  encoding the functoriality of relative matrix factorizations together with the Thom-Sebastiani theorem expressing the compatible action of quadratic bundles.
\end{construction}


\subsubsection{Packaging \texorpdfstring{$\affineline{}$}{A¹}-isotopy coherences, strong functoriality and duality.}
\label{sectionisotopycoherences}

\begin{notation}\label{labelalgebraicsimplices}
  Recall from \crossref{reminderexplicitmodelusualA1localization} the functor of algebraic simplices $\simplex^\bullet_\basefield \colon \Delta \to \dAff_\basefield$, defined for each $[n] \in \Delta$ as the subspace $\simplex^n \subseteq \affine^{n+1}$ cut out by the equation $\sum_{i=0}^{n} x_i=1$. For each $n$, the composition with the projection onto the first $n$-coordinates $\simplex^n \subseteq \affine^{n+1} \to \affine^n$ is an isomorphism. For a derived inf-scheme $S$, we denote by $\simplex^n_S \coloneqq \simplex^n_\basefield \times S$.
\end{notation}

\begin{construction}\label{MFfunctorpackisotopycoherences}
  By \cref{definitioninfschemes}, the de Rham functor $(-)_\derham \colon \dSt_\basefield \to \dSt_\basefield$ of \cref{formalstacks} restricts to a functor $(-)_\derham \colon \dInfSch_{\basefield} \to \dInfSch_{\basefield}$. Composing the functor $\myuline{\MF}^{\infty, \acts}$ of \cref{constructionMFwithactionfunctorialityTS} with the functor
  \begin{equation}\label{functordeltatimesderham}\Delta^\op \times \dInfSch_\basefield^\op \to \dInfSch_\basefield^\op,  \hspace{1cm} ([n], S) \mapsto \simplex^n_{S_\derham} \coloneqq S_\derham \times \simplex_\basefield^n
  	\end{equation}
  we obtain a new \ifunctor
\[
    \begin{lrbox}{\subdiagram}
      \begin{tikzcd}[row sep=tiny, column sep=tiny, cramped]
        \dInfSch_{\simplex^n_{S_\derham} \times \smash{\fAffLine}}^{\op} \ar{r}&\PrLSt[\MF^\infty(\simplex^n_{S_\derham},0)]
      \\
        \B\orthogonal^{\oplus,\op}_{|_\dInf}(\simplex^n_{S_\derham})\ar[phantom,start anchor=center, end anchor=center]{u}{\circlearrowleft} \ar{r} &\PrLStinvertible[\MF^\infty(\simplex^n_{S_\derham})] \ar[phantom,start anchor=center, end anchor=center]{u}{\circlearrowleft}
      \end{tikzcd}
    \end{lrbox}
    \myuline{\MF}^{\infty, \acts}_{\simplex\times  \derham} \colon
    \begin{tikzcd}[column sep=tiny, row sep=tiny, baseline=(X.base)]
      |[alias=X]| \Delta^\op \times \dInfSch_\basefield^\op \ar{r} & \Fun^\lax(\Delta^1, \catmod)
    \\
      ([n], S) \ar[mapsto]{r}  & \left[\usebox{\subdiagram}\right]
    \end{tikzcd}
  \]
 
\end{construction}

For the purpose of isotopy invariance we will only need to consider families of $\LG$-pairs that are constant in the $\simplex$-direction.

\begin{construction}\label{freeization}
  Since $[0]$ is an initial object in $\Delta^\op$, we have $\Fun(\Delta^\op, \inftycats) \simeq \Fun((\Delta^\op)_{[0]/}, \inftycats)$ and therefore every simplicial object $F$ in $\inftycats$ is canonically pointed by $F(0)$ via the unique degeneracy map $[n] \to [0]$. The functoriality of essential images -- see \cite[\href{https://kerodon.net/tag/05E7}{Prop.\,05E7}]{kerodon} -- allows us to produce an \ifunctor
  \[
    (-)^\constant \colon \Fun(\Delta^\op, \inftycats) \to \Fun(\Delta^\op, \inftycats)
  \]
  sending a simplicial object $[n] \mapsto F(n)$ to the simplicial object $[n] \mapsto \image(F(0)\to F(n))$ (see also \crossref{notationconstantobjects}).
  This construction preserves finite products and comes with a canonical functorial morphism $F^\constant \to F$.
\end{construction}

\begin{construction}\label{packagefree}
	Composing the functor $(-)^\constant$ of \cref{freeization}  with the functor   $\myuline{\MF}^{\infty, \acts}_{\simplex\times  \derham}$ of \cref{MFfunctorpackisotopycoherences}
we obtain a new \ifunctor
  \[
    \begin{lrbox}{\subdiagram}
      \begin{tikzcd}[row sep=tiny, column sep=tiny, cramped]
        \dInfSch_{\simplex^n_{S_\derham} \times \smash{\fAffLine}}^{\constant,\op} \ar{r}& \left(\PrLSt[\MF^\infty(\simplex^n_{S_\derham},0)]\right)^{\constant}
      \\
        \B\orthogonal^{\oplus,\op}_{|_\dInf}(\simplex^n_{S_\derham})^{\constant} \ar[phantom,start anchor=center, end anchor=center]{u}{\circlearrowleft} \ar{r} & \left(\PrLStinvertible[\MF^\infty(\simplex^n_{S_\derham})]\right)^{\constant} \ar[phantom,start anchor=center, end anchor=center]{u}{\circlearrowleft}
      \end{tikzcd}
    \end{lrbox}
    \myuline{\MF}^{\infty, \acts}_{\simplex\text{-}\constant\times \derham} \colon
    \begin{tikzcd}[column sep=0, row sep=tiny, baseline=(X.base)]
      |[alias=X]| \Delta^\op \times \dInfSch_\basefield^\op \ar{r} & \Fun^\lax(\Delta^1, \catmod)
    \\
      ([n], S)\ar[mapsto]{r}& \left[\usebox{\subdiagram}\right]
    \end{tikzcd}
  \]

\end{construction}

\begin{proposition} \label{finalMFfunctorpack}
  The functor $\myuline{\MF}^{\infty, \acts}_{\simplex\text{-}\constant\times \derham}$ of \cref{packagefree} factors by the full subcategory $\Fun(\Delta^1, \catmod)\subseteq \Fun^{\lax}(\Delta^1, \catmod)$. In particular, by adjunction it amounts to an \ifunctor
  \[
    \myuline{\MF}^{\infty, \acts}_{\simplex\text{-}\constant\times \derham} \colon \dInfSch_\basefield^\op \times \Delta^1 \times \Delta^\op \to \catmod.
  \]
  \begin{proof}
    By construction, we can forget that we land on $\simplex$-\emph{constant} objects in the $\Prl$-side (\cf  \cref{freeization}). All we need is to show that for every morphism $\rho \colon T \to S$ in $\dInfSch_\basefield$ and every morphism $\sigma \colon [m]\to [n]$ in $\Delta$, the induced morphism $\xi \coloneqq (\rho, \sigma) \colon T_\derham \times \simplex^m_\basefield \to S_\derham \times \simplex^n_\basefield$ in $\dInfSch_\basefield$ satisfies base change: \ie the restricted lax natural transformation $\theta_{\xi}$
    \[
      \begin{tikzcd}[row sep=large, column sep=large]
        \dInfSch_{\simplex^n \times S_\derham \times \smash{\fAffLine}}^{\thick,\constant,\op} \ar{r}{\MF^\infty_{\simplex^n \times S_\derham}} \ar{d}[swap]{- \underset{\simplex^n \times S_\derham}{\times} \simplex^m \times T_\derham}
      &
        \PrLSt[\MF^\infty(\simplex^n \OLDtimes S_\derham,0)] \ar{d}{\mathrlap{- \underset{\MF^\infty(\simplex^n \times S_\derham,0)}{\otimes} \MF^\infty(\simplex^m \times T_\derham,0)}} \ar[Rightarrow]{dl}{\theta_\xi}
      \\
        \dInfSch_{\simplex^m \times T_\derham \times \smash{\fAffLine}}^{\thick,\constant,\op} \ar{r}[swap]{\MF^\infty_{\simplex^m \times T_\derham}}
      &
        \PrLSt[\MF^\infty(\simplex^m \OLDtimes T_\derham,0)]
      \end{tikzcd}
    \]
    is an equivalence.
    The class of morphisms $\xi$ such that $\theta_\xi$ is invertible is stable under composition. We may thus assume $\xi$ is either $(\Id, \sigma) \colon (S_\derham, [n]) \to (S_\derham, [m])$ or $(\rho, \Id) \colon (S_\derham, [n]) \to (T_\derham, [n])$.
    The case of $\xi = (\rho, \Id)$ follows straightforwardly from \cref{basechangeMF2}.

    Now, the pullback functor $\simplex^n \times - \colon \dInfSch_{S_\derham \times \smash{\fAffLine}} \to \dInfSch_{\simplex^n \times S_\derham \times \smash{\fAffLine}}^{\constant}$ is by definition essentially surjective.
    In particular, the transformation $\theta_{(\Id, \sigma)}$ is invertible if and only if $\theta_{(\Id, \pi_n)}$ and $\theta_{(\Id, \pi_m)}$ are, where $\pi_p \colon [p] \to [0]$ denotes the unique morphism.
    The case $\theta_{(\Id, \pi_p)}$ follows trivially from \cref{thomsebastianiMF} (see also \cref{basechangeformatrixfactorizations}).
  \end{proof}
\end{proposition}

\pagebreak
\begin{observations}{}\label{recovernotationsforsheavesofcategoriesfromintro}
	\item \label{landsindualizable} Thanks to \cref{relativeTSoverderham},  the functor $\myuline{\MF}^{\infty, \acts}_{\simplex\text{-}\constant\times \derham}$ of \cref{finalMFfunctorpack} lands in \emph{dualizable modules}, \ie for each $([n],S) \in \Delta^\op  \times \dInfSch_\basefield^\op$ , $\myuline{\MF}^{\infty, \acts}_{\simplex\text{-}\constant\times \derham}$ takes values in $\PrLStdual[\MF^\infty(S_\derham\OLDtimes \simplex^n_\basefield)]$. This is stable under base change and under the action of invertible modules.
	\item
	Notice that since $S_\derham\times \simplex^n_\basefield$ is 1-affine (\cf \cref{1affinessreminder}), the \ifunctor $\Delta^\op\times \dInfSch_\basefield^\op\to \monoidalcats$  sending $(S, [n])\mapsto \PrLStdual[\MF^\infty(S_\derham\OLDtimes \simplex^n_\basefield)]$ by \cref{1affinebeta}-\cref{1affinebetadualizable} coincides symmetric monoidally with  the stack $\Delta^\op\times \dInfSch_\basefield^\op\to \monoidalcats$  sending $([n],S)\mapsto  \shcatsqcohtwoper(S_\derham\times \simplex^n_\basefield)$ of \cref{definitionofgoodstackofsheavesofcategories}.
	
	\item In the same way, the functor $([n],S)\mapsto \PrLStinvertible[\MF^\infty(S_\derham\OLDtimes\simplex^n_\basefield  )]$ coincides with $([n],S)\mapsto\Azumayatwoper(S_\derham\times\simplex^n_\basefield)$.
	
	\item Moreover, the action of $\PrLStinvertible[\MF^\infty(S_\derham\OLDtimes\simplex^n_\basefield)]$ on $\PrLStdual[\MF^\infty(S_\derham\OLDtimes \simplex^n_\basefield)]$ coincides functorially with that of $\Azumayatwoper$ on $\shcatsqcohtwoper$ by tensor product.
 
 \item \label{rewritingasmodulesinsheaves} The \ifunctor  $\myuline{\MF}^{\infty, \acts}_{\simplex\text{-}\constant\times \derham}$ of \cref{finalMFfunctorpack} can be re-interpreted via the equivalences of \icategories
 \begin{align*}
 	\Fun(\Delta^\op\times \dInfSch_\basefield^\op, &\Fun(\Delta^1, \catmod)) \\ 
 	 &\simeq \Fun(\Delta^1, \Fun(\Delta^\op , \Mod(\Fun( \dInfSch_\basefield^\op, \inftycats))))
 \end{align*}
 where we see $\Fun( \dInfSch_\basefield^\op, \inftycats)$ as a Cartesian symmetric monoidal category. The last equivalence follows from \cite[4.5.1.6, 2.4.2.5, Def. 2.4.2.1]{lurie-ha} describing module objects as monoid objects.
\end{observations}

\begin{construction}\label{constructionmotherearth}
	Following \cref{recovernotationsforsheavesofcategoriesfromintro}, the functor $\myuline{\MF}^{\infty, \acts}_{\simplex\text{-}\constant\times \derham}$ of \cref{finalMFfunctorpack} encodes an equivariant natural transformation of functors $\Delta^\op\to \Mod(\Fun(\dInfSch_\basefield^\op, \inftycats))$
	\begin{equation}\label{motherearth}
	\begin{tikzcd}[row sep=small]
		\stackofinfschemes^{\boxplus}_{\fAffLine}((-)_\derham\times \simplex^\bullet)^\constant\ar{r} & \shcatsqcohtwoper((-)_\derham\times \simplex^\bullet)^\constant\\
		\B\orthogonal^{\oplus,\op}_{|_\dInf}((-)_\derham\times \simplex^\bullet)^{\constant} \ar[phantom,start anchor=center, end anchor=center]{u}{\circlearrowleft} \ar{r}& \Azumayatwoper((-)_\derham\times \simplex^\bullet)^\constant \ar[phantom,start anchor=center, end anchor=center]{u}{\circlearrowleft}
	\end{tikzcd}
	\end{equation}
	where $\stackofinfschemes^{\boxplus}_{\fAffLine}$ is the functor of \cref{constructionMFtoFunlax}.
\end{construction}

\subsubsection{Packaging Kn\"orrer periodicity as 2-torsion}
\label{sectionpackagingtwotorsion}
We finally reach the last stage of our construction, addressing the remaining feature of matrix factorization categories not yet discussed, namely, Kn\"orrer periodicity of \cref{MFquadselfinverse}.

\begin{definition}\label{2torsionstackdefinition}
  We define the stack of of \emph{2-torsion} 2-periodic Azumaya algebras  $\Azumayatwopertwotorsion=\left(\Azumayatwoper\right)^{\twotorsion}$ as the homotopy fiber of the square map $(-)^{\otimes_2}:\Azumayatwoper\to \Azumayatwoper$, where $\Azumayatwoper:\dSt_\basefield^\op\to \monoidalcats$ is the abelian group stack of \cref{definitionofgoodstackofsheavesofcategoriesinvertible}.
\end{definition}

\begin{lemma}\label{factorizationthrough2torsion2}
  Fix a square root of $-1$ in $\basefield$. Then the map of symmetric monoidal stacks $\dInfSch_\basefield^\op\to \monoidalcats$
  \[
    \B\orthogonal^{\oplus,\op}_{|_\dInf} \to \Azumayatwoper
  \]
  factors through 2-torsion $\B\orthogonal^{\oplus,\op}_{|_\dInf} \to \Azumayatwopertwotorsion$.
  In particular, the equivariant map of \eqref{motherearth} lifts to an equivariant map
  \[
    \begin{tikzcd}[row sep=small]
      \stackofinfschemes^{\boxplus}_{\fAffLine}((-)_\derham\times \simplex^\bullet)^\constant\ar{r} & \shcatsqcohtwoper((-)_\derham\times \simplex^\bullet)^\constant\\
      \B\orthogonal^{\oplus,\op}_{|_\dInf}((-)_\derham\times \simplex^\bullet)^{\constant} \ar[phantom,start anchor=center, end anchor=center]{u}{\circlearrowleft} \ar{r}& \Azumayatwopertwotorsion((-)_\derham\times \simplex^\bullet)^\constant. \ar[phantom,start anchor=center, end anchor=center]{u}{\circlearrowleft}
    \end{tikzcd}
  \]
  encoded by an \ifunctor $\myuline{\MF}^{\infty, \acts, \twotorsion}_{\simplex-\constant\times \derham}\colon \Delta^1\to \Fun(\Delta^\op, \Mod(\Fun(\dInfSch_\basefield^\op, \inftycats)))$.
  \begin{proof}
    This is a reformulation of \cref{MFquadselfinverse}, applied in conjunction with \cref{constructionmotherearth}.
  \end{proof}
\end{lemma}

\subsection{Restriction of \texorpdfstring{$\MF^\infty$}{MF∞} to the Darboux stack}\label{actionquad}\label{matrixfactorizationasasectionofstack}

The goal of this section is to use $\myuline{\MF}^{\infty, \acts, \twotorsion}_{\simplex-\constant\times \derham}$ (\cf \cref{factorizationthrough2torsion2} ) to exhibit the functor of matrix factorizations $\MF^\infty$ as a map from the Darboux stack (\cf \crossref{definitiondarbouxstackasstacks}) to the stack of sheaves of crystals of 2-periodic dg-categories incorporating  the action of quadratic bundles, duality, Kn\"orrer periodicy and the isotopy coherences -- see \cref{MFoutofDarb} below. Most of the work has been done in \cref{axiomatizationofMFinvariants}.  It remains now to establish the comparison with the constructions appearing in \cite{hennionholsteinrobaloI}.

\begin{reminders}{}
  \item \label{stackofstacksLGconv}
    Recall from \crossref{liouvtoLGmonoidal} the symmetric monoidal stack
    \[
      \stackfactor^{\LG, \boxplus} \colon \Fun(\Delta^1,\dSt_\basefield)\to \monoidalcats
    \]
    sending $S \to S'$ to the \icategory of factorizations $S \to Y \to S'$ with $Y$ equipped with a function $Y \to \affineline{\basefield}$.
  \item \label{naturaltransformationep}
    Recall from \crossref{stackofstacksnofactorizations} and \crossref{notationfunctorpderham} the functors \vspace{-0.5em}
    \[
      \hspace{3em}\begin{tikzcd}[row sep={1.6em,between origins}]
        & \mathllap{e \colon{}} \dSt \ar{r} & \Fun(\Delta^1, \dSt), & S \mapsto (\varnothing \to S)
      \\
        \mathllap{p_X \colon{}} \dSt_{X_\derham} \ar{rr} && \Fun(\Delta^1, \dSt), &  S \mapsto \Big(S \fiberproduct{X_\derham} X \to S\Big)
      \\
        \mathllap{e_X \colon{}} \dSt_{X_\derham} \ar{r}{\forget} & \dSt \ar{r}{e} & \Fun(\Delta^1, \dSt), &
      \end{tikzcd}
    \]
    together with the canonical natural transformation $e_X \to p_X$. For a symmetric monoidal stack $\rmF \colon \Fun(\Delta^1, \dSt)^\op \to \monoidalcats$ we denote by $\rmF_{p_X}$, $\rmF_{e_X}$ the respective restrictions.

  \item \label{factLiouvilledInf} Recall from \crossref{constructionLiouville4mars} the symmetric monoidal stack \[\stackfactor^{\Liouv,\boxplus}: \Fun(\Delta^1,\dSt_\basefield)^\op\to \monoidalcats\] sending $S \to S'$ to the space of factorizations $S \to Y \to S'$ equipped with a section of the Liouville complex (\cf \crossref{Liouvillecomplex}). In particular, $Y$ comes with a function $f$ and following \crossref{liouvtoLGmonoidal} forgetting all Liouville data except the function induces a symmetric monoidal map of stacks
  \[
    \stackfactor^{\Liouv,\boxplus}\to\stackfactor^{\LG,\boxplus}\rlap{.}
  \]

\item \label{quadraticactionFactLiouv}	Let $\B\orthogonal^{\oplus}_{X_\derham} $ denote the restriction of the monoid stack $\B\orthogonal$ of \cref{mapofstackquadmonoidal} to $\dSt_{X_\derham}$. By \crossref{orthogonaltoliouv}, taking formal completion along the zero section defines a symmetric monoidal morphism of stack $\dSt^\op_\basefield \to \monoidalcats$
\[\B\orthogonal^{\oplus}_{X_\derham}\to \stackfactor^{\Liouv,\boxplus}_{p_X} \]

\item \label{quadraticactionmapfrompaper1} The natural transformation $e_X \to p_X$ of \cref{naturaltransformationep} induces a symmetric monoidal morphism of categorical stacks $\dSt_{X_\derham}^\op \to \monoidalcats$:
\[
\stackfactor^{\LG,\boxplus}_{p_X} \to
\stackfactor^{\LG,\boxplus}_{e_X}
\]
which we can compose with the maps from \cref{factLiouvilledInf,quadraticactionFactLiouv}:
\begin{equation}\label{symmetricmonoidalcompositionofpaper1}
	\B\orthogonal^{\oplus}_{X_\derham}\to\stackfactor^{\Liouv,\boxplus}_{p_X} \to \stackfactor^{\LG,\boxplus}_{p_X} \to
	\stackfactor^{\LG,\boxplus}_{e_X}.
\end{equation}
\end{reminders}

The main observation that allows us to relate the constructions in this paper to the ones appearing in \cite{hennionholsteinrobaloI} is the following:

\begin{observation}\label{firstcomparisionlemmatopaper1}
  By \crossref{constructionLiouville4marsmonoidal} and \crossref{stackofstacksnofactorizations}, the composition
  \[
    \begin{tikzcd}[row sep=-0.3em]
      \dSt_\basefield^\op \ar{r}{e} & \Fun(\Delta^1, \dSt_\basefield)^\op \ar{r}{\stackfactor^{\LG, \boxplus}} &[1em] \monoidalcats\phantom{{} \subset \dSt_\basefield^\op}
    \\
      S \ar[mapsto]{r} & (\varnothing \to S)
    \end{tikzcd}
  \]
  coincides with the \ifunctor $S\mapsto (\dSt_{/S})_{S\times \affineline{\basefield}}^{\boxplus}$
  of \cref{relativeconvoverS} applied to $\calC=\dSt_\basefield^\op$ and $A=\affineline{\basefield}$.
\end{observation}

\begin{notations}{}
	\item  Denote by $\dInfSch_{X_\derham} \subseteq \dSt_{X_\derham}$ the full subcategory of derived stacks $Y \to X_\derham$ with $Y$ a derived inf-scheme over $\basefield$.
	
	\item For a stack $\rmF:\dSt_{X_\derham}^\op \to \monoidalcats$, we denote by $\rmF_{|_{\dInf}}$ its restriction along the inclusion $\dInfSch_{X_\derham} \subseteq \dSt_{X_\derham}$.
	
	\item \label{matchofnotationsfrompaperone}We denote by $\stackfactor^{\LG,\dInf,\boxplus}_{e_X}\subseteq \left(\stackfactor^{\LG,\boxplus}_{e_X}\right)_{|_{\dInf}}$ the full symmetric monoidal substack spanned by pairs $(Y\to S,f)$ where $Y$ is a derived inf-scheme over $\basefield$ and the function is pro-nilpotent (\cf \cref{definitionpronilpotent}).
	
	\item We denote by $\stackfactor^{\Liouv,\dInf, \boxplus}_{p_X}:\dInfSch_{X_\derham}^\op \to \monoidalcats$ the full symmetric monoidal substack given by the fiber product
	\[
	\begin{tikzcd}[row sep=small, column sep=small]
		\stackfactor^{\Liouv,\dInf, \boxplus}_{p_X}\ar{r}\ar{d}\tikzcart&\ar{d}\left(\stackfactor^{\Liouv,\boxplus}_{p_X}\right)_{|_\dInf}\\
	\stackfactor^{\LG,\dInf,\boxplus}_{e_X}	\ar[hookrightarrow]{r}&\left(\stackfactor^{\LG,\boxplus}_{e_X}\right)_{|_{\dInf}}
	\end{tikzcd}
	\]
\end{notations}

The bridge between the constructions in this paper and the ones of \cite{hennionholsteinrobaloI} is given by the following:

\begin{observation}
  \label{bridgepaper1and2}
  Since taking  formals completion along the zero section of a quadratic bundle yields a derived inf-scheme with a pro-nilpotent function, the map of symmetric monoidal stacks $(\B\orthogonal^{\oplus}_{X_\derham})_{|_\dInf} \to \left(\stackfactor^{\LG,\boxplus}_{e_X}\right)_{|_\dInf}$ of \cref{quadraticactionmapfrompaper1} factors through the full substack $\stackfactor^{\LG,\dInf,\boxplus}_{e_X}$ of \cref{matchofnotationsfrompaperone}, corresponding to an \ifunctor
  \begin{equation}\label{inffunctorofpaper1}
    \dInfSch_\basefield^\op \to \Fun(\Delta^1, \monoidalcats).
  \end{equation}
  It now follows from \cref{firstcomparisionlemmatopaper1} that the following diagram commutes
  \[
    \begin{tikzcd}
      \dInfSch_{X_\derham}^\op \ar["\eqref{inffunctorofpaper1}", rounded corners, to path={|- ([yshift=.7em]\tikztotarget.north)[near end]\tikztonodes -- (\tikztotarget)}]{rr} \ar{r}[swap]{\textrm{forget}} &    \dInfSch_{\basefield}^\op \ar{r}[swap]{\myuline{\rmA}} & \Fun(\Delta^1, \monoidalcats)
    \end{tikzcd}
  \]
  with $\myuline{\rmA}$ classifying the morphism $\B\orthogonal_{|_\dInf}^{\oplus} \to \stackofinfschemes^{ \boxplus}_{\fAffLine}$ of \cref{mapofstackquadmonoidal}.
\end{observation}

\begin{construction}\label{naturaltransformationmodulespaper1} Following  \cref{bridgepaper1and2}, the symmetric monoidal morphism \eqref{symmetricmonoidalcompositionofpaper1} induces\footnote{via the functor  $\Fun(\Delta^1, \monoidalcats)\simeq \CAlg(\catmod) \to \catmod$ of \cite[3.4.1.3]{lurie-ha}.} an \ifunctor
	\begin{equation}\label{encodingactionpaper1}\dInfSch_{X_\derham}^\op\to \Fun(\Delta^1, \catmod)\end{equation}
	which as in \cref{rewritingasmodulesinsheaves} we can read as a natural transformation 
	\[
	\begin{tikzcd}[row sep=tiny, column sep=scriptsize]
		\stackfactor^{\Liouv,\dInf, \boxplus}_{p_X}  \ar{r}& 	\stackfactor^{\LG,\dInf,\boxplus}_{e_X}=\left(\stackofinfschemes^{ \boxplus}_{\fAffLine}\right)_{|_{X_\derham}}
		\\
		\left(\B\orthogonal^{\oplus}_{X_\derham}\right)_{|_\dInf} \ar[phantom,start anchor=center, end anchor=center]{u}{\circlearrowleft} \ar{r}{\Id} & \left(\B\orthogonal^{\oplus}_{\dInf}\right)_{|_{X_\derham}} \ar[phantom,start anchor=center, end anchor=center]{u}{\circlearrowleft}
	\end{tikzcd}\]
	where $(-)_{|_{X_\derham}}$ denotes the restriction along the forgetful functor \[\dInfSch_{X_\derham}^\op\to \dInfSch_{\basefield}.\]
\end{construction}

\begin{construction}\label{intermediatestep66}
  Combining \cref{definitioninfschemes} and \cref{formalstacks} and since $(X_\derham)_\derham\simeq X_\derham$, the de Rham functor $(-)_\derham \colon \dSt_\basefield \to \dSt_\basefield$ restricts to a functor $(-)_\derham \colon \dInfSch_{X_\derham} \to \dInfSch_{X_\derham}$.
  We consider a variant of the functor \eqref{functordeltatimesderham} of  \cref{MFfunctorpackisotopycoherences} given by
  \[
    \Delta^\op \times \dInfSch_{X_\derham}^{\op} \to \dInfSch_{X_\derham}^{\op},  \hspace{0.2cm} ([n], S\to X_\derham) \mapsto (\simplex^n_{\basefield}\times S_\derham \to \ast \times X_\derham\simeq X_\derham).
  \]
  By construction, it is compatible with \eqref{functordeltatimesderham} via the forgetful functor $\dInfSch_{X_\derham}\to \dInfSch_{\basefield}$. Composing with the \ifunctor \eqref{encodingactionpaper1} of \cref{naturaltransformationmodulespaper1}, followed by the constant simplicial procedure of \cref{freeization} (and using \cref{rewritingasmodulesinsheaves}), we obtain a new \ifunctor
  \[
    \Delta^1\to \Fun(\Delta^\op, \Mod(\Fun(\dInfSch_{X_\derham}^\op, \inftycats)))
  \]
  encoding the action on simplicially constant objects:
  \[
    \begin{tikzcd}[row sep=tiny, column sep=scriptsize]
      \stackfactor^{\Liouv,\dInf, \simplex-\constant\times \derham}_{p_X}  \ar{r}& 	\stackfactor^{\LG,\dInf,\constant\times \derham}_{e_X}=\left(\dInfSch_{\simplex^\bullet\times (-)_\derham \times \smash{\fAffLine}}^{\constant,\op}\right)_{|_{X_\derham}}
    \\ 
      \left(\B\orthogonal^{\oplus}_{X_\derham}\right)_{|_\dInf}(\simplex^\bullet\times(-)_\derham)^{\constant} \ar[phantom,start anchor=center, end anchor=center]{u}{\circlearrowleft} \ar{r}{\Id} & \left(\B\orthogonal^{\oplus}_{\dInf}\right)_{|_{X_\derham}}(\simplex^\bullet\times (-)_\derham)^{\constant} \ar[phantom,start anchor=center, end anchor=center]{u}{\circlearrowleft}
    \end{tikzcd}
  \]
  where the target coincides with the source of $\myuline{\MF}^{\infty, \acts, \twotorsion}_{\simplex-\constant\times \derham}$ (\cf \cref{factorizationthrough2torsion2}) restricted along the forgetful functor $\dInfSch_{X_\derham}^\op \to \dInfSch_\basefield$. Therefore, composing with the restriction  of  $\myuline{\MF}^{\infty, \acts, \twotorsion}_{\simplex-\constant\times \derham}$,  we obtain a new \ifunctor
  \[
    \myuline{\MF}^{\infty, \acts, \twotorsion, \Liouv}_{\simplex-\constant\times \derham}\colon	\Delta^1\to \Fun(\Delta^\op, \Mod(\Fun(\dInfSch_{X_\derham}^\op, \inftycats)))
  \]
  encoding the functor of matrix factorizations
  \begin{equation}\label{motherearth2}
    \begin{tikzcd}[row sep=tiny, column sep=scriptsize]
      \stackfactor^{\Liouv,\dInf, \constant\times \derham}_{p_X}  \ar{r}& \shcatsqcohtwoper((-)_\derham\times \simplex^\bullet)^\constant
    \\
      \left(\B\orthogonal^{\oplus}_{X_\derham}\right)_{|_\dInf}(\simplex^\bullet\times (-)_\derham)^{\constant} \ar[phantom,start anchor=center, end anchor=center]{u}{\circlearrowleft} \ar{r} & \Azumayatwopertwotorsion((-)_\derham\times \simplex^\bullet)^\constant\rlap{.} \ar[phantom,start anchor=center, end anchor=center]{u}{\circlearrowleft}
    \end{tikzcd}
  \end{equation}
\end{construction}

We now assume $X$ is a $(-1)$-shifted symplectic derived Deligne-Mumford $X$ equipped with its canonical exact structure $\lambda$ (\cf \crossref{canonicalexactstructure}).

\begin{reminder}\label{reminderstackliouvilleexact}
	 Recall from \crossref{fixedexactformstack} the stack \[\stackfactor^{\Liouv, \lambda, \boxplus}_{p_X} \colon \dSt_{X_\derham}^\op \to \monoidalcats\] classifying Liouville data together with a compatibility with the fixed exact structure $\lambda$, together with a forgetful functor $\stackfactor^{\Liouv,\lambda}_{p_X} \to \stackfactor^{\Liouv }_{p_X}$. 
	 
	 Following \crossref{orthogonaltoliouv}, the action of $\B\orthogonal^{\oplus}_{X_\derham}$ on $\stackfactor^{\Liouv, \boxplus}_{p_X}$ of \cref{mapofstackquadmonoidal} lifts to an action on the stack $\stackfactor^{\Liouv,\lambda}_{p_X}$, assembling as an \ifunctor
  \begin{equation}\label{EquivariantDarbinclusion0}
    \begin{lrbox}{\subdiagram}
      \begin{tikzcd}[row sep=tiny, column sep=tiny, cramped]
        \stackfactor^{\Liouv,\lambda}_{p_X}(S) \ar{r} & \stackfactor^{\Liouv}_{p_X}(S)
      \\
        \B\orthogonal^{\oplus}_{X_\derham}(S)\ar[phantom,start anchor=center, end anchor=center]{u}{\circlearrowleft} \ar{r} & \B\orthogonal^{\oplus}_{X_\derham}(S) \ar[phantom,start anchor=center, end anchor=center]{u}{\circlearrowleft}
      \end{tikzcd}
    \end{lrbox}
    \begin{tikzcd}[row sep=0pt, column sep=small]
      \dSt_{X_\derham}^\op \ar{r}&\Fun(\Delta^1, \catmod)&
    \\
      S\ar[mapsto]{r} & \left[\usebox{\subdiagram}\right]\rlap{.}
    \end{tikzcd}
  \end{equation}
\end{reminder}

\begin{construction}\label{finalstepfinally}
	Repeating the steps of \cref{intermediatestep66} with the stack $\stackfactor^{\Liouv, \lambda, \boxplus}_{p_X}$ of \cref{reminderstackliouvilleexact} and composing the maps  \eqref{EquivariantDarbinclusion0} and \eqref{motherearth2} we obtain a new \ifunctor
	\[
		\myuline{\MF}^{\infty, \acts, \twotorsion, \Liouv, \lambda}_{\simplex-\constant\times \derham}\colon \Delta^1\to \Fun(\Delta^\op, \Mod(\Fun(\dInfSch_{X_\derham}^\op, \inftycats)))
    \]
	encoding the functor of matrix factorizations
	\[
		\begin{tikzcd}[row sep=tiny, column sep=scriptsize]
			\stackfactor^{\Liouv,\lambda, \dInf, \constant\times \derham}_{p_X}  \ar{r}& 	\shcatsqcohtwoper((-)_\derham\times \simplex^\bullet)^\constant
			\\
			\left(\B\orthogonal^{\oplus}_{X_\derham}\right)_{|_\dInf}(\simplex^\bullet\times (-)_\derham)^{\constant} \ar[phantom,start anchor=center, end anchor=center]{u}{\circlearrowleft} \ar{r} & \Azumayatwopertwotorsion((-)_\derham\times \simplex^\bullet)^\constant\rlap{.} \ar[phantom,start anchor=center, end anchor=center]{u}{\circlearrowleft}
		\end{tikzcd}
	\]
\end{construction}

We now restrict to the small étale site of $X$ and finally recall the definition of the stack of quadratic bundles with flat connection and the Darboux stack:

\begin{reminders}{}
  \item Recall from \crossref{stackquaddef} the monoid stack of non-degenerate quadratic bundles on X with an étale locally trivial flat connection \[\stackQuadnabla_X \colon (X_\et^\daff)^\op \to \CAlg(\inftygpd).\] 
  By definition it is the substack of the restriction
  \[
    \begin{tikzcd}[row sep=0, cramped,baseline=(X.base)]
      |[alias=X]| X_\et^{\daff,\op} \ar{r}{(-)^X_\derham} & \dInfSch_{X_\derham}^\op \ar{r}{\B\orthogonal^{\oplus}_{X_\derham, \dInf}} &[1em] \CAlg(\inftygpd)\subseteq \monoidalcats
      \\
      S \ar[mapsto]{r} & S_\derham
    \end{tikzcd}
  \]
  spanned by locally trivial sections.

  \item By \crossref{definitiondarbouxstackasstacks}, the
  Darboux stack
  \[
    \Darbstack^{\lambda}_X \colon (X_\et^\daff)^\op \to \inftygpd
  \]
  is obtained by extracting the maximal $\infty$-groupoid of the substack of the restriction $\left(\stackfactor^{\Liouv, \lambda}_{p_X}\right)_{|_\dInf}\circ (-)^X_\derham$, classifying factorizations $S \to \formalU \to S_{\derham}$ where $S \to \formalU$ is a reduced equivalence and $\formalU$ is a smooth formal scheme equipped with a Lagrangian Liouville data compatible with $\lambda$.
  By \crossref{actionondarbouxcharts} the action of $\left(\B\orthogonal^{\oplus}_{X_\derham}\right)_{|_\dInf}$ on $\left(\stackfactor^{\Liouv, \lambda}_{p_X}\right)_{|_\dInf}\circ (-)^X_\derham$ restricts to an action of $\stackQuadnabla_X$ on $\Darbstack^{\lambda}_X$ which we assemble now as part of an \ifunctor
  \begin{equation}\label{EquivariantDarbinclusion1}
    \begin{lrbox}{\subdiagram}
      \begin{tikzcd}[row sep=tiny, column sep=tiny, cramped]
        \Darbstack^{\lambda}_X(S) \ar{r} & \stackfactor^{\Liouv, \lambda}_{p_X}(S_\derham)
      \\
        \stackQuadnabla_X(S) \ar[phantom,start anchor=center, end anchor=center]{u}{\circlearrowleft} \ar{r} & \B\orthogonal^{\oplus}_{X_\derham}(S_\derham) \ar[phantom,start anchor=center, end anchor=center]{u}{\circlearrowleft}
      \end{tikzcd}
    \end{lrbox}
    \begin{tikzcd}[row sep=0pt, column sep=small]
      (X_\et^\daff)^\op \ar{r}&\Fun(\Delta^1, \catmod)&\\
      S\ar[mapsto]{r} & \left[\usebox{\subdiagram}\right]\rlap{.}
    \end{tikzcd}
  \end{equation}
\end{reminders}

\begin{notation}\label{notationstackofcategoriesonetalesite}
  We denote by $\shcatsqcohtwoperfunctornabla_X$ the stack (\cf \cref{definitionofgoodstackofsheavesofcategories,twoperdgcatdescent}) on $X_\et^\dAff$ defined by sending $S$ to the maximal \igroupoid  $\shcatsqcohtwoper(S_\derham)^{\simeq}$. Similarly, we denote by $\Azumayatwopertwotorsionnabla_X$ the stack on  $X_\et^\dAff$ sending $S$ to the maximal \igroupoid  $\Azumayatwopertwotorsion(S_\derham)^{\simeq}$.
\end{notation}

\begin{remarks}{}
  \item The stack $\shcatsqcohtwoperfunctornabla_X$ has a more explicit description:
  for every $S \in X_\et^\daff$ we have
  \[
    \shcatsqcohtwoperfunctornabla_X(S) \coloneqq \shcatsqcohtwoper(S_\derham) \simeq \PrLStdual[\DQCoh(S_\derham)\ofbetaLaurent] \simeq \PrLStdual[\IndCoh(S_\derham)\ofbetaLaurent].
  \]
  This equivalence is furthermore functorial in $S$.
  Indeed, the equivalence follows from \cref{1affinebeta,YdR1affine}.
  Its functoriality stems from base change: see \cite[Lem.\,3.2.4]{gaitsgory_affineness} and \cite[\S 2.4.1]{gaitsgoryrozenblyum:crystals}.

  \item As $\DQCoh(S_\derham)$ is equivalent to the \icategory of complexes of $\calD$-modules on $S$, objects $\shcatsqcohtwoperfunctornabla_X(S)$ can be interpreted as presentable dg-categories with an (hypercomplete) action of the category of $2$-periodic complexes of $\calD$-modules on $S$.
\end{remarks}

\begin{proposition}\label{MFoutofDarb}
  Fix a square root of $-1$ in $\basefield$.
  Let $X$ a $(-1)$\nobreakdashes-shifted symplectic Deligne--Mumford stack equipped with its canonical exact structure $\lambda$ (see \crossref{canonicalexactstructure}).
  There is a morphism in $\Mod(\Sh(X_\et^\daff, \inftygpd)^{\times})$, which we denote by $\myuline{\MF}^\infty$:
  \[
  \begin{tikzcd}[row sep=small]
    \Darbstack^\lambda_X \ar{r}{\myuline{\MF}^\infty} & \shcatsqcohtwoperfunctornabla_X\\
    \stackQuadnabla_X \ar[phantom,start anchor=center, end anchor=center]{u}{\circlearrowleft} \ar{r}{\myuline{\MF}^\infty} & \Azumayatwopertwotorsionnabla_X\rlap{.} \ar[phantom,start anchor=center, end anchor=center]{u}{\circlearrowleft}
  \end{tikzcd}
  \]
  \begin{proof}
    Composing  the functor	$\myuline{\MF}^{\infty, \acts, \twotorsion, \Liouv, \lambda}_{\simplex-\constant\times \derham}$ of \cref{finalstepfinally} with the functor $X_\et^\daff\to  \dInfSch_{X_\derham}$ sending $S\mapsto S_\derham$ and extracting maximal $\infty$-groupoids we obtain a new \ifunctor
    \[
      \Delta^1 \to \Fun(\Delta^\op, \Mod(\Fun((X_\et^\daff)^\op, \inftygpd)))
    \]
    Evaluating along the zero simplex $[0]\in \Delta^\op$ we obtain the \ifunctor  
    \[
      \Delta^1 \to \Mod(\Fun((X_\et^\daff)^\op, \inftygpd))
    \]
    encoding a natural transformation
    \[
      \begin{tikzcd}[row sep=tiny, column sep=scriptsize]
        \stackfactor^{\Liouv,\lambda}_{p_X}((-)_\derham)  \ar{r}& \shcatsqcohtwoperfunctornabla_X
      \\
        \B\orthogonal^{\oplus}_{X_\derham}((-)_\derham) \ar[phantom,start anchor=center, end anchor=center]{u}{\circlearrowleft} \ar{r} & \Azumayatwopertwotorsionnabla_X \ar[phantom,start anchor=center, end anchor=center]{u}{\circlearrowleft}
      \end{tikzcd}
    \]
    which we can compose with \eqref{EquivariantDarbinclusion1} to obtain \eqref{finalmapfromdarbouxstacklevelzero}. Notice that since all $\Darbstack_X$, $\stackQuadnabla_X$, $\shcatsqcohtwoperfunctornabla_X$ and $\Azumayatwopertwotorsionnabla_X$ are stacks (see \crossref{stackquaddef}, \crossref{definitiondarbouxstackasstacks} and \cref{notationstackofcategoriesonetalesite}) and since the inclusion of sheaves in presheaves preserves products, the diagram lands in the full subcategory
    \[
      \Mod(\Sh(X_\et^\daff, \inftygpd)^{\times}) \subseteq \Mod(\Fun((X_\et^\daff)^\op, \inftygpd)).\qedhere
    \]
  \end{proof}
\end{proposition}

\begin{definition}\label{definitionofcliffordstack}
  Fix a square root of $-1$ and let $X$ be a derived Deligne--Mumford stack. We define the \emph{Clifford stack} $\cliffordstacknabla_X$ as the stack in \igroupoids given by the stack-theoretic image of the morphism of prestacks on $X_\et^\daff$
  \[
    \myuline{\MF}^\infty \colon \stackQuadnabla_X \to \Azumayatwopertwotorsionnabla_X.
  \]
  We hence get a morphism in $\Mod(\Sh(X_\et^\daff, \inftygpd)^{\times})$ also denoted by $\myuline{\MF}^\infty$:
  \begin{equation}\label{finalmapfromdarbouxstacklevelzero}
    \begin{tikzcd}[row sep=small]
      \Darbstack^\lambda_X \ar{r}{\myuline{\MF}^\infty} & \shcatsqcohtwoperfunctornabla_X\\
      \stackQuadnabla_X \ar[phantom,start anchor=center, end anchor=center]{u}{\circlearrowleft} \ar{r}{\myuline{\MF}^\infty} & \cliffordstacknabla_X\rlap{.} \ar[phantom,start anchor=center, end anchor=center]{u}{\circlearrowleft}
    \end{tikzcd}
  \end{equation}
\end{definition}

\section{Gluing Matrix Factorizations}
\label{sectiongluing}

In this section we formulate and prove our main theorem (see \cref{theoremgluingMF} below).  We fix the data required in the statement of \crossref{theoremcontractibility}, namely, a $(-1)$-shifted symplectic derived Deligne-Mumford $X$ equipped with its canonical exact structure $\lambda$ (\cf \crossref{canonicalexactstructure}).

\subsection{Gluing an Azumaya-twisted \texorpdfstring{$\MF^\infty$}{MF∞} up to \texorpdfstring{$\affineline{}$}{A¹}-isotopies}
\label{sectionA1invarianceofMFrevisited}

The goal of this section is to formulate our main theorem. In \cref{sectionhomotopygroupscliff} and \cref{sectioncategoricalorientationdataandpin} below we discuss the three obstruction classes that allow us to \emph{untwist} the gluing. 

\begin{construction}\label{takingsimplicialcolimit}
	Let $\calC^{\times}$ be a cartesian symmetric monoidal \icategory admitting sifted colimits. Since sifted colimits are compatible with finite products, the colimit functor induces an \ifunctor $\colim_{\Delta^\op} : \Fun(\Delta^\op, \Mod(\calC))\to \Mod(\calC)$.  We  observe that the \ifunctor $\myuline{\MF}^{\infty, \acts, \twotorsion, \Liouv, \lambda}_{\simplex-\constant\times \derham}$ of \cref{finalstepfinally} takes values in the full subcategory of stacks
	\[
\Sh\left(\dInfSch_{X_\derham}, \inftycats\right)	\subseteq\Fun(\dInfSch_{X_\derham}^\op, \inftycats)
	\]
and apply this discussion by considering colimit in stacks 
    \[
      \begin{tikzcd}[column sep=0.3in]
        \Delta^1\ar{r} & \Fun(\Delta^\op, \Mod(\Sh\left(\dInfSch_{X_\derham}, \inftycats\right))) \ar{r}{\colim} & \Mod(\Sh\left(\dInfSch_{X_\derham}, \inftycats\right))
      \end{tikzcd}
	 \]
 encoding an equivariant map
 \[
 	\begin{tikzcd}[row sep=tiny, column sep=scriptsize]
 		\displaystyle \colim_{\Delta^\op}\,\stackfactor^{\Liouv,\lambda, \dInf, \constant\times \derham}_{p_X}  \ar{r}& \displaystyle \colim_{\Delta^\op}\,\shcatsqcohtwoper((-)_\derham\times \simplex^\bullet)^\constant
 		\\
 	\displaystyle \colim_{\Delta^\op}\,	\left(\B\orthogonal^{\oplus}_{X_\derham}\right)_{|_\dInf}(\simplex^\bullet\times (-)_\derham)^{\constant} \ar[phantom,start anchor=center, end anchor=center]{u}{\circlearrowleft} \ar{r} & \displaystyle \colim_{\Delta^\op}\,\Azumayatwopertwotorsion((-)_\derham\times \simplex^\bullet)^\constant\rlap{.} \ar[phantom,start anchor=center, end anchor=center]{u}{\circlearrowleft}
 	\end{tikzcd}
 \]
Following \crossref{mappingspacesofTareA1localizations-pi0} we know that each simplicial object is constant on $\pi_0$-presheaves so that this colimit does not alter the $\pi_0$'s.  
\end{construction}

\begin{definitions}{ Following \cref{takingsimplicialcolimit}, we define\footnote{See also \crossref{constructionofquotientA1stacksgeneral} and \crossref{Liouvilleuptoisotopy}}:}
	\item $\shcatsqcohtwoperfunctornablagrpisotopy_X$ as the restriction of $\colim_{\Delta^\op}\,\shcatsqcohtwoper((-)_\derham\times \simplex^\bullet)^\constant$ (as above) to the small affine étale site of $X$.
	\item $\cliffordstacknablaisotopy_X$ as the full substack of the restriction of $\colim_{\Delta^\op}\,\Azumayatwopertwotorsion((-)_\derham\times \simplex^\bullet)^\constant$ to the small étale site of $X$ with $\homotopysheaf_0$ given by $\homotopysheaf_0(\cliffordstacknabla_X)$.
\end{definitions}

\begin{corollary}
  \label{homotopyfactorization}
  By the universal property of the isotopic localization, we have a canonical commutative square in $\Mod(\Sh(X_\et^\daff, \inftygpd)^{\times})$ 
  \[
    \begin{tikzcd}
      \left(\stackQuadnabla_X\acts \Darbstack^\lambda_X\right) \ar{d}\ar{r}{\myuline{\MF}^\infty}& \left(\cliffordstacknabla_X\acts \shcatsqcohtwoperfunctornabla_X\right) \ar{d}\\
      \left(\QA\acts \DA\right)\ar[dashed]{r}[swap]{\myuline{\MF}^\infty_{\affineline{}}}& \left(\cliffordstacknablaisotopy_X\acts \shcatsqcohtwoperfunctornablagrpisotopy_X\right)\rlap{.}
    \end{tikzcd}
  \]
  In particular, it descends  the (stacky) étale quotients (\cf \crossref{theoremcontractibility}):
  \[
    \constantsheaf \ast_X \simeq \quotA \to \quot{\shcatsqcohtwoperfunctornablagrpisotopy_X}{\cliffordstacknablaisotopy_X}.
  \]
\end{corollary}

The above statement can be informally rephrased as follows: Every $(-1)$-symplectic Deligne--Mumford stack $X$ has a \emph{twisted} section of $\shcatsqcohtwoperfunctornablagrpisotopy_X$ constructed from matrix factorizations.
In the next section we exhibit the necessary data required to resolve this twist.

\subsection{Homotopy groups of the Clifford stack}
\label{sectionhomotopygroupscliff}
The next theorem explains why only three obstruction classes appear in \cref{theoremB}.
\begin{theorem}\label{homotopygroupsofcliff}
  Let $X$ derived Deligne--Mumford stack.
  We have
  \[
    \homotopysheaf_n\left(\cliffordstacknabla_X\right) \simeq
    \left\{
    \begin{array}{lll}
      \quot{\integers}{2} & \text{if } n = 0 & (\MF^\infty(X_\deRham, 0) \text{ and } \MF^\infty(\affineline{X_\deRham}, x^2)) \\
      \quot{\integers}{2} & \text{if } n = 1 & (\Id \text{ and } [1]) \\
      \mu_2 \simeq \quot{\integers}{2} & \text{if } n = 2 & (\pm 1 \colon \Id \to \Id) \\
      0 & \text{for } n \geq 3.
    \end{array}
    \right.
  \]
  Note that $\cliffordstacknabla_X$ is a group-stack, so this is independent of the choice of a base-point. In particular, by \cite[6.5.2.9]{lurie-htt}, $\cliffordstacknabla_X$ satisfies étale hyperdescent.

  Moreover, the Postnikov tower of $\cliffordstacknabla_X$ splits, so that $\cliffordstacknabla_X \simeq \quot{\integers}{2} \times \B \quot{\integers}{2} \times \EM(\quot{\integers}{2}, 2)$.
\end{theorem}

\begin{lemma}\label{lemmahomotopysheavesofcliff}
  Denote by $A$ the restriction $A \coloneqq \Azumayatwopertwotorsion_{|\Aff_\basefield} \colon \Aff_\basefield^\op \to \inftygpd$.
  The group stack $\Omega_1 A$ is equivalent to $\quot{\integers}{2} \times \B \quot{\integers}{2}$.
  \begin{proof}
    \newcommand{\picardstack}{\myuline{\Picard}\left(\kbetaLaurent\right)}
    \newcommand{\thecliff}{A}
    We first compute its homotopy sheaves $\homotopysheaf_n$. They are computed as the $2$\nobreakdashes-torsion in
    \[
      \Omega_{\unit}\shcatsqcohtwoper \simeq \picardstack
    \]
    where the RHS is the Picard $\infty$\nobreakdashes-stack $S \mapsto \Picard\left( \DQCoh(S) \ofbetaLaurent \right)$ (see \cite{AkhilMathew955} for its descent properties).
    In particular, the fiber sequence of \cref{2torsionstackdefinition} gives a long exact sequence
    \newcommand{\snakearrow}{\ar[dll, rounded corners, to path={ -- ([xshift=2ex]\tikztostart.east)
            |- ($(\tikztostart)!0.5!(\tikztotarget)$)
            -| ([xshift=-2ex]\tikztotarget.west)
            -- (\tikztotarget)}]}
    \[
      \begin{tikzcd}[row sep=small]
        \cdots \ar{r} & \homotopysheaf_{3}\left(\thecliff \right) \ar{r} & \homotopysheaf_{2}\left(\picardstack\right) \ar{r}{(-)^{\otimes 2}} & \homotopysheaf_{2}\left(\picardstack\right)
        \snakearrow
        \\
        & \homotopysheaf_{2}\left(\thecliff \right) \ar{r} & \homotopysheaf_{1}\left(\picardstack\right) \ar{r}{(-)^{\otimes 2}} & \homotopysheaf_{1}\left(\picardstack\right)
        \snakearrow
        \\
        & \homotopysheaf_{1}\left(\thecliff \right) \ar{r} & \homotopysheaf_{0}\left(\picardstack\right) \ar{r}{(-)^{\otimes 2}} & \homotopysheaf_{0}\left(\picardstack\right). & \phantom{\cdots}
      \end{tikzcd}
    \]
    By definition, for a symmetric monoidal \icategory $\calC$, we have
    \[
      \pi_n(\Picard(\calC)) \simeq
      \begin{cases}
        \pi_0(\End(\unit))^{\times} & \text{if } n = 1 \\
        \pi_{n-1}(\End(\unit)) & \text{if } n \geq 2.
      \end{cases}
    \]
    In particular, $(-)^{\otimes 2}$ induces an automorphism (multiplication by $2$) on $\pi_n(\Picard(\calC))$ for $n \geq 2$, as the group structure is linear (and we are working in characteristic $0 \neq 2$).
    We get $\homotopysheaf_{n}\left(\thecliff\right) \simeq 0$ for $n \geq 3$ and the above long exact sequence becomes:
    \[
      \begin{tikzcd}[row sep=small]
        0 \ar{r} & \homotopysheaf_{2}\left(\thecliff\right) \ar{r} & \Gm{} \ar{r}{(-)^{2}} & \Gm{}
        \snakearrow
        \\
        & \homotopysheaf_{1}\left(\thecliff\right) \ar{r} & \homotopysheaf_{0}\left(\picardstack\right) \ar{r}{(-)^{\otimes 2}} & \homotopysheaf_{0}\left(\picardstack\right).
      \end{tikzcd}
    \]
    In particular, we find $\homotopysheaf_{2}\left(\thecliff\right) \simeq \mu_{2}$.
    Moreover, the square morphism $(-)^2 \colon \Gm{} \to \Gm{}$ is a sheaf epimorphism so that we have an exact sequence
    \[
      \begin{tikzcd}[row sep=small]
        0 \ar{r} & \homotopysheaf_{1}\left(\thecliff\right) \ar{r} & \homotopysheaf_{0}\left(\picardstack\right) \ar{r}{(-)^{\otimes 2}} & \homotopysheaf_{0}\left(\picardstack\right).
      \end{tikzcd}
    \]
    It remains to compute $\homotopysheaf_{1}\left(\thecliff\right)$.
    For any $S \in \Aff_\basefield$, we have an exact sequence as in \cite[Rmk.\,2.4.8]{AkhilMathew955}
    \[
      0 \to \pi_0\Picard(S) \to \pi_0\Picard(S\ofbetaLaurent) \to \quot{\integers}{2} \to 0.
    \]
    It induces an exact sequence of sheaves $0 \simeq \homotopysheaf_0\myuline{\Picard} \to \homotopysheaf_0\left(\picardstack\right) \to \quot{\integers}{2} \to 0$ and thus $\homotopysheaf_0\left(\picardstack\right) \simeq \quot{\integers}{2}$.
    We deduce $\homotopysheaf_1(A) \simeq \quot{\integers}{2}$.
    Note that by construction, the generator of $\quot{\integers}{2}$ corresponds to the shift $[1]$ under this isomorphism.

    It remains to show the Postnikov tower splits: the shift functor $[1]$ defines a section of $\Omega_1\thecliff \to \homotopysheaf_0(\thecliff) \simeq \quot{\integers}{2}$ as abelian group stacks, thus proving the claim.
  \end{proof}
\end{lemma}

\begin{proof}[Proof of \cref{homotopygroupsofcliff}]
  \newcommand{\picardstack}{\myuline{\Picard}\left(X_\deRham\ofbetaLaurent\right)}
  By definition, the quadratic bundle $(\affineline{X_\deRham}, x^2)$ defines a sheaf isomorphism $\quot{\integers}{2} \to \homotopysheaf_0\left(\cliffordstacknabla_X\right)$.
  The higher homotopy sheaves are computed using \cref{lemmahomotopysheavesofcliff}:
  denote by $A^\nabla \colon \Aff_\basefield^\op \to \inftygpd$ the stack $S \mapsto \Azumayatwopertwotorsion(S_\deRham)$.
  By \cref{lemmahomotopysheavesofcliff}, we have $\Omega^1 \cliffordstacknabla_X \simeq \Omega_1 A^\nabla \simeq \Map\left((-)_\deRham, \quot{\integers}{2}\right) \times \Map\left((-)_\deRham, \B \quot{\integers}{2}\right)$.
  The natural transformation $\Id \to (-)_\deRham$ induces morphisms $\xi \colon \Map\left((-)_\deRham, \quot{\integers}{2}\right) \to \quot{\integers}{2}$ and $\zeta \colon \Map\left((-)_\deRham, \B\quot{\integers}{2}\right) \to \B \quot{\integers}{2}$.
  The morphism $\xi$ is obviously an isomorphism and we will prove that so is $\zeta$.
  In fact, because $\xi$ is an isomorphism, it suffices to prove that for any $S$, any $\quot{\integers}{2}$-bundle on $S_\deRham$ is étale locally trivial.
  Let us fix such a bundle $E$.
  
  Up to localizing, we can assume $S$ is a closed subscheme of a smooth formal thickening $\calY$.
  As $S_\deRham$ and $\calY_\deRham$ agree, we have $\Picard(S_\deRham) \simeq \Picard(\calY_\deRham)$.
  In particular, the bundle $E$ amounts to a $2$-torsion line bundle with connection on $\calY$.
  Up to further localization, we can assume the underlying bundle (without its connection) is trivial.
  Since the space of connections on the trivial line bundle is linear, its $2$-torsion vanishes and hence $E$ is itself trivial.
\end{proof}

\begin{lemma}\label{lemmaBZ2A1local} 
  The classifying stack $\B \quot{\integers}{2}$ is $\affineline{}$-local.
  \begin{proof}
    Let $A$ be a local ring and $S \coloneqq \Spec A$. Let $P$ be a principal $\quot{\integers}{2}$-bundle on $S \times \affineline{}$.
    Denote by $U_1 \amalg \cdots \amalg U_n \to \affineline{}$ an open cover such that $P$ is trivial on each $S \times U_i$.
    As each $U_{ij} \coloneqq U_i \cap U_j$ is connected, the transition functions $S \times U_{ij} \to \quot{\integers}{2}$ are necessarily constant.
    Their cocyclicity implies they are all equal and the bundle $P$ is in fact trivial.
    In particular, for any $S \in \Aff_\basefield$, and principal $\quot{\integers}{2}$-bundle on $S \times \affineline{}$ is trivial locally on $S$.

    This holds for derived schemes as well: for any $S \in \dAff_\basefield$, we have $\Map(S \times \affineline{}, \B \quot{\integers}{2}) \simeq \Map(S_\red \times \affineline{}, \B \quot{\integers}{2})$ (since $\B \quot{\integers}{2}$ is étale, it is equivalent to its de Rham stack).
    As a consequence, the trivial bundle yields a stack epimorphism
    \[
      \alpha \colon * \to \stackMap\left(\affineline{}, \B\quot{\integers}{2}\right).
    \]
    The result follows as the nerve of $\alpha$ is equivalent to the nerve of $* \to \B \quot{\integers}{2}$.
  \end{proof}
\end{lemma}

\begin{corollary}\label{corollarycliffordisisotopyinvariant}
  Let $X$ be a derived Deligne--Mumford stack. Then the canonical map  
  \[\cliffordstacknabla_X\to \cliffordstacknablaisotopy_X\]
  is an equivalence of monoid stacks on the small étale site of $X$.
  \begin{proof}
    Follows from \cref{lemmahomotopysheavesofcliff,lemmaBZ2A1local}, using \crossref{lemmamappingspacesA1isotopiclocalizationsareA1localizations}.
  \end{proof}
\end{corollary}

\subsection{Categorical orientation data and Pin structures}
\label{sectioncategoricalorientationdataandpin}
In this section we describe the three classes appearing in \cref{theoremgluingMF}. Recall the universal orientation obstruction from \crossref{defuniversalorientation}, defined for any $(-1)$\nobreakdashes-shifted exact symplectic Deligne--Mumford stack $X$:
\[
  \universalobstruction_X \colon \constantsheaf \ast_X \lra \B\QA.
\]

\begin{definition}\label{defcategoricalorientation}
  The \emph{categorical orientation obstruction} is the composition
  \[
    \categoricalobstruction_X \colon \constantsheaf \ast_X \lra \B\QA \lra \B \cliffordstacknablaisotopy_{X} \simeq \B \cliffordstacknabla_X.
  \]
  The Postnikov splitting of \cref{homotopygroupsofcliff} implies the categorical orientation obstruction amounts to three orientation classes
  \begin{align*}
    \beta_1 & \in \cohomology^1_\et\left(X, \quot{\integers}{2}\right), \\
    \beta_2 & \in \cohomology^2_\et\left(X, \quot{\integers}{2}\right) \\
    \text{and }\beta_3 & \in \cohomology^3_\et\left(X, \quot{\integers}{2}\right).
  \end{align*}
\end{definition}

\begin{lemmalist}{\label{describebetas}Consider the classes $\beta_1$, $\beta_2$ and $\beta_3$ of \cref{defcategoricalorientation}.}
  \item The class $\beta_1$ always vanishes. It thus classifies the (trivial) bundle of (locally constant) $\quot{\integers}{2}$\nobreakdashes-valued functions on $X$.
  \item The class $\beta_2$ classifies the gerbe of square roots of the determinantal bundle $\det \tangent_X$ (thus coincides with the obstruction in the case of vanishing cycles).
  \item The class $\beta_3$ arises from Pin (or Spin) structures and is a delooping of the second Stiefel--Whitney class (associated to quadratic bundles).
  We may describe  the $2$-gerbe it classifies as the 2-gerbe of \emph{spinorial theories}, denoted by $\spintheories_X$ and defined informally as:
        \[
          \spintheories_X(S) = \left\{ \sigma \colon \stackLagDist_S \to \EM\left(\quot{\integers}{2}, 2\right) ~\Big|~ \sigma(L \oplus Q) = \sigma(L) + \spinornorm(Q) \right\}
        \]
        Here $\stackLagDist_S$ denotes the stack of Lagrangian distributions on $S$ while $Q$ is a quadratic bundle and $\spinornorm(Q)$ is the spinor norm of $Q$, \ie the gerbe of Pin structures on $Q$, see \cite{bassCliffordAlgebrasSpinor1974}.
\end{lemmalist}
\begin{proof}
  The monoidal morphism $\MF^\infty \colon \stackQuadnabla_X \to \cliffordstacknabla_X$ together with \cref{homotopygroupsofcliff} associates classes $w_i(Q) \in \cohomology^i_\et(X, \quot{\integers}{2})$, $i \in \{0,1,2\}$ to any section $Q$ of $\stackQuadnabla_X$.
  By construction, we see that $w_0(Q)$ corresponds to the rank of $Q$ (mod $2$), $w_1(Q)$ classifies the determinant of $Q$ and $w_2(Q)$ classifies the gerbe of complex Pin structures\footnote{or Spin structures if $w_1(Q)$ is assumed to vanish.} on $Q$ (so the classes $w_1$ and $w_2$ are Stiefel--Whitney classes associated to $Q$).
  The description of $\beta_3$ now follows by unravelling the definition.

  Denote by $\calG_X$ the $\mu_2$\nobreakdashes-gerbe over $X$ of square roots of the determinantal bundle $\det(\cotangent_X)$.
  For $S \in X^\et$ and $i \colon S \hookrightarrow (\formalU, f)$ in $\Darbstack_X(S)$, the induced fiber sequence
  \[
    \begin{tikzcd}
      i^*\cotangent_\formalU \ar{r} & \cotangent_X \simeq \tangent_X[1] \to i^* \tangent_\formalU[1]
    \end{tikzcd}
  \]
  yields an isomorphism $\det(i^* \cotangent_\formalU)^2 \simeq \det(\cotangent_X)$, and thus a section $\rmD(\formalU, f) \in \calG_X(S)$.
  This assembles into an equivariant morphism of stacks over $X^\et$:
  \[
    \begin{tikzcd}[row sep=small]
      \Darbstack_X^\lambda \ar{r}{\rmD} & \calG_X \\
      \stackQuadnabla_X \ar{r}{\det} \ar[phantom, start anchor=center, end anchor=center]{u}{\circlearrowleft} & \B\mu_{2,X} \ar[phantom, start anchor=center, end anchor=center]{u}{\circlearrowleft}[right]{\rlap{ gerbe band action}} \rlap{.}
    \end{tikzcd}
  \]
  In particular, the morphism (which corresponds to the class $\beta_2$)
  \[
    \begin{tikzcd}
      \displaystyle \constantsheaf{\ast}_X \simeq \quotA \ar{r} & \B\QA \ar{r}{\B\det} & \EM(\mu_2, 2)
    \end{tikzcd}
  \]
  is the morphism classifying the $\mu_2$-gerbe $\calG_X$.
  Similarly, the dimension defines an equivariant morphism $(\formalU, f) \mapsto \dim \formalU \text{ mod } 2$
  \[
    \begin{tikzcd}[row sep=small]
      \Darbstack_X^\lambda \ar{r}{\dim} & \quot{\integers}{2} \\
      \stackQuadnabla_X \ar{r}{\rank} \ar[phantom, start anchor=center, end anchor=center]{u}{\circlearrowleft} & \quot{\integers}{2} \ar[phantom, start anchor=center, end anchor=center]{u}{\circlearrowleft}
    \end{tikzcd}
  \]
  showing that $\beta_1$ classifies the (trivial) bundle of $\quot{\integers}{2}$-values functions.
\end{proof}

\begin{definition}
  Let $X$ be a $(-1)$-symplectic derived Deligne--Mumford stack.
  Categorical orientation data for $X$ is a trivialization of $\categoricalobstruction_X$.
  In view of \cref{describebetas}, orientation data $\tau_X$ consists of a triple $(r_X, \sfK_X^{\quot{1}{2}}, \sigma_X)$:
  \begin{itemize}
   \item A (locally constant) function $r_X \colon X \to \quot{\integers}{2}$,
   \item A square root $\sfK_X^{\quot{1}{2}}$ of the canonical bundle $\sfK_X \coloneqq \det(\cotangent_X)$ and
   \item A spinorial theory $\sigma_X \in \spintheories_X(X)$.
  \end{itemize}
\end{definition}

\begin{observation}
  The critical locus of a function $\dCrit(f)$ (for $f \colon \formalU \to \fAffLine$) comes naturally equipped with canonical orientation data denoted by $\tau_{\dCrit(f)}^\mathrm{can}$ (see 
  \crossref{remarkorientationondCrit}).
  Its first two components can be easily described as $r_{\dCrit(f)}^{\mathrm{can}} = \dim \formalU$ and $\sfK_{\dCrit(f)}^{\quot{1}{2},\mathrm{can}} = \det(\cotangent_\formalU)_{|\dCrit(f)}$.

  Describing the third component, the spinorial theory $\sigma_{\dCrit(f)}^\mathrm{can}$, is a little more involved. The  $(-1)$-shifted symplectic form on $\tangent_X$ can be seen as a $(+1)$-shifted quadratic form, with a class in the Grothendieck-Witt group $\mathsf{GW}^{[1]}$ (see \eg \cite{Schlichting2017}). Given a Lagrangian distribution $L \to \tangent_X$, the difference between (the classes of) $L$ and $(\tangent_\formalU)_{|_S}$ in a suitable $\rmK$-theory group yields a class $c$ in the Grothendieck--Witt group $\mathsf{GW}^{[0]}$ of quadratic bundles on $X$, using the algebraic Bott sequence of \eg \cite[Thm.\,6.1]{Schlichting2017} for $n=0$.
  The gerbe $\sigma_{\dCrit(f)}^\mathrm{can}(L)$ is then the spinor norm of $c$.
\end{observation}

We can finally formally state and prove our main theorem:

\begin{theorem}\label{theoremgluingMF}
  Let $X$ be a $(-1)$\nobreakdashes-shifted symplectic derived Deligne--Mumford stack.
  Assume we have orientation data $\tau_X$ consisting of
  \begin{itemize}
    \item A continuous function $X \to \quot{\integers}{2}$ controlling the parity of the dimension of the local Darboux models -- corresponding to a trivialization of the zero class $\beta_1 = 0 \in \cohomology^1(X, \quot{\integers}{2})$;
    \item A square root of the canonical bundle $\canonicalbundle_X$ -- corresponding to a trivialization of a class $\beta_2 \in \cohomology^2(X, \quot{\integers}{2})$;
    \item A spinorial theory\footnote{See \cref{defcategoricalorientation}} on $X$ (akin to a shifted Pin structure) -- corresponding to the trivialization of a class $\beta_3 \in \cohomology^3(X, \quot{\integers}{2})$.
  \end{itemize}
  There exists a global section $\sfM_X$ of $\shcatsqcohtwoperfunctornablagrpisotopy_X$ such that for any local Darboux chart $X \supset S \simeq \dCrit(f) \hookrightarrow \formalU \to_f \affineline{}$ such that $\tau_{X|_S} \simeq \tau_S^{\mathrm{can}}$, we have
  \[
    \sfM_X(S) \simeq \left[\MF^\infty(\formalU, f)\right],
  \]
  where $[\MF^\infty(\formalU, f)]$ denotes the class of $\MF^\infty(\formalU, f)$ in $\shcatsqcohtwoperfunctornablagrpisotopy_X(S)$.
\end{theorem}

\begin{remarks}{In the situation of the theorem}
  \item Fix a local Darboux chart
  \[
    X \supset S \simeq \dCrit(f) \hookrightarrow \formalU \to^f \affineline{}.
  \]
  The difference between the chosen orientation $\tau_{X|_S}$ and the canonical $\tau_S^{\mathrm{can}}$ yields a function $S \to \quot{\integers}{2}$, a $\quot{\integers}{2}$-bundle and a $\quot{\integers}{2}$-gerbe.
  By \cref{homotopygroupsofcliff} that such a triple characterizes a section $\calA$ of $\cliffordstacknabla_X(S)$.
  The class $\sfM_X(S)$ is then equivalent to $[\MF^\infty(\formalU, f) \otimes_{\IndCoh(S_\deRham)\ofbetaLaurent} \calA]$.
  \item \label{invariantnontrivial} The invariant $\sfM_X$ we constructed is non-trivial in general. For instance, applying pointwise the $\affineline{}$-invariant functor of periodic homology followed by the Euler characteristic, we can recover the Behrend function $\nu_X \colon X \to \integers$.
  Indeed, let $\sfF \subset \shcatsqcohtwoperfunctornablagrp_X$ denote the full substack whose sections are sheaves $\calC \colon S \mapsto \calC(S)$ such that for any $S \in X_\et^\aff$, the category $\calC(S)$ is compactly generated.
  Denote by $\sfF^{\affineline{}}$ the corresponding full substack of $\shcatsqcohtwoperfunctornablagrpisotopy_X$.
  By construction, Euler characteristic of periodic homology of compact objects induces a morphism $\chi_{\mathsf{HP}} \colon \sfF \to \Fun(X, \integers)$ to the sheaf of $\integers$-valued functions.
  It canonically factors through $\homotopysheaf_0 \sfF = \homotopysheaf_0 \sfF^{\affineline{}}$.

  \cref{compactgenerationMF} implies that $\sfM_X$ has values in $\sfF^{\affineline{}}$ and \cite[Theorem 1.1]{1212.2859} that composing with $\chi_{\mathsf{HP}}$ yields Behrend's (constructible) function $\nu_X \colon X \to \integers$.
  \item Without orientation data, the $(-1)$-shifted symplectic Deligne--Mumford stack $X$ still comes with a \emph{twisted} section of $\shcatsqcohtwoperfunctornablagrpisotopy_X$.
  Note also that the space of choices of orientation data forms a gerbe étale over $X$, over which such a gluing always canonically exists.
\end{remarks}

\begin{proof}[{Proof of \cref{theoremgluingMF}}]
  By \cref{homotopyfactorization}, we have a global section of $\constantsheaf \ast_X \to \quot{\shcatsqcohtwoperfunctornablagrpisotopy_X}{\cliffordstacknablaisotopy_X}$.
  \cref{corollarycliffordisisotopyinvariant} implies $\cliffordstacknablaisotopy_X \simeq \cliffordstacknabla_X$.
  Consider now the solid commutative diagram
  \[
    \begin{tikzcd}
      & \shcatsqcohtwoperfunctornablagrpisotopy_X \ar{d} \ar{r} \tikzcart & \constantsheaf \ast_X \ar{d}{\mathrm{trivial}}
    \\
      \constantsheaf \ast_X \ar{r} \ar[dashed]{ru}{\sfM_X} \ar["\categoricalobstruction_X", rounded corners, to path={|- ([yshift=-.7em]\tikztotarget.south)[near end,swap]\tikztonodes -- (\tikztotarget)}]{rr} & \quot{\shcatsqcohtwoperfunctornablagrpisotopy_X}{\cliffordstacknabla_X} \ar{r} & \B \cliffordstacknabla_X\rlap{.}
    \end{tikzcd}
  \]
  The orientation data $\tau_X \colon \categoricalobstruction_X \sim 0$ thus provides a lift $\sfM_X$ as announced.
\end{proof}

\begin{remark}
  We conclude with an hypothetical interpretation of the class $\beta_3$ assuming we have a trivialization of $\beta_2$ (\ie a square root of $\det(\cotangent_X)$).
  Indeed, in this situation, the first Chern class of $\tangent_X$ is even and it provides an integral lift of its second Chern character:
  \[
    \chern_2(\tangent_X) \in \cohomology^4(X, \integers).
  \]
  On the other hand, the symplectic structure $\tangent_X \simeq \cotangent_X[-1]$ yields:
  \[
    \chern_2(\tangent_X) = \chern_2(\cotangent_X[-1]) = -\chern_2(\cotangent_X) = -\chern_2(\tangent_X).
  \]
  In particular $2 \chern_2(\tangent_X)$ is nullhomotopic in $\mathsf{C}^\bullet(X, \integers)$.
  From the triangle
  \[
    \mathsf{C}^\bullet(X, \integers) \lra^2 \mathsf{C}^\bullet(X, \integers) \lra \mathsf{C}^\bullet(X, \quot{\integers}{2}),
  \]
  the datum of $\chern_2(\tangent_X)$ and the nullhomotopy $2 \chern_2(\tangent_X) \sim 0$ amounts to a class in $\gamma \in \cohomology^3(X, \quot{\integers}{2})$.
  We believe the two classes $\beta_3$ and $\gamma$ coincide.
  In particular the class $\beta_3$ would be the obstruction to the nullhomotopy $2 \chern_2(\tangent_X) \sim 0$ being even.
\end{remark}

\printbibliography
\makeatletter\@input{refs.tex}\makeatother
\end{document}